\documentclass[12pt, a4paper]{amsart}
\usepackage{setspace}
\usepackage[margin=2.6cm]{geometry}

\usepackage[centertags]{amsmath} 
\usepackage{amsthm}
\usepackage{amssymb}
\usepackage{mathtools} 
\usepackage{enumitem}

\usepackage{mathrsfs}

\usepackage[english]{babel}
\usepackage[utf8x]{inputenc}

\usepackage[all]{xy}
\usepackage[bookmarks,breaklinks]{hyperref}

\usepackage[centertags]{amsmath}
\usepackage{amsthm}
\theoremstyle{plain}

\newtheorem{thm}{Theorem}[subsection]

\newtheorem{cor}[thm]{Corollary}
\newtheorem{fact}[thm]{Fact} 

\newtheorem{lemma}[thm]{Lemma}
\newtheorem{prop}[thm]{Proposition}

\theoremstyle{definition}
\newtheorem{defn}[thm]{Definition}
\newtheorem{rem}[thm]{Remark}

\newcommand{\A}{\mathbb{A}}
\newcommand{\Q}{\mathbb{Q}}
\newcommand{\R}{\mathbb{R}}
\newcommand{\Z}{\mathbb{Z}}

\newcommand{\HH}{\mathbb{H}}
\newcommand{\Hh}{\prescript{'}{}{H}}
\newcommand{\Kk}{\prescript{'}{}{K}}
\newcommand{\Mm}{\prescript{'}{}{M}}
\newcommand{\Dd}{\prescript{'}{}{D}}
\newcommand{\C}{\operatorname{\mathbb{C}}}
\newcommand{\CP}{\operatorname{\mathbb{CP}}}
\renewcommand{\P}{\operatorname{\mathbb{P}}}
\newcommand{\diam}{\operatorname{diam}}

\newcommand{\ev}{\operatorname{ev}}
\newcommand{\coev}{\operatorname{coev}}
\newcommand{\bu}{\bullet}
\renewcommand{\Im}{\operatorname{Im}}
\newcommand{\St}{\operatorname{St}}
\newcommand{\Stbar}{\operatorname{\overline{St}}}
\newcommand{\vol}{\operatorname{vol}}

\newcommand{\what}{\widehat}
\newcommand{\id}{\operatorname{id}}

\newcommand{\Ew}{\prescript{}{W}{E}}

\newcommand{\EG}[1]{\prescript{#1}{G}{E}}

\newcommand{\const}{\operatorname{const}}

\newcommand{\eps}{\varepsilon}
\newcommand{\an}{\operatorname{an}}
\newcommand{\red}{\operatorname{red}}
\newcommand{\Sk}{\operatorname{Sk}}

\newcommand{\Cells}{\operatorname{Cells}}

\renewcommand{\d}{\delta}

\newcommand{\Oo}{\operatorname{\mathcal O}}

\newcommand{\Aa}{\operatorname{\mathscr A}}
\newcommand{\Cc}{\operatorname{\mathscr C}}

\newcommand{\Ff}{\operatorname{\mathscr F}}

\newcommand{\Jj}{\operatorname{\mathscr J}}
\newcommand{\Pp}{\operatorname{\mathscr P}}
\newcommand{\Xx}{\operatorname{\mathscr{X}}}

\newcommand{\Cone}{\operatorname{Cone}}

\newcommand{\SL}{\operatorname{SL}}
\newcommand{\GL}{\operatorname{GL}}
\newcommand{\G}{\mathbb{G}}

\newcommand{\Gys}{\operatorname{Gys}}
\newcommand{\res}{\operatorname{res}}
\newcommand{\gr}{\operatorname{gr}}

\newcommand{\sgn}{\operatorname{sgn}}

\newcommand{\set}[1]{\{\,#1\,\}}
\newcommand{\suchthat}[2]{\{\ #1\ \mid\ #2\ \}}

\newcommand{\Pic}{\operatorname{Pic}}
\newcommand{\Aff}{\operatorname{Aff}}

\renewcommand{\sp}{\operatorname{sp}}
\newcommand{\Trop}{\operatorname{Trop}}
\newcommand{\Log}{\operatorname{Log}}

\newcommand{\Ker}{\operatorname{Ker}}

\newcommand{\Hom}{\operatorname{Hom}}

\newcommand{\6}{\partial}

\newcommand{\la}{\langle}
\newcommand{\ra}{\rangle}

\newcommand{\ul}{\underline}
\newcommand{\wo}{\setminus}

\numberwithin{equation}{section}

\begin{document}

\bibliographystyle{alpha}

\begin{center}
  \large \textbf{Combinatorial part of the cohomology of the nearby
    fibre}\\[2ex]
\end{center}

\begin{center}
  {\normalsize Dmitry Sustretov\footnote{ has
      received funding from the European Union's Horizon 2020 research
      and innovation programme under the Marie Sklodowska-Curie grant
      agreement No.~843100 (NALIMDIF) } \\[2ex]
  }

{\small 
\hspace{0.15\linewidth}
\begin{minipage}[t]{0.85\linewidth}
  \begin{center} {\bf Abstract}
  \end{center}  
  \footnotesize Let $f: X \to S$ be a unipotent degeneration of
  projective complex manifolds over a disc such that the reduction of
  the central fibre $Y=f^{-1}(0)$ is simple normal crossings, and let
  $X_\infty$ be the canonical nearby fibre. Building on the work of
  Kontsevich, Tschinkel, Mikhalkin and Zharkov, I introduce a sheaf of
  graded algebras $\Lambda^\bu$ on the dual intersection complex of
  $Y$, denoted $\Delta_X$. I show that there exists a map
  $H^q(\Delta_X, \Lambda^p) \to \gr^W_{2p} H^{p+q}(X_\infty, \Q)$,
  where $W$ is the monodromy weight filtration, which is injective
  whenever there exists a class $\omega \in H^2(Y)$ which is
  combinatorial and Lefschetz, a certain technical condition.  When
  $f$ is a Type~III Kulikov degeneration of $K3$ surfaces, the sheaf
  $\Lambda^1$ recovers the affine structure with singularities of
  Engel and Friedman on $\Delta_X$. In this case, I show that
  existence of such class follows from the existence of a positive
  $d''$-closed $(1,1)$-superform or supercurrent in the sense of
  Lagerberg on $\Delta_X$. The latter is established in the case of
  simple affine structure singularities in \cite{hessian}, in fact,
  the cohomology of sheaves $\Lambda^p$ coincides with the full nearby
  fibre cohomolgy then.
\end{minipage}
}
\end{center}
\setcounter{tocdepth}{1}
\tableofcontents

\section{Introduction}

\subsection{Overview of the results}

Given a reduced proper complex algebraic variety $X$ such that its
irreducible components are smooth and intersect normally, the
combinatorics of intersections of the components determines part of
the singular cohomology of $X$. This object appears, for example, in
\cite{abw13} under the name ``combinatorial part of the cohomology''
and is identified with the top weight piece of the associated graded
of the cohomology of $X$, with respect to the weight filtration
defined by Deligne \cite{dhodge2}. In a similar vein, recently
\cite{cgp1,cgp2} studied the top weight cohomology of the moduli space
of curves using tools from tropical geometry, which help describe the
combinatorics of the intersections of the components of the boundary
divisor in the Deligne-Mumford compactification. Finally, if
$f: X \to S$ is a proper morphism from a complex manifold to a disc,
smooth away from the snc divisor $Y:=f^{-1}(0)$, then by work of
Steenbrink \cite{ste76, ste95} weight 0 part of the limiting mixed
Hodge structure on the cohomology of $X_t$ for $t$ sufficiently small
is isomorphic to the singular cohomology of the \emph{dual
  intersection complex} of $Y$, a CW-complex encoding the
intersections of its irreducible components (a more functorial version
of this statement was obtained by Berkovich \cite{berk09}).

In this paper I extend the results of Steenbrink and Berkovich,
defining certain constructible sheaves $\Lambda^p$ of $\Q$-vector
spaces on the dual intersection complex for each integer $p \geq 0$,
such that cohomology of $\Lambda^p$ can be related to the cohomology
of $X_t$, recovering parts of cohomology in even weights of the
limiting mixed Hodge structure on $X_t$. These results are motivated
by the non-archimedean approach of Kontsevich and Sobelman to the SYZ
conjecture \cite{ks06}, where singular affine structures on the dual
intersection complexes of $Y$ for certain \emph{minimal} degenerations
play an important role. In fact, in the case of minimal Type~III
degenerations of K3 surfaces constructed by Kulikov, Persson and
Pinkham \cite{kul77,pp}, such singular affine structures have been
studied in \cite{ghk15, eng18, aet19}, and it turns out that the sheaf
$\Lambda^1$ defines the same singular affine structures for this class
of degenerations. For more detailed discussion, see
Section~\ref{sec:discuss}.

In the following denote $f:X \to S$ a proper morphism from a complex
manifold $X$ to a disc $S$, such that the reduction of the central
fibre $Y=\sum_{i \in I} N_i Y_i$ is a simple normal crossings
divisor. We assume $\Q$ coefficients for cohomology everywhere.
Denote $X_\infty$ the canonical nearby fibre and let $N = \log T$,
where $T: H^\bu(X_\infty,\R) \to H^\bu(X_\infty)$ is the monodromy
endomorphism, which is assumed to be unipotent. Call a class
$\omega \in H^2(Y)$ \emph{cohomologically K\"ahler} if
$\omega|_{Y_i} \in H^2(Y_i)$ is K\"ahler for all $i \in I$.  Call the
intersections of the form $Y_{i_1} \cap \ldots \cap Y_{i_k}$ the
\emph{strata} of the divisor $Y$ and denote $Y^{(p)}$ the union of
strata of $Y$ of codimension $p$ in $X$. Denote $\Hh^{\bu}(Y^{(p)})$
the subspace spanned by cycle classes of strata. The dual intersection
complex of $Y$ is a CW complex whose $k$-cells correspond to
codimension $k+1$ strata and are glued into the boundary of those
$k+1$ cells that correspond to the substrata (see
Section~\ref{sec:dual} for precise definition). I denote the dual
intersection complex of $Y$ as $\Delta_X$, underlining the fact that I
am only interested in the situation when $Y$ is smoothable.

I define (Section~\ref{sec:sheaves}) certain sheaves of graded
$\Q$-algebras $\Lambda^\bu$ on $\Delta_X$, $\Lambda^0=\Q$. There
always exists a sheaf $A^1$ which fits into the exact sequence
$$
  0 \to \Lambda^0 \to A^1 \to \Lambda^1 \to 0.
$$
In general, the graded component $\Lambda^1$ does not determine
$\Lambda^p, p > 1,$ even locally. When $\Lambda^p$ is determined by
$\Lambda^1$ in a neighbourhood of a face $\sigma$
(Definition~\ref{regular}), let us say that $\Lambda^\bu$ is regular
at the face $\sigma$. If $\Lambda^\bu$ is regular at all faces
$\sigma$ of $\Delta_X$ then higher graded components of the sheaf
$A^\bu$ can be defined, and they fit into an exact sequence
\begin{equation}
  0 \to \Lambda^p \to A^{p+1} \to \Lambda^{p+1} \to 0.  \label{seq:lambda}
\end{equation}

The first main result of the paper says that under a standard
requirement the cohomology of sheaves $\Lambda^p$ maps to the
cohomology classes of $X_\infty$ of even weight.

\vspace{2ex}
\noindent \textbf{Theorem~A}.  Assume that there exists a
cohomologically Kähler class $\omega \in H^2(Y)$. Then for every
$p,q$ there exists a morphism
$$
H^q(\Delta_X, \Lambda^p) \to \gr^W_{2p} H^{p+q}(X_\infty).
$$
Moreover, for any class $v \in H^q(\Delta_X, \Lambda^p)$ the image of
$v$ in $H^{q+1}(\Delta_X, \Lambda^{p-1})$ under the coboundary
morphism associated to the short exact sequence (\ref{seq:lambda})
coincides with $Nv$.

\vspace{2ex}
If $M$ is a complex K\"ahler manifold of dimension $n$, a class
$\omega \in H^2(M)$ is called a \emph{Lefschetz class} if
\begin{itemize}
\item the Lefschetz operator $L_\omega(x) = x \cup
  \omega$ on the cohomology algebra induces isomorphisms
$$
L_\omega^{n-i}: H^i(M) \to H^{2n-i}(M) 
$$
(\emph{Lefschetz property});
\item the Hodge decomposition on $H^k(M)$ is orthogonal with respect
  to the form
  $$
  \psi(x,y) = i^k \int_M x \wedge y \wedge \omega^{n-k}
  $$
  and $ i^{p-q-k} (-1)^{\frac{(p+q)(p+q-1)}{2}} \psi$ is positive
  definite on $\Ker L_\omega^{n-p-q+1} \cap H^{p,q}(M)$.
  (\emph{Hodge-Riemann bilinear relations}).
\end{itemize}
A classic theorem of Hodge theory states that a class $\omega$ is a
Lefschetz class if it is a K\"ahler class (see, for example,
\cite[Section~6.3.2]{voisin}). If $f: M \to N$ is a semi-small map to
a complex variety and 
$\omega$ is the first Chern class of an ample line bundle, then any
positive multiple of $f^*\omega$ is Lefschetz \cite{dcm02}.

We call a class $\omega \in H^{2j}(Y)$ \emph{combinatorial} if its
restriction to any irreducible component $Y_i$ is a linear combination
of cycle classes of strata.

\vspace{2ex}
\noindent \textbf{Theorem~A'}.  If there exists a combinatorial
Lefschetz class $\omega \in H^2(Y)$ then the morphism constructed in
theorem~A is injective and for $p > q$
$$
N^{p-q}: H^q(\Delta_X, \Lambda^p) \to H^p(\Delta_X, \Lambda^q)
$$
is an isomorphism. Moreover,
$\dim H^{n-q}(\Delta_X, \Lambda^{n-p}) = \dim H^q(\Delta_X,
\Lambda^p)$.

\vspace{2ex} For proof of Theorems~A and A' see Theorems~\ref{thm-a}
and \ref{hodge}.

One can observe easily that central fibres of degenerations of
curves or of degenerations of Abelian varieties with toric reduction
admit combinatorial Lefschetz classes. In the latter case,
$\Hh^i(Y^{(j)}) = H^i(Y^{(j)})$ for all $i,j$.

The theorem~A' is intended to be applied in the case of maximally
unipotent degenerations, i.e. when $T$ has a Jordan block of size
$n+1$, though I expect that the requirement on $\omega$ can be
reformulated in lower unipotency rank case to still yield the same
result.

It seems plausible that the cup product on cohomology of $\Lambda^p$
on $\Delta_X$ is compatible with the cup product on the cohomology of
$X_\infty$, and in particular, the equality of dimensions in
Theorem~A' follows from a Lefschetz property for the class in
$H^1(\Delta_X, \Lambda^1)$ that corresponds to a combinatorial
Lefschetz class in $H^2(Y)$. The proof of this statement is at present
hindered by the fact that there is no known formula for the cup
product on the Steenbrink complex that computes the cohomology of the
nearby fibre and that is used crucially in the proof of Theorem~A.

The superforms is a generalization of differential forms on a real
vector space to the context of affine geometry, introduced by
Lagerberg \cite{lag11}, that parallels the notion of $(p,q)$-forms in
complex analysis. The definition of superforms has been extended to
the setting of tropical varieties by Chambert-Loir and Ducros, see
\cite{gubler16, cld}. The tropical cohomology groups of Mikhalkin and
Zharkov \cite{mz14, ikmz} can be computed with superforms using a
resolution similar to the Dolbeault complex \cite{jss19}.
The cohomology of the sheaves $\Lambda^p$ can be computed with the
help of a similar resolution of superforms.  If the sheaf
$\Lambda^\bu$ is regular at every face of $\Delta_X$ then furthermore
the morphism $N$ can be lifted to the level of superforms.

Recall that a unipotent snc degeneration of K3 surfaces $f: X \to S$
is called a \emph{Kulikov degeneration} if $K_X = 0$.  In this case
the components of the central fibre are rational surfaces obtained
from toric surfaces by a finite number of blow-ups (this number is
called \emph{charge}), which can be interpreted as operations that
introduce singularities to the standard affine structure on the fan of
the toric surface. The existence of a combinatorial Lefschetz class on
the central fibre $Y = f^{-1}(0)$ can be characterized in this case
with the help of positive superforms or supercurrents.

\vspace{2ex}
\noindent \textbf{Theorem~B}.  Let $f:X \to S$ be a Kulikov
degeneration of K3 surfaces of Type~III.  Then $\Lambda^1$ is the
push-forward of the sheaf of parallel 1-forms from the complement of
the finite set of the singularities, with respect the affine structure
defined in \cite{ghk15}, \cite{eng18}.

Moreover,
\begin{enumerate}
\item If all irreducible components of $\Delta_X$ have charge 1 or 0
  then
$$
H^1(\Delta_X, \Lambda^1) \hookrightarrow \gr^W_2 H^2(X_\infty)
$$
is surjective, and
$$
H^q(\Delta_X, \Lambda^p) \to \gr^W_{2p} H^{p+q}(X_\infty)
$$
is an isomorphism for all other values of $p,q$;
\item If $\Delta_X$ admits a $\Lambda^1$-torsor with a convex PL
  metric or, equivantly, a $d''$-closed positive $(1,1)$-superform or
  supercurrent then this morphism is injective.
\end{enumerate}

The statements of Theorem~B are proved as
Propositions~\ref{lambda1-from-affine}, \ref{form-plmetric},
\ref{kahler-convex}, and \ref{k3-coho}.

In the article \cite{hessian}, I show that if the singularities of the
affine structure on $\Delta_X$ have unipotent monodromy (which is the
case for the affine structure on $\Delta_X$ when all components of the
central fibre have charge at most 1) then there exists a positive
$(1,1)$-supercurrent on $\Delta_X$.  In particular, the main result of
\cite{hessian} implies that in the case of Type~III Kulikov
degenerations the map
$$
H^q(\Delta_X, \Lambda^p) \to \gr^W_{2p} H^{p+q}(X_\infty)
$$
in the Theorem~A is an isomorphism for all $p, q \geq 0$.

\subsection{State of the art}

Let us make a brief overview of known results about degenerations that
are similar to Theorems~A and A'.

If a family $X^* \to S^*$ of complex varieties over a punctured disc
factors into an embedding $X \hookrightarrow \P^m \times S^*$ and the
projection on $S^*$ then one can associate to it yet another CW
complex, its \emph{tropical limit} $\Trop(X)$, see
\cite[Section~3.1]{ikmz}, \cite{payne} via logarithmic and
non-archimedean tropicalization maps. Geometrically, this complex can
be interpreted as a dual intersection complex of the union of the
central fibre of the closure of $X$ in a certain family of toric
varieties over $S$ and the toric boundary. It is canonically embedded
into a certain compactification of $\R^m$ homeomorphic to an
$m$-simplex.

A tropical limit is in particular a tropical variety. A homology
theory for tropical varieties is developed in \cite{ikmz}. The
coefficients of this theory are given by certain cosheaves $\Ff_p$,
which are defined in terms of the local polyhedral structure of
$\Trop(X)$. The main result of \cite{ikmz} shows that under a certain
combinatorial condition of \emph{tropical smoothness}, spaces
$\Hom(H_q(\Trop(X), \Ff_p)$ are isomorphic to
$\gr^W_{2p}H^{p+q}(X_\infty)$ and $X_\infty$ has no cohomology in odd
weights. In particular, the monodromy is maximally unipotent. It is
also clear from the proof of the theorem that all cohomology classes
in $H^\bu(X_\infty)$ have type $(p,p)$, or that the limit mixed Hodge
structure of $X_\infty$ is of Hodge-Tate type. The latter property is
much stronger than maximal unipotency (though, for example, is known
to coincide with it for hyperk\"ahler manifolds, see
\cite{soldatenkov}), which serves as an illustration to how
restrictive the requirement of tropical smoothness is.

The tropical homology (and associated cohomology) theory enjoys many
properties that make it similar to the (co)homology of complex
manifolds, for example, tropical (co)homology groups enjoy Poincar\'e
duality \cite{jss19} on tropical manifolds and satisfy a version of
Lefschetz theorem on $(1,1)$-classes \cite{jrs17}. On the other hand,
there are strange pathological phenomena already for curves: the group
$H^1(\Trop(X), \Ff_1)$ can be infinite-dimensional, see
\cite[Section~4]{jell19}.

Let us also mention the toric degenerations of Gross and Siebert which
come by definition with singular affine structure on the dual
intersection complex. In this case cohomology of the certain sheaves
related to the mentioned singular affine structure coincides with the
hypercohomology of the complex of log-differential forms
\cite[Section~3]{gs2}. The relationship between these cohomology
groups and the nearby fibre cohomology have been studied in
\cite{ruddat2010log}.

The approach presented in this paper $\Lambda^p$ can be seen as
interpolating between the tropical cohomology and the Gross-Siebert
approach (see also \cite{yamamoto} for results that bridge the
two). On the one hand, one is not restricted to the degenerations that
satisfy tropical smoothness condition and the embedding of the
degeneration into $\P^m$ is not important. On the other hand, the
definition of sheaves $\Lambda^p$ works on any degeneration with
smooth total space and central fibre with snc support.

\subsection{Motivation and discussion}
\label{sec:discuss}

The study of degenerations of complex varieties has recently received
much interest in connection to mirror symmetry and in particular the
non-archimedean approach to the SYZ conjecture proposed by Kontsevich
and Soibelman. The main object of study of this research program is a
family of polarized Calabi-Yau manifolds $X^*$ of complex dimension
$n$ fibered over a punctured disc $S^*$ and having a maximally
unipotent monodromy. Such a family gives rise to a projective
Calabi-Yau variety $\Xx$ over the field of germs of meromorphic
fuctions $\C\{\{t\}\}$ and its Berkovich analytification $\Xx^{\an}$,
a locally ringed topological space, whose points are absolute values
on the residue fields of scheme-theoretic points of $\Xx$ which
restrict to the natural non-archimedean absolute value on
$\C\{\{t\}\}$. Kontsevich and Soibelman \cite{ks06} defined a certain
canonical subset $\Sk(\Xx^{\an}) \subset \Xx^{\an}$ and conjectured
that it is homeomorphic to a manifold of dimension $n$, and that
$\Xx^{\an}$ admits a retraction onto $\Sk(\Xx^{\an})$, such that the
general fibre, away from a codimension 2 subset of $\Sk(\Xx^{\an})$ is
isomorphic to a fibration in non-archimedean tori
$\set{ |x_1|=\ldots=|x_n|=1 } \subset (\G_m^n)^{\an}$. This retraction
can then be used to construct a \emph{mirror family} $\check{X}^*$ by
taking a Legendre dual of the singular integral affine structure on
$\Sk(\Xx^{\an})$ induced by the torus fibration, and reconstructing
the variety $\check{\Xx}$ so that
$\Sk(\check{\Xx}^{\an}) = \Sk(\Xx^{\an})$ and a retraction of
$\check{\Xx}$ onto $\Sk(\check{\Xx})$ induces the dual singular affine
structure.

This singular affine structure is important in the metric version of
the picture above. Let $L$ denote the realatively ample line bundle on
$X^*$, then by a theorem of Yau $X_t$ admits a Ricci-flat K\"ahler
metric with the fundamental form $\omega_t \in c_1(L_t)$. Kontsevich
and Soibelman conjecture that as $t$ tends to 0,
$$
(X_t, \dfrac{\omega_t}{\sqrt{\diam(X_t)}}) \to B
$$
in the sense of Gromov-Hausdorff, where $B$ is a metric space that
satisfies the following properties:
\begin{enumerate}
\item $B$ is homeomorphic to $\Delta_X$ (Conjecture~3 \cite{ks06});
\item there exists a dense open subset $B^{sm}$, such that $B \wo
  B^{sm}$ has Hausdorff codimension at least 2, and such that $B^{sm}$
  is an oriented Riemannian manifold of dimension $n$; 
\item $B^{sm}$ has an integral affine structure;
\item the metric $g$ on $B^{sm}$  satisfies a real
  Monge-Ampere equation, i.e. $g$ is given in some local affine
  coordinates $x_1, \ldots, x_n$ by
  $$
  g_{ij} = \dfrac{\partial^2 F}{\partial x_i \partial x_j} \qquad
  \det(g_{ij}) \equiv const,
  $$
  for some smooth function $F$.
\end{enumerate}
The space $B$ is also conjectured to be a sphere of real dimension $n$
if $X\to S$ is a family of Calabi-Yau in the strict sense
(i.e. $H^{i,0}(X_t) = 0$, $0 < i <n$), and to $\CP^{n/2}$ if fibres
$X_t$ are hyper-k\"ahler.

The set $\Sk(\Xx^{\an})$ is defined as the minimum locus of a certain
\emph{weight function} on $\Xx^{\an}$ (see \cite{mn15} for its
properties), but in fact this set can be defined without reference to
non-archimedean geometry at all. Whenever $X^*$ admits an extension to
a smooth manifold $X$ over the disc $S \supset S^*$ so that the
central fibre $Y$ has a strictly normal crossing support, there exists
a canonical embedding $\Delta_X \hookrightarrow \Xx^{\an}$ and
retraction $\Xx^{\an} \to \Delta_X$ (in fact, even strong homotopy
retraction) \cite{berkovich99, thuillier07}. The weight function can
be defined in terms of the order of vanishing of the extension of a
holomorphic volume form on the irreducible components $Y$.

In fact, one can weaken the requirement on the pair $(X,Y)$, namely,
Nicaise and Xu have realized that if $(X,Y_{\red})$ is a dlt pair, $X$
is $\Q$-factorial and $K_X+Y_{\red}$ is semi-ample (\emph{good minimal
  dlt degeneration}) then $\Sk(\Xx^{\an})=\Delta_X$ \cite{nicaise16}
and $\Xx^{\an}$ has homotopy type of $\Delta_X$, though there is no
canonical retraction any more. Crepant modifications of $X$ preserve
the homeomorphism type of $\Delta_X$ but alter its subdivision into
simplices. In \cite{nicxuyu19} the authors show that such retraction
exists in codimension 1 on $\Delta_X$ for good minimal dlt
degenerations.  See \cite{kx16, klsv18, mauri20} for further results
on dual intersection complexes of degenerations.

There has been considerable progress in proving properties $(ii)-(iv)$
of the Gromov-Hausdorff limit, see \cite{tosatti20} for a survey, but
property $(i)$ is only known for degenerations of Kummer surfaces
\cite{goto22}. In the case of degenerations of hyperk\"ahler
manifolds, the main technique, pioneered by Gross and Wilson is to
start with a fibration of a fixed hyperk\"ahler manifold into Abelian
varieties, then obtain the degeneration $X \to S$ via hyperk\"ahler
rotation of the complex structure. In \cite{oo18} the authors manage
to extend this technique to any degeneration of K3 surfaces using
approximation and careful study of the moduli space. See also
\cite{sus18} for a description of Gromov-Hausdorff limits of curves
with abelian differentials in terms of the minimum loci of the weight
function.

If we assume the property $(i)$ above, the conjectures $(ii)-(iv)$
imply that for each choice of polarisation there exists a set of
affine structures on $\Delta_X$ with respect to which the limit metric
satisfies the real Monge-Ampere equation. One can ask the following
question:

\noindent\textbf{Question A}: can these structures be described algebraically,
purely in terms of the polarisation $L$ ?

\vspace{2ex}
The definitions and results presented in this paper are intended as a
first step in the study of this question. 

The cohomology of $\Lambda^p \otimes \R$ can be computed using the
sheaves of superforms defined by Lagerberg \cite{lagerberg12}, see
\cite{jss19, dacl12, gubler16} for the construction in the tropical
geometry. A Lagerberg superform on a real vector space $V$ is a real
form on $V \oplus V$ which is invariant under translations along the
second summand. In particular, a $(p,q)$-Lagerberg form on a base of a
torus fibration can be lifted to a $(p+q)$-form on the total space of
the fibration. In \cite{bj17} Boucksom and Jonsson consider
\cite[Section~2]{bj17} for any snc degeneration $X \to S$ a map
$\Log_X$ to $\Delta_X$, defined on a neighbourhood of the special
fibre. If $X$ is maximally unipotent, the generic fibre of this map is
a complex torus of the same dimension as the skeleton. In fact, this
map is not canonical and depends on choices of coordinates near the
strata of the central fibre, but the error resulting from different
choices is $O(1/\log |t|)$. Similarly, the choice of local coordinates
on $\Delta_X$, i.e. sections of $\Lambda^1$, provides a choice of a
trivialization of the torus fibration (also up to a bounded error), so
Lagerberg forms can be locally lifted to forms on nearby fibres $X_t$,
gluing using a partition of unity should allow to lift the forms
globally, again, up to a controlled error.

I expect this lifting to be compatible with the morphism constructed
in Theorem~A in the limit as $t \to 0$. The data of a metric $g$ from
the statement of conjectures $(ii)-(iv)$ is equivalent to a
$(1,1)$-superform $\omega$ on $\Delta_X$ such that $\omega^n=\mu$,
where $\mu$ is the Euclidean measure on $\Delta_X$. It is natural to
expect that the lifting of $\omega$ should approach the Ricci-flat
K\"ahler metrics in the cohomology class $c_1(L_t)$ as $t \to 0$, and
so the cohomology class of $\omega$ should map under the morphism from
Theorem~A to the class of $c_1(L)$ in $\gr^W_2 H^2(X_\infty)$.

I suggest therefore that tackling the following question might be
useful in answering the Question~A:

\vspace{2ex}
\noindent\textbf{Question~B}: for a given maximally unipotent family
of polarizedCalabi-Yau manifolds over a punctured disc what is the set
of extensions $X \to S$ such that the first Chern class of the
polarisation lies in the image of the morphism constructed in
Theorem~A?

\subsection{Structure of the paper}

The background information on dual intersection complexes, Steenbrink
complex and Hodge-Lefschetz modules is recalled in
Section~2. Section~3 introduces the sheaves $\Lambda^\bu$ and
proves that their definion is independent of the subdivision of
$\Delta_X$ induced by the blow-ups of the strata. Certain complexes
quasi-isomorphic to the even rows of the monodromy weight spectral
sequence associated to Steenbrink complex are constructed in
Section~4.  They are then used to prove Theorems~A and A' in
Section~5. The singular affine structure on $\Delta_X$ for Type~III
Kulikov degenerations defined in \cite{ghk15, eng18} is related to
$A^1$ in Section~6, where Theorem~B is also proved.

\vspace{2ex}
\noindent\textbf{Acknowledgements}. I would like to thank Edouard Balzin,
Nero Budur, Julien Grivaux, Grisha Papayanov, and Robin
van der Veer for useful discussions.

\section{Background}

Let $f: X \to S$ be a proper morphism from a complex manifold $X$ to a
disc $S \subset \C$. Let $S^*$ be the punctured disc and let
$\tilde{S} \xrightarrow{exp} S^*$ be its universal cover. The
restriction of $f$ to $f^{-1}(S^*)$ is a locally trivial fibration by
Ehresmann's theorem and therefore the space
$X_\infty= X \times_S \tilde S$, the \emph{canonical nearby fibre}, is
homotopy equivalent to $X_t$ for any $t \neq 0$. The fundamental group
$\pi_1(S^*) \cong \Z$ acts on $X_\infty$ by deck transformations, and
we denote the action of the generator represented by a
counter-clockwise loop $H^\bu(X_\infty, \Q)$ as $T$. The operator $T$
is quasi-unipotent by a theorem of Borel \cite[Lemma~4.5,
Theorem~6.1]{schmid}. 

Assume that $f$ is a smooth morphism away from $Y = f^{-1}(0)$. Let $k$
be the projection $X_\infty \to X$ and let $i: Y \to X$ be the
embedding of the central fibre. Given a ring $R$ the nearby cycles
complex $\psi R := i^*Rk_* R$ has the property
$$
H^\bu(\psi R)_x \cong \varprojlim_{\eps \to 0} H^\bu(F_{x, \eps}, R)
$$
where $F_{x,\eps} = B_{x,\eps} \cap f^{-1}(0)$ and $B_{x,\eps}$ is an
$\eps$-ball centered at $x \in Y$ for a sufficiently small $\eps$
($F_{x,\eps}$ is called the \emph{Milnor fibre} of $f$ at
$x$). Moreover,
$$
\HH^\bu(Y, \psi R) \cong H^\bu(X_\infty, R)
$$

If $f$ is smooth away from the central fibre
$Y= \sum_{i \in I} Y_i E_i$ and the latter is a divisor with simple
normal crossings support then we will call such morphism an \emph{snc
  degeneration}.  If additionally the monodromy operator
$T$ is unipotent, we call
it a \emph{unipotent snc degeneration}.

\subsection{Dual intersection complexes}
\label{sec:dual}

Denote
$\ul\sigma$ the set of vertices of a simplex $\sigma$. Given a face
$\sigma$ and a vertex $i \in \ul\sigma$ there exists a unique face,
which we will denote $\6_i \sigma$, such that
$\6_i \sigma \subset \sigma, \ul{\6_i \sigma} \cup \{i\} = \ul\sigma$.

If $\sigma$ is a face of a simplicial complex, we will denote
$\St(\sigma)$ the \emph{open star} of $\sigma$: the union of interiors
of all simplices that contain it, or just $\sigma$ if $\sigma$ is a
vertex. Clearly, $\St(\sigma) = \cap_{i\in\ul\sigma} \St(i)$. We will
denote as $\Stbar(\sigma)$ the \emph{closed star} of a simplex
$\sigma$, that is, the union of simpleces that contain $\sigma$.  This
set has the natural structure of a simplical complex, and we will
denote $\Stbar^0(\sigma)$ its set of vertices.  The dual intersection
complex of an snc divisor generally has a structure of a
$\Delta$-complex, see \cite[Ch.~0]{hatcher} for definition.

Let $X$ be a smooth variety and let $D = \sum_{i=1}^m N_i D_i$ be a
divisor in $X$ with snc support; the connected components of finite
intersections $D_{i_1} \cap \ldots \cap D_{i_k}$ are called the
\emph{strata} of $D$. Let the set of $k$-dimensional cells of
$\Delta(D)$ be in bijective correspondence with the set of the strata
of $D$ of codimension $k+1$ in $X$. Suppose that the $k$-skeleton of
$\Delta(D)$ is already defined. Since the divisor $D$ has snc support
for any stratum $Z \subset D_{i_1} \cap \ldots \cap D_{i_k}$ and any
$l, 1 \leq l \leq k$ there exists a unique stratum
$$
Z_l \subset \bigcap_{j \neq l} D_{i_j}
$$
that contains $Z$.  The $k+1$-cell corresponding to $Z$ is glued to
the $k$-skeleton is such a way that the cells corresponding to $Z_l$
are in its boundary.  Thus in a dual intersection complex $\Delta_X$
the faces $\tau \in \St(\sigma)$ correspond to the strata $Y_\tau$
that contain $Y_\sigma$ and the faces $\tau \in \Stbar(\sigma)$
correspond to the strata $Y_\tau$ that intersect non-trivially with
$Y_\sigma$. It is noticed in \cite[Definition~8]{dfkx} that if $(X,Y)$
is a dlt pair, then the above procedure is well defined for the union
$Y^{=1}$ of the irruducible components of $Y$ with multiplicity 1.

We identify a cell $\sigma$ with the simplex
$$
\{ (x_i)_{i \in \ul\sigma} \in \R \ul\sigma \mid \sum_{i \in \ul\sigma}
N_i x_i = 1 \}
$$
We denote the affine space spanned by $\sigma$ as $\la \sigma \ra$ and
its tangent space as $T\la \sigma \ra$. In particular, each cell has
an integral affine structure associated to it, but even if $\Delta_X$
is a manifold, it is not always possible to glue these affine
structures to obtain an affine structure on $\Delta_X$.

Let $f: X \to S$ be an snc degeneration with the central fibre $Y$,
let $Z \subset Y$ be a closed connected submanifold. Let
$\pi: X' \to X$ be the the blow-up of $X$ in $Z$, let
$Y'= (f \circ \pi)^{-1}(0) = Y'_{\neq 0} \cup Y_0$ where $Y'_0$ is the
exceptional divisor and $Y'_{\neq 0}=\pi_*^{-1}(Y)$ is the strict
transform of $Y$. The dual intersection complex $\Delta(Y')$ then
admits a natural PL morphism to $\Delta_X$ and can be described as
follows.

If $Z$ is a stratum of $Y$ then the dual intersection complex of
$\Delta(Y')$ is obtained from $\Delta_X$ by the subdivision of the
face corresponding to $Z$ and the morphism $\Delta(Y') \to \Delta_X$
is a homeomorphism. In particular, after blowing sufficiently many
strata of $Y$ one obtains the dual intersection complex $\Delta(Y')$
of the special fibre $Y'$ of the blow-up can be ensured to be a
simplicial complex homeomorphic to $\Delta_X$.

If $Z$ is contained in a smallest stratum $Y_\sigma$ then $\Delta(Y')$
can be described as follows \cite[Paragraph~9]{dfkx}. The dual
intersection complex $\Delta(Y'_0 \cap Y'_{\neq 0})$ is the join
$\sigma \ast \Delta(Z \cap (\cup_{i \notin \ul\sigma} Y'_i))$. Since
each stratum of $\Delta(Y'_0 \cap Y'_{\neq 0})$ is sent by $\pi$ to a
stratum that intersects $Z$ non-trivially we have a map of $\Delta$
complexes
$$
\Delta(Y'_0 \cap Y'_{\neq 0}) \to \Stbar(\sigma)
$$
The dual intersection complex $\Delta(Y')$ is obtained by glueing the
cone over $\Delta(Y'_0 \cap Y'_{\neq 0})$ with $\Stbar(\sigma)$ via
this map. The map $\Delta(Y') \to \Delta_X$ is the map above on
$\Delta(Y'_0 \cap Y'_{\neq 0})$ and identity otherwise.

Following \cite[Sect.~1]{kashiwara}, we call the sheaves on a
$\Delta$-complex $\Sigma$ such that their restrictions to the
interiors of the faces of $\Sigma$ are constant sheaves
$\Sigma$-\emph{constant sheaves}. They admit the following
combinatorial description (p.~327 of \cite{kashiwara}):

\begin{lemma}
  \label{const-shv}
  The category of $\Sigma$-constant sheaves with values in an Abelian
  category $\Aa$ is equivalent to the category of functors from the
  poset of the cells ordered by inclusion $\Cells(\Sigma)$ to $\Aa$.
\end{lemma}

\begin{proof}
  Let $\Ff$ be a $\Delta$-constant sheaf. Consider the correspondence
  $\sigma \mapsto H^0(\St(\sigma), \Ff)$, together the restriction
  morphisms $H^0(\St(\sigma), \Ff) \to H^0(\St(\tau), \Ff)$ for all pairs of
  faces $\sigma \subset \tau$ it defines a functor $\Cells \to \Aa$.

  Conversely, if $F: \Cells(\Sigma) \to \Aa$ is a functor, define the
  presheaf $\Pp_F$ by putting $\lim_{D_U} F$ into correspondence to an
  open $U$, where $D_U$ is the diagram of cells and inclusions that
  intersect $U$ non-trivially. The sheafification of $\Pp_F(U)$ is a
  $\Delta$-constant sheaf.
\end{proof}

We further consider simplicial complexes with orientation induced by
an order on the verctices.

For any ordered set $I=\set{i_1, \ldots, i_k}, i_1 < \ldots < i_k$ and
any $i = i_l \in I$ we will denote
$$
\sgn(i, I) = \dfrac{i_1 \wedge \ldots \wedge i_k}{i_l \wedge i_1
  \wedge \ldots \wedge \what{i_l} \wedge \ldots i_k} =  (-1)^{l-1}.
$$

The open stars $\St(i)$ of vertices are contractible and cover
$\Delta_X$, so given a functor $F: \Cells(\Sigma) \to \Aa$,
 the Cech complex with respecto to the cover $\set{ \St(i) }_{i \in
  \Delta^0_X}$ computes the cohomology of $\Ff$:
$$
C^k(\Sigma, F) = \bigoplus_{|\ul\sigma|=k+1} F(\sigma),
\qquad
d( (a_\sigma) ) = \sum_\sigma \sum_{\sigma = \6_i\tau} \sgn(i, \ul\tau)
F(\sigma \subset \tau)(a_\sigma),
$$
since
$$
F(\sigma) = H^0(\St(\sigma), F).
$$

\subsection{Gysin and residue morphisms}

For any pair of faces $\sigma \supset \tau$ denote $\iota^\sigma_\tau$
the embeddings $Y_\sigma \to Y_{\tau}$.

We will use signed and unsigned Gysin and restriction maps which we
denote as follows
$$
\begin{array}{rclrl}
\res^\sigma_\tau & = &  (\iota^\tau_\sigma)^* & : & H^\bu(Y_\sigma)
                                                    \to H^\bu(Y_\tau) \\ 

  \Gys^\sigma_\tau &  = &   (\iota^\sigma_\tau)_! & : & H^\bu(Y_\sigma)
                                                    \to
                                                        H^{\bu+2m}(Y_\tau),
                                                        m=\dim \sigma -
                                                        \dim \tau \\
  \gamma^\sigma_i & = &  \sgn(i,\ul\sigma) \Gys^\sigma_{\d_i \sigma}
                                              & : &  H^\bu(Y_\sigma)
                                                    \to
                                                    H^{\bu+2}(Y_{\6_{i_l}
                                                    \sigma}) \\
  \gamma^{(k)} & = & \sum\limits_{|\ul\sigma|=k} \sum\limits_{i\in \ul\sigma} \gamma^\sigma_i&
                                                                                        : & H^\bu(Y^{(k)}) \to H^{\bu+2}(Y^{(k-1)}) \\
  \rho^{(k)}& = & \sum\limits_{|\ul\sigma|=k} \sum\limits_{\tau: \6_i
                  \tau=\sigma} \sgn(i, \ul\tau)\res^\sigma_\tau & : & H^\bu(Y^{(k)}) \to H^{\bu}(Y^{(k+1)})\\  
\end{array}
$$
One checks that
$$
\gamma^{(k-1)}\circ \gamma^{(k)}=0, \qquad \rho^{(k+1)} \circ \rho^{(k)} = 0.
$$
We will drop the superscript in $\res^\sigma_\tau$ when it is clear
from context.

\begin{lemma}
  \label{res-gys}
  Consider a commutative diagram of manifolds
  $$
  \xymatrix{
    Y \ar[r]^f & X \\
    W \ar[r]^F \ar[u]^G & Z \ar[u]^g \\
  }
  $$
  where all morphisms are embdennings and $Y$ and
  $Z$ meet transversely. Then
  $$
  g^*f_! \alpha = F_! G^* \alpha
  $$
  for any $\alpha \in H^*(Y)$.
\end{lemma}

Though this statement is folklore and well-known, we sketch a proof
for the benefit of the reader.

\begin{proof}
  Interpreting Gysin maps as compositions of Thom isomorphisms and
  restriction maps of local cohomology (see,
  e.g. \cite[VIII.2.4]{iversen}) we get the diagram
  $$
  \xymatrix{
    H^\bu(Y) \ar[r]_{\sim}^{\cup \tau_Y} \ar[d]^{G^*}& H^{\bu+l}_Y(X) \ar[r]^{r_X}
    \ar[d]^{g^*} & H^{\bu+l}(X) \ar[d]^{g^*} \\
    H^\bu(W) \ar[r]_{\sim}^{\cup \tau_W}  & H^{\bu+l}_W(Z) \ar[r]^{r_Z} & H^{\bu+l}(Z) \\
  }
  $$
  where $\tau_Y, \tau_W$ are Thom classes of normal bundles of $Y$ in
  $X$ and $W$ in $Z$, respectively, and the left square commutes by
  the naturality of Thom isomorphism. By interpreting local cohomology as
  the cohomology of the Thom space of the normal bundle and maps $r_X$
  and $r_Z$ as Pontryagin-Thom collapse, and using the fact that the
  square in the statement of the theorem is Cartesian, we observe that
  the right square in the diagram above commutes too.
\end{proof}

\begin{cor}
  \label{cor-res-gys}
  Let $f: X \to S$ be an snc degeneration, then for any stratum
  $Y_\sigma$, any $i,j \in \ul\sigma, i \neq j$ and any class
  $a \in H^r(Y_{\6_i \sigma})$,
  $$
  \res_{\6_j \sigma} \Gys^{\6_i \sigma}_j a = \Gys^{\6_j \6_i
    \sigma}_i \res_{\6_j \6_i \sigma} a. 
  $$
\end{cor}

\begin{lemma}
  \label{sum-res-gys}
  Let $f: X \to S$ be an snc degeneration, $Y=\sum_{i\in I} N_i Y_i$
  be the central fibre and assume that $\Delta_X$ is a simplicial
  comlpex. Let $\sigma \in \Delta_X$, then for any
  $a \in H^r(Y_\sigma)$
  $$
  \sum_{i \in \ul\sigma} N_i \res_\sigma \Gys^{\sigma}_{\6_i \sigma} a
  + \sum_{\tau: \6_j \tau=\sigma} N_j \Gys^{\tau}_\sigma
  \res_\tau a = 0
  $$
\end{lemma}

\begin{proof}
  By Lemma~\ref{res-gys}, Projection Formula \cite[7.3]{iversen} and
  \cite[7.5]{iversen} we have 
  \begin{multline*}
    0 = a \cup \sum_{j \in I} N_j (\iota^j)^* (\iota^j)_!  1_{Y_j} =
    \sum_{j \in \ul\sigma} N_j a \cup (\iota^\sigma)^* (\iota^j)_!
    1_{Y_\sigma} + \sum_{\tau: \6_j \tau=\sigma} N_j a  \cup (\iota^\tau_\sigma)_!1_{Y_\tau} = \\
    = \sum_{i \in \ul\sigma} N_j\res_\sigma \Gys^{\sigma}_{\6_i
      \sigma} a +
    \sum_{\tau: \6_j \tau = \sigma}N_j \Gys^{\tau}_\sigma \res_{Y_\tau} a, \\
  \end{multline*}
  (the first expression vanishes since the divisor
  $Y = \sum_{j \in I} N_j Y_j$ is principal) and therefore
  $$
  \sum_{j \in I} N_j (\iota^j)^* (\iota^j)_!  1_{Y_j} = 0.
  $$
  
\end{proof}

\subsection{Weight spectral sequence of Steenbrink}
\label{s:steenbrink}

Let $f: X \to S$ be a unipotent snc degeneration. We will recall below
the definition of the limit mixed Hodge structure on the cohomology of
$X_\infty$ and a spectral sequence due to Steenbrink that computes its
weight filtration.

Let $V$ be a vector space and $N$ be a nilpotent operator on $V$,
$N^{m+1}=0$, $N^m \neq 0$. Define inductively the following $2m$-step
filtration on $V$. First put
$W_0 := \Im N^m, W_{2m-1} := \Ker N^m, W_{2m} = V$. Supposing that
$W_0, \ldots, W_l$ and $W_{2m-l-1}, \ldots, W_{2m}$ are defined and
$m > l + 1$, put
$$
\begin{array}{lll}
W_{l+1} & := & \Im N^{m-l-1} \cap W_{2m-l-1}, \\
W_{2m-l-2} & := & \Ker (N^{m-l-1}: W_{2m-l-1}/W_l \to W_{2m-l-1}/W_l).
\end{array}
$$

\begin{defn}[Weight monodromy filtration]
  The filtration $W(M)_\bu$ associated to the nilpotent operator
  $N = \log T: H^\bu(X_\infty) \to H^\bu(X_\infty)$ in the way
  described above is called the \emph{weight monodromy filtration}.
\end{defn}

Steenbrink \cite{ste76, ste95} has defined certain complexes
$A^\bu_{\C}, A^\bu_{\Q}$ quasi-isomorphic to $\psi \C, \psi \Q$,
respectively, endowed with an increasing filtration $W_\bu$ such that
$A^\bu_{\Q}$ is filtered quasi-isomorphic to $A^\bu_{\C}$, and such
that
$$
W(M)_m H^\bu(X_\infty) = \Im (\HH^\bu(Y, W_m A^\bu_{\Q}) \to
\HH^\bu(Y,  A^\bu_{\Q}))
$$
under assumption that $Y$ is ``cohomologically K\"ahler'', i.e. that
there exists a class $\omega \in H^2(Y)$ that restricts to a
K\"ahler class on every irreducible component of $Y$.

Given a double complex $(C^{\bu,\bu}, d', d'')$, we will denote
$(sC^\bu,d)$ the total complex:
$$
sK^r = \bigoplus_{p+q=r} C^{p,q}, \qquad d = d' + d''.
$$

The spectral sequence associated to the filtration $W$ on $A^\bu_{\Q}$
degenerates on $E_2$.  Following \cite[Section 2]{gna90} we describe
its first sheet $\Ew_1^{\bu,\bu}$. For any $k \geq 1$ denote $Y^{(k)}$
the disjoint union of $k$-fold intersections of the irreducible
components $Y_i$ and denote
$$
\begin{array}{rcl}
  K^{i,j,k}(Y) & = & \left\{\begin{array}{ll}
                              H^{i+j-2k+n}(Y^{(2k-i+1)}, \Q),  & \textrm{ if } k
                              \geq \max\{0, i\},\\
                              0 & \textrm{ else, }
                              \end{array}\right.\\
K^{i,j}(Y) & = & \bigoplus_k K^{i,j,k}.\\
\end{array}
$$

Now denote
$$
\begin{array}{ll}
  d': K^{i,j,k} \to K^{i+1,j+1,k+1}, & d'=\rho^{(2k-i+1)},  \\
  d'': K^{i,j,k}\to K^{i+1,j+1,k}, & d''=-\gamma^{(2k-i+1)}, \\
  N: K^{i,j,k}\to K^{i+2,j,k}, & N=\id. \\
\end{array}
$$
Then
$$
\begin{array}{l}
  \Ew_1^{-r, q + r} = \HH^q(X_\infty, \gr^W_{q+r} A^\bu_{\Q})
  \Longrightarrow H^q(X_\infty, \Q)\\
  \Ew_1^{-r, q + r} = K^{-r, q-n}(Y)\\
\end{array}
$$
and $d_1 = d' + d''$. Since the terms $K^{i,j,k}$ only depend on the
the special fibre, we will further refer to them as $K^{i,j,k}(Y)$ and
to the terms of the first page of the spectral sequence as
$\Ew_1^{i,j}(Y)$.  Note that $\Ker N$ is a sub-spectral sequence of
$\Ew_1^{\bu,\bu}(Y)$ that computes (the associated graded of) the
cohomology of $Y$.

Let $Y=\sum_{i=1}^m Y_i$ be an snc divisor in a proper complex
manifold $X$. Recall that the weight spectral sequence of Deligne
\cite[3.2.4.1]{dhodge2}, \cite[4.7]{ps08}
$$
E_1^{i,j} = H^{2i+j}(Y^{(-i)}) 
\Longrightarrow H^{i+j}(X \wo Y),
$$
where $Y^{(0)}=X$ and the differential
$d^{i,j}_1: H^{2i+j}(Y^{(-i)}) \to H^{2i+j+2}(Y^{(-i-1)})$ is
$-\gamma^{(-i)}$, degenerates on the first page. We will denote
$M^\bu_p(X,Y)$ the complexes standing in the rows of the first page of
this spectral sequence after reindexing:
$$
M^q_p(X, Y) = E_1^{q-p,2p} = H^{2q}(Y^{(p-q)}).
$$

If $Y$ is the central fibre of an snc degeneration $f: X \to S$, then
for any $p$ denote $D_p^{\bu,\bu}(Y)$ the double subcomplex of
$K^{\bu,\bu,\bu}(Y)$ with the terms
$$
D^{q,r}_p(Y) = K^{-p+q+r, p+q+r-n, q}(Y) = H^{2r}(Y^{(p+q-r+1)}).
$$
In particular, 
$$
sD_p^q =  K^{-p+q,p+q-n} = \Ew_1^{-p+q,2p}(Y) = \bigoplus_{r=0}^{\min\{p,q\}} H^{2r}(Y^{(p+q-2r+1)}).
$$
The differentials $d', d''$ on $D_p^{\bu,\bu}$ have degrees
$(1,0), (0,1)$, respectively.  Denoting the stupid filtration on $M$
as $\sigma_{\leq m}$,
$$
\sigma_{\leq m}M_p^q(X,Y) = \left\{
    \begin{array}{ll}
      M_p^q(X,Y), & \textrm{ if } i \leq m,\\
      0, & \textrm{ otherwise, }
  \end{array}
  \right.
$$
one can think about the double complex $D^{\bu,\bu}_p$ as the complex
$$
0 \to \sigma_{\leq p}M^\bu_p(X,Y) \xrightarrow{\rho^{(p)}} \sigma_{\leq
  p}M^\bu_{p+1}(X,Y) \xrightarrow{\rho^{(p+1)}} \sigma_{\leq
  p}M^\bu_{p+2}(X,Y) \to \ldots
$$
in the category of complexes.

\subsection{Hodge-Lefschetz modules}

Recall that a rational Hodge structure of weight $n$ is a rational
vector space $V$ together with a decomposition
$V \otimes \C \cong \bigoplus_{p+q=n} V^{p,q}$ such that $\dim
V^{p,q}= \dim V^{q,p}$ for all $p,q$. A polarization on $V$ is a
symmetric form $\psi: V \otimes V \to \Q$ such that for its linear
continuation to $V \otimes \C$ the following holds:
\begin{enumerate}
\item $\psi(V^{p,q},V^{p',q'})=0$ for $p \neq p', q \neq q'$;
\item $i^{p-q} \psi(x, \bar x) > 0$ for all $x \in V^{p,q}$.
\end{enumerate}

Let $L=\bigoplus_{i,j \in \Z} L^{i,j}$ be a bigraded vector space such
that $L^{i,j}$ are endowed with real Hodge structures of weight
$p+i+j$ and let
$l_1: L^{i,j} \to L^{i+2,j}(1), l_2: L^{i,j} \to L^{i, j+2}(1)$ be
morphisms of Hodge structures such that
$$
l_1^i: L^{-i,j} \xrightarrow{\sim} L^{i,j}\qquad
l_1^j: L^{i,-j} \xrightarrow{\sim} L^{i,j}
$$
are isomorphisms. The tuple $(V,l_1, l_2)$ is called a \emph{bigraded
  Hodge-Lefschetz module of weight $p$}. A form
$\psi: L \otimes L \to \Q(-p)$ is called \emph{polarization} of $V$ if
\begin{enumerate}
\item $\psi$ restricts to morphisms of Hodge structures on
  $L^{-i,-j} \otimes L^{i,j}$;
  
\item $\psi(l_1 x, y) + \psi(x, l_1 y) = \psi(l_2 x, y) + \psi(x, l_2 y) = 0,$
\item  the restriction of $\psi(-, l_1^i l_2^j -)$ to
$L^{-i,-j} \cap \Ker l_1 \cap Ker l_2$ is a polarization of a Hodge
structure.
\end{enumerate}

\begin{rem}
  We define Hodge-Lefschetz modules over reals following the
  convention of \cite{gna90}, though we will only apply the results
  of \cite{gna90} to rational modules (tensoring with $\R$).
\end{rem}

\begin{fact}[Proposition~3.6, \cite{gna90}]
  \label{hl-module}
  If $\omega \in H^2(Y)$ restricts to a Lefschetz class on all strata
  $Y_\sigma$ then $(K^{\bu,\bu} \otimes \R, N, L_\omega)$ is a
  bigraded Hodge-Lefschetz module.
\end{fact}

A map $d: V^{\bu,\bu} \to V^{\bu+1,\bu+1}$ is a called
\emph{differential} on a polarized Hodge-Lefschetz module
if
\begin{enumerate}            
\item $d \circ d = 0$;
\item $[d, l_1] = [d, l_2] = 0$;
\item $\psi(dx,y)=\psi(x, dy)$.
\end{enumerate}

\begin{fact}[Proposition~3.5, \cite{gna90}]
  \label{hl-diff}
  If $\omega$ restricts to Lefschetz classes on all strata $Y_\sigma$
  then $d=d'+d''$ is a differential on the Hodge-Lefschetz module
  $(\Kk^{\bu,\bu}(Y) \otimes \R, N, L_\omega)$.
\end{fact}

A real Hodge-lefschetz module $L$ is naturally endowed with an
action of a Lie group
$$
\sigma: \SL(2,R) \times \SL(2,R) \to \GL(L) \qquad
d\sigma\left(
  \left(
    \begin{array}{cc}
      0 & 1 \\
      0 & 0 \\
    \end{array}
  \right), 0
\right)=l_1
\textrm{ and }
d\sigma\left(0,
  \left(
    \begin{array}{cc}
      0 & 1 \\
      0 & 0 \\
    \end{array}
  \right)
\right)=l_2
$$
Let
$$
w = \left(
  \begin{array}{cc}
    0 & 1 \\
    -1 & 0\\
  \end{array}
\right),
\qquad w_2 = (w, w) \in \SL(2,\R) \times \SL(2, \R)
$$
and define operators
$$
d^t = w_2^{-1} d w_2 \qquad \Box  = d d^t + d^t d
$$              
              
\begin{fact}[Theor\`eme~4.5, \cite{gna90}]
  \label{harmonic}  
  $\Ker \Box$ is naturally isomorphic to $H^\bu(V,d)$.
\end{fact}

Facts~\ref{hl-diff} and \ref{harmonic} together imply that every
cohomology class of $sD^{\bu}(Y)$ contains a unique representative
that belongs to $\Ker \Box$.

\subsection{Superforms}

We recall here the definitions of the superforms on a real vector
space introduced by Lagerberg following the ideas of Berndtsson. We
refer to the articles \cite{lag11, lag12} for the basic facts about
superforms and their integrals. See also \cite{gubler16, cld} for the
adaptation of the theory of superforms to the setting of tropical and
non-archimedean geometry.

Let $V$ be a real vector space. Consider the direct sum of two copies
of $V$, which we will denote $\bar V=V^0 \oplus V^1$ and define
differential superforms on $V$ to be the space of differential forms
on $\bar V$ which are invariant under translations along vectors from
$V^1$. Denote the sheaf of $n$-superforms on $V$ as $\Aa^n$.

The direct sum decompositionof $\bar V$ lifts to the tangent space,
denote
$$
J: TV^0 \oplus TV^1 \to TV^1 \oplus TV^0
$$ the morphism that
permutes the direct summands. The morphism $J$ then also acts on the
superforms. Denote the subsheaf of the superforms on $V$ that are
pull-backs of the real forms on $V^0$ as $\Aa^{1,0}$ and define
sheaves $\Aa^{0,1}$ and $\Aa^{p,q}$ as follows: for every open
$U \subset V$ set
$$
\Aa^{0,1}(U) = J \Aa^{1,0}(U) \qquad \Aa^{p,q}(U) = \bigwedge^p \Aa^{1,0}(U)
\otimes_{\Cc^\infty(U)} \bigwedge^q \Aa^{0, 1}(U) \subset \Aa^{p+q}(U).
$$
The sections of $\Aa^{p,q}$ are called \emph{$(p,q)$-superforms}. We
have
$$
\Aa^{p,q}(U) \cong \Cc^\infty(U) \otimes_{\R} \bigwedge^p T^* V^0
\otimes_{\R} \bigwedge^q T^* V^1,
$$
where the isomorphism is defined as follows. Let $x_1, \ldots, x_n$ be
a set of affine functions on $U$ such that $d'x_1, \ldots, d'x_n$ define
a basis of $T^*_x V \cong T^*V$ for any $x \in U$. Then any
$(p,q)$-superform $\alpha$ on $U$ is given by an expression
$$
\alpha = \sum_{|I|=p, |J|=q} f_{IJ}(x_1, \ldots, x_n) d'x_I \wedge d''x_J
$$
where $d'x_I = d'x_{i_1} \wedge \ldots \wedge d'x_{i_k}$ for any index
set $I=\set{i_1, \ldots, i_k}, i_1 < \ldots < i_k,$ and similarly for
$d''x_J$. Here, $f_{IJ} \in \Cc^{\infty}(U), d'x_I \in T^*V, d''x_J
\in T^*V$.

By definition, $J$ is an isomorphism of sheaves between $\Aa^{p,q}$
and $\Aa^{q,p}$.  Denote the map of sheaves induced by the de Rham
differential as
$$
d': \Aa^{p,q} \to  \Aa^{p+1, q}
$$
and define $d''= J \circ d' \circ J$. Then $d''$ maps $\Aa^{p,q}$ to
$\Aa^{p,q+1}$. One easily checks that $d',d'',J$ are given by the
following formulas

\begin{eqnarray*}\label{diff-super}
d'\alpha & = & \sum_{|I|=p, |J|=q} \dfrac{\partial f_{IJ}}{\partial x_i} d'x_ii
               \wedge d'x_I \wedge d''x_J\\
d''\alpha & = & (-1)^p \sum_{|I|=p, |J|=q} \dfrac{\partial f_{IJ}}{\partial x_i} d'x_ii
               \wedge d'x_I \wedge d''x_j \wedge d''x_J \\
J\alpha & = & (-1)^{pq} \sum_{|I|=p, |J|=q} f_{IJ} d'x_J \wedge d''x_I\\
\end{eqnarray*}

Denote $\ev: TV \otimes T^* V \to \R$ the natural pairing between
tangent vectors and covectors. The coevaluation morphism $\coev: \R
\to TV \otimes T^* V$ is the unique morphism that makes the
compositions of the following maps
$$
\begin{array}{l}
T^* V \xrightarrow{\id \otimes \coev}  T^*V \otimes (TV \otimes T^*V)
\cong (T^*V \otimes TV) \otimes T^*V \xrightarrow{\ev \otimes
  \id} T^*V\\
TV \xrightarrow{\coev \otimes \id}   (TV \otimes T^*V) \otimes TV 
\cong TV \otimes (T^* V \otimes TV) \xrightarrow{\id \otimes
  \ev} TV\\  
\end{array}
$$
identities.

\begin{defn}[Monodromy morphism]
  Let $V$ be a real vector space. For any open set $U \subset V$
  define the morphism $N: \Aa^{p,q}(U) \to \Aa^{p-1,q+1}(U)$ to be the
  composition of the morphisms
  $$
  \begin{array}{ll}
    C^\infty(U) \otimes \wedge^p T^*V \otimes \wedge^q T^*V & \to
    C^\infty(U) \otimes \wedge^p T^*V \otimes (TV \otimes T^* V)
    \otimes \wedge^q T^*V  \\
    & \to C^\infty(U) \otimes (\wedge^p T^*V \otimes TV) \otimes (T^*
    V
    \otimes \wedge^q T^*V) \\
    & \to C^\infty(U) \otimes \wedge^{p-1} T^*V \otimes \wedge^{q+1}
    T^*V\\
  \end{array}
  $$
  where the first map is given by the coevaluation map.
\end{defn}
The definition of the morphism $N$ is due to Yifeng Liu~\cite{liu17}).

The explicit formula for the  morphism $N$ is the following:
$$
N(\sum_{|I|=p, |J|=q} f_{IJ} d'x_I \wedge d''x_J) = \sum_{|I|=p,
  |J|=q} \sum_{i \in I} (-1)^{p-1}\sgn(i,I) d'x_{I \wo \{i\}} \wedge
d''x_i \wedge d''x_J.
$$

\section{Sheaves $\Lambda^p$}
\label{s:lambda-p}

In this section we fix a unipotent snc degeneration $f: X \to S$ with
the central fibre $Y = \sum\limits_{i\in I} N_i Y_i$, where $I$ is a
finite ordered set. For any face $\tau$ of the dual intersection
complex of $Y$ we will denote $N_\tau = \prod_{i \in \ul\tau} N_i$.

To simplify the exposition, we will assume that all intersections
$Y_{i_1} \cap \ldots \cap Y_{i_k}$ are connected and smooth, and so
$\Delta_X$ is a simplicial complex. With a little bit of effort the
definitions in this section can be generalized to
$\Delta$-complexes. We justify this choice by the fact that any
$\Delta$-complex can be turned into a simplical complex after a
subdivision corresponding to some blow-ups of the strata of $Y$, and
the the definition of the sheaves $\Lambda^\bu$ is invariant under
subdivisions by Proposition~\ref{subdiv}.

\subsection{Definition of $\Lambda^p$}
\label{sec:sheaves}

For any face $\sigma$ denote
$\lambda^\bu(\sigma) = \bigwedge^\bu \Q\ul\sigma/\Jj^\bu(\sigma)$,
where $\Jj^\bu(\sigma)$ is the homogeneous ideal generated by
$\sum_{i \in \ul\sigma} i$. For each $p$ the vector space
$\lambda^p(\sigma)$ is naturally isomorphic to the vector space of
translation invariant differential $p$-forms with rational
coefficients on the affine space $\la \sigma \ra$, and the space of
derivations of $\bigwedge^\bu \Q\ul\sigma$ that preserves
$\Jj^\bu(\sigma)$ is identified with the tangent space
$T\la \sigma \ra$.  Under this interpretation, for any pair of faces
$\sigma \subset \tau$ there exists a natural restriction of
differential forms map
$a \in \lambda^\bu(\tau) \mapsto a|_{\sigma} \in
\lambda^\bu(\sigma)$. In concrete terms, it can be defined on the
primitive tensors as follows:
$$
(i_1 \wedge \ldots \wedge i_l)|_{\sigma} = \left\{\begin{array}{ll}
    i_1 \wedge \ldots \wedge i_l &, \textrm{ if } \forall l\  i_l \in
                                   \ul\sigma\\
    0 & \textrm{ otherwise}.                                              \end{array}\right.
$$
Now for any face $\sigma$ define
$$
\bar \Lambda^\bu(\sigma) = \suchthat{ (f_\tau) \in \bigoplus_{\tau
    \supset \sigma} \lambda^\bu(\tau) }{ \forall \tau \supset \alpha
  \supset \sigma, \ f_\tau|_\alpha = f_\alpha }
$$
This defines a functor: the morphism associated to an inclusion
$\sigma \supset \tau$ is the natural projection
$\bar \Lambda^\bu(\sigma) \hookrightarrow \bar \Lambda^\bu(\tau)$
coming from the inclusion $\Stbar(\tau) \subset \Stbar(\sigma)$. We
will denote it $a \mapsto a|_{\Stbar(\tau)}$.

The space $\Lambda^p(\sigma)$ is spanned by primitive tensors of the
form $i_1 \wedge \ldots \wedge i_p$

For any face $\sigma, |\ul\sigma| \geq p$ define the map
$c^p_\sigma: \bar \Lambda^p(\sigma) \to H^2(Y_\sigma) \otimes
\lambda^{p-1}(\sigma)$ on the primitive tensors by the formula
$$
c^p_\sigma(i_1 \wedge \ldots \wedge i_p) = \sum_{l=1}^p (-1)^{l-1}
c_1(\Oo(Y_{i_l}))|_{Y_\sigma} \otimes (i_1 \wedge \ldots \wedge
\what{i_l} \wedge \ldots \wedge i_p)|_{\sigma}.
$$

The maps $c^p_\sigma$ are coordinate components of the differential of
the complex $K^0(\sigma, \Lambda^p)$ that will be introduced in
Section~\ref{sec:k-diff}. They are well-defined by Lemma~\ref{ind-k}.

\begin{lemma}
  \label{c-deriv}
  For any face $\sigma$ the map $c^\bu_\sigma$ is a derivation, that
  is, for any $a \in \bar \Lambda^p, b \in \bar \Lambda^q$
  $$
  c^{p+q}_\sigma(a \wedge b) = c^p_\sigma(a) \wedge b + (-1)^p a
  \wedge c^q_\sigma(b).
  $$
\end{lemma}

\begin{proof}
  Immediate from the definition.
\end{proof}

\begin{defn}
  \label{lambdap-defn}
  For any $p \geq 1$ define functors $\Lambda^p$ as follows. For
  all faces $\sigma, |\sigma| \geq k$ let
  $$
  \Lambda^p(\sigma) = \suchthat{ \alpha \in \bar \Lambda^p}{ \forall
    \tau \supset \sigma, \alpha|_\tau \in \Ker c^p_\tau } ,
  $$
  and denote the sheaves on $\Delta$ corresponding to these functors
  by Lemma~\ref{const-shv} as $\Lambda^p$. We also agree to denote
  the constant sheaf $\Q$ as $\Lambda^0$.
\end{defn}

\begin{lemma}
  For any $p > 1$ let
  $$
  U_p = \bigcup_{|\ul\tau| \geq p} \St(\tau).
  $$
  and $j: U_p \hookrightarrow \Delta$ be the open embeddings. Then
  $$
  \Lambda^p = j_* j^* \Lambda^k.
  $$
\end{lemma}

\begin{proof}
  Follows immediately from Definition~\ref{lambdap-defn}, as
  $$
  \Ker c^p_\sigma = \bar \Lambda^p(\sigma)
  $$
  for all $\sigma$ such that $|\ul\sigma|<p$, since
  $\lambda^{p-1}(\sigma)=0$ for such faces $\sigma$.
\end{proof}

\begin{lemma}
  For all $p,q$ for which the corresponding sheaves are defined
  $$
  \Lambda^p \wedge \Lambda^q \subset \Lambda^{p+q}.
  $$
\end{lemma}

\begin{proof}
  This amounts to showing that for any face $\alpha$ and any faces
  $\sigma, \tau \subset \alpha$,
  $$
  a \in \Lambda^p(\sigma), b \in \Lambda^q(\tau) \Rightarrow
  a|_{\Stbar(\alpha)} \wedge b|_{\Stbar(\alpha)} \in
  \Lambda^{p+q}(\alpha),
  $$
  which follows immediately from Lemma~\ref{c-deriv}. 
\end{proof}

\subsection{Regularity and sheaves $A^p$}

Consider the sheaf of  piece-wise affine functions on $\Delta_X$
that restrict to affine functions on the faces $\tau \supset \sigma$
--- such functions are determined by their values at the vertices:
$$
\bar A^1(\sigma) = \set{ f: \Stbar^0(\sigma) \to \Q }
$$
Consider further a subsheaf of $\bar A^1$ that corresponds to the
functor
$$
A^1(\sigma) = \suchthat{ f \in \bar A^1(\sigma)}{ \sum_{i \in
    \Stbar^0(\sigma)} N_i f(i) \cdot c_1(\Oo(Y_i))|_{Y_\sigma} = 0}
$$
This sheaf have been defined in a note \cite{kt} of Konstevich and
Tschinkel (where it is called the sheaf of ``linear functions'',
Sect.~6.2), see also \cite{tony, jsr-epi}. One can observe easily that
for any face $\sigma$ we have a functorial exact sequence
$$
0 \to \Lambda^0(\sigma) \to A^1(\sigma) \to \Lambda^1(\sigma) \to 0.
$$
that gives rise to the exact sequence of sheaves
$$
  0 \to \Lambda^0 \to A^1 \to \Lambda^1 \to 0.
$$
One can wonder under which circumstances an analogous sequence exists
in higher degrees.

Take a face $\sigma \subset \Delta_X$. There exists a natural map given
by the evaluation function
$$
e_\sigma: \St(\sigma) \to H^0(\St(\sigma), A^1)^* \qquad x \mapsto [f
\mapsto f(x)]
$$
which factors through the subspace
$$
T(\sigma) := \suchthat{ F \in H^0(\St(\sigma), A^1)^*}{ \forall c\in
  \R\, F(c) = c }.
$$
The latter is clearly a torsor under the vector subspace
$$
H^0(\St(\sigma), \Lambda^1)^* \cong \suchthat{ F \in H^0(\St(\sigma),
  A^1)^*}{ F(c) = 0} \subset H^0(\St(\sigma), A^1)^*.
$$
The elements of $A^1(\sigma)$ are tautologically functions on
$T(\sigma)$. The function $e_\sigma$ is affine on the simplices of
$\Stbar(\sigma)$, therefore $e_\sigma(\Stbar(\sigma))$ is a simplicial
complex in $T(\sigma)$ and $e_\sigma(\St(\sigma))$ is the star of the
simplex $e_\sigma(\sigma)$.

\begin{defn}[Regularity]
  \label{regular}
  Let $\Delta$ be a simplicial complex endowed with a constructible
  sheaf $\Lambda^\bu$ of graded algebras and a subsheaf $A^1$ of the
  sheaf of piece-wise affine functions on $\Delta$ such that
  $A^1/\const = \Lambda^1$. We say that $\Lambda^\bu$ is
  \emph{regular} at a face $\sigma$ if
  $$
  \Lambda^p(\sigma) = e_\sigma^* \bigwedge^p \Lambda^1_{T(\sigma)}.
  $$
\end{defn}

If $\Lambda^\bu$ is regular at a face $\sigma$ we can define for any
$p \geq 0$ the sheaf
$$
A^p(\sigma) = e_\sigma^* \bigwedge^p A^1(\sigma)
$$
that fits into the exact sequence
$$
0 \to \Lambda^p(\sigma) \to A^{p+1}(\sigma) \to \Lambda^{p+1}(\sigma)
\to 0
$$
If $\Lambda^\bu$ is regular at every face
$\sigma \subset \Delta_X$ then $A^p$ defines a functor and a sheaf on
$\Delta_X$, and we have an exact sequence of sheaves
\begin{equation}  
  0 \to \Lambda^p \to A^{p+1} \to \Lambda^p \to 0.
  \label{seq:lambda'}
\end{equation}

\subsection{Invariance under subdivision}

Let $f: X' \to X$ be a blow-up of a stratum $Y_\sigma$. The dual
intersection complex $\Delta_X$ is naturally homeomorphic to
$\Delta(Y')$ and we will identify the underlying topological spaces
from now on; let $\xi$ denote the homeomorphism. The simplicial
complex $\Delta(Y')$ is obtained from $\Delta_X$ by adding the vertex
$e$ corresponding to the exceptional divisor to $\sigma$ and
subdividing all faces that contain $\sigma$. For every face
$\tau \not\supset \sigma$ we denote the corresponding face of
$\Delta(Y')$ as $\tau$ too. If $\tau$ is a face of $\Stbar(\sigma)$
that does not contain $\sigma$ then there exists a unique face of
$\Delta(Y')$, that we will denote $\tau'$, that contains $\tau$ and
the vertex $e$. We will allow $\tau$ to be empty for uniformity, in
which case we agree that $\tau'=e$.

The face $\sigma$ is identified with the intersection of the positive
octant in the vector space $\R^m$,
$m=|\ul \sigma|, \sigma=\set{i_1, \ldots, i_m}$, with the affine
subspace cut out by the equation
$$
\sum_{k = 1}^\infty N_{i_k} x_k = 1
$$
and the vertex $e$ with the point 
$\left( \dfrac{1}{N_\sigma}, \ldots, \dfrac{1}{N_{\sigma}}
\right)$.

Let $\beta$ be a face in $\Stbar(\sigma)$ and let
$\eta \subseteq \beta$ be a face of $\Delta(Y')$ such that
$f(Y_\eta) = Y_\beta$. In particular,
$\St(\eta) \subset \St(\beta)$; denote this open embedding
$\xi_{\beta,\eta}$. There are two options: either
$\beta \not\supset \sigma$ and $\eta=\beta$, or $\beta \supset \sigma$
and $\eta = \alpha'$ for some $\alpha \subset \beta$. In both cases
$\xi_{\beta,\eta}^*: \bar A_X^1(\beta)=H^0(\St(\beta), \bar A^\bu_X)
\hookrightarrow \bar A_{X'}^1(\eta)=H^0(U_\eta, \bar A^\bu_{X'})$ is
injective. Indeed, if an element of $\bar A^1_X(\beta)$ is represented
by a piece-wise affine function $f: \Stbar^0(\beta) \to \R$, then to
describe its restriction to $\Stbar^0_{X'}(\eta)$ it suffices to
calculate the value of $f$ in $e$, which is
$$
f(e) = \dfrac{1}{\sum_{i \in \ul\sigma}N_i} \sum_{i \in \ul\sigma} f(i) = \dfrac{1}{N_e} \sum_{i \in \ul\sigma} f(i) 
$$

\begin{lemma}
  \label{h2-of-strata}
  For any stratum $Y'_\eta \subset Y'_e$,
  $$
  H^2(Y'_\eta) \cong H^2(Y_\sigma) \oplus \Q c_1(\Oo(Y'_e))|_{Y'_\eta}
  $$
\end{lemma}

\begin{proof}
  The exceptional divisor $Y'_e$ is a
  projective bundle over the stratum $Y_\sigma$. Since $Y$ is an snc
  divisor, the stratum $Y'_\tau$ is a projective subbundle of
  $Y'_e$. The statement of the Lemma then follows from 
  \cite[Lemma~7.32]{voisin}. 
\end{proof}

\begin{prop}
  \label{subdiv}
  Let $\beta$ beta a face in $\Stbar(\sigma)$ and let
  $\eta \subseteq \beta$ be a face of $\Delta(Y')$ such that
  $\St(\eta) \subset \St(\beta)$. Denote
  $r_{\beta,\eta}: \lambda^{\bu}(\beta) \to
  \lambda^{\bu}(\eta)$. Then the diagram
  $$
  \xymatrix{ \bar \Lambda^p_X(\beta) \ar[d]^{\xi^*_{\beta,\eta}}
    \ar[rr]^{c^p_\beta} & &
    H^2(Y_\beta) \otimes \lambda^{p-1}(\beta) \ar[d]^{f^* \otimes r_{\beta,\eta}}
    \\
    \bar \Lambda^p_{X'}(\eta) \ar[rr]^{c^p_\eta} & & H^2(Y_\eta) \otimes \lambda^{p-1}(\eta) }
  $$
  is commutative, and
  $\Ker c^p_\eta = \xi_{\beta,\eta}^* \Ker c^p_\beta$. In particular,
  the definition of sheaves $A^p$ and $\Lambda^p$ is invariant under
  the subdivisions induced by blow-ups of the strata of $Y$.
\end{prop}

\begin{proof}
  Take a basis element
  $i \in \bar \Lambda^1(\beta)$, then
  $$
  \xi_{\beta,\eta}^*(i) = i + \dfrac{N_i}{N_e}.
  $$  
  The diagram is commutative for $p=1$ since
  $$
  c^1_\beta(\xi^*_{\beta, \eta}(i)) = c_1(i + \dfrac{N_i}{N_e}e) =
  N_i (c_1(\Oo_{X'}(Y'_i)) + c_1(\Oo_{X'}(Y'_e)))=
  f^* N_i c_1(\Oo_X(Y_i)) = f^* c^1_\eta (i).
  $$
  The commutativity of the diagram for $p > 1$ follows by induction
  from Lemma~\ref{c-deriv} and commutativity of the diagram
  $$
  \xymatrix{ \bar \Lambda^p(\beta) \ar[r] \ar[d]^{\xi_{\beta,\eta}^*}
    &
    \lambda^p(\beta) \ar[d]^{r_{\beta,\eta}} \\
    \bar \Lambda^p(\eta) \ar[r] & \lambda^p(\beta) \\
  }
  $$
  where the horizontal morphisms are restriction maps, for all
  $p \geq 1$.

  The second statement follows from the fact that
  $\bar \Lambda^1(\eta) = \xi^*_{\beta,\eta} \bar \Lambda^1(\beta)
  \oplus \Q e$, from Lemma~\ref{h2-of-strata} and from the inclusion
  $$
  c^1_\beta(\xi^*_{\beta,\eta} \bar \Lambda^1(\beta)) \subset f^*
  H^2(Y_\beta) \subset H^2(Y_\eta),
  $$
  that has just been established.
\end{proof}

In view of Proposition~\ref{subdiv} the following definition is
consistent.

\begin{defn}[Sheaves $\Lambda^\bu$ in the normal crossings
  case]
  \label{sheaves-non-snc}
  For any degeneration $f: X \to S$ such that reduction of its central
  fibre is normal crossings but not necessarily simple normal
  crossings, define the sheaves $\Lambda^\bu$ on $\Delta_X$ to be
  those of a modification $X'$ of $X$ obtained by a sequence of
  blow-ups of the strata, and such that the central fibre is strictly
  normal crossings and $\Delta_{X'}$ is a simplicial complex.
\end{defn}

\subsection{Toric strata}

Let $X$ be a toric variety acted upon by a torus
$T = \C[M], M \cong \Z^n$, let $N = \Hom(M, \Z)$ and let
$\Sigma \subset N \otimes \R$ be the toric fan of $X$. Let
$t: X \to \A^1$ be an equivariant morphism, assume that
$Y=t^{-1}(0)$ is an snc divisor and denote the corresponding linear
function on $N$ as $\tilde{t}$.  The dual intersection complex of
$t^{-1}(0)$ can be identified with the polyhedral complex
$\tilde{t}^{-1}(1)$. We will denote the intersections of cones $\sigma$
with $\tilde{t}^{-1}(1)$ as $\sigma$ too.

\begin{prop}
  \label{toric}
  For any toric stratum $X_\sigma \subset Y$,
  $\Lambda^\bu(\sigma)$ is the algebra of translation-invariant
   differential forms on $\Stbar(\sigma) \subset \tilde{t}^{-1}(1)$ with
   rational coefficients.
\end{prop}

\begin{proof}
  The group of Cartier $\Q$-divisors on $X_\Sigma$ is isomorphic to
  the group of piece-wise linear functions with rational coefficients
  on $\Sigma$ that restrict to a linear function on each cone
  $\eta \subset \Sigma$; the linear functions among them correspond to
  the principal divisors. The stratum $X_\sigma$ is a toric variety in
  its own right, with the dense torus $T_\sigma =
  \C[\sigma^\perp]$. We may identify
  $\Hom((\sigma \cap \Q)^\perp, \Q)$ with $N/\la \sigma \ra$, where
  $\la \sigma \ra$ is the vector space spanned by by $\sigma \cap
  \Q$. Given a $\Q$-divisor $D$ on $X_\Sigma$ corresponding to a piece-wise
  linear function $\tilde{t}_D: N \otimes \R \to \R$, the restriction of
  its linear equivalence class $X_\sigma$ can be decribed as
  follows. There exists a linear function
  $\tilde{t}'_D: N \otimes \R \to \R$ such that $\tilde{t}_D - \tilde{t}'_D$ vanishes
  on $\la \sigma \ra$, this function is therefore a pull-back of a
  piece-wise linear function with rational coefficients
  $\eta: N/\la \sigma \ra \to \R$; the function $\eta$ corresponds to
  the linear equivalence class of the restriction of the linear
  equivalence class of $D$ to $X_\sigma$.

  The elements of $\bar A^1(\sigma)$ by definition correspond to
  piece-wise linear functions with rational coefficients on the cone
  over $\Stbar(\sigma)$, or to affine functions with rational
  coefficients on $\tilde{t}^{-1}(1)$.  It then immediately follows
  from the description of the restriction of toric divisors to
  $X_\sigma$ that
  $$
  \Lambda^1(\sigma) = \Aff(\Stbar(\sigma))/\const.
  $$
  By Lemma~\ref{c-deriv} we have
  $$
  \bigwedge^p \Lambda^1 \subset \Lambda^p
  $$
  To conclude, it will suffice show for any codimension 1 face
  $\tau \subset \Stbar(\sigma)$ that
  $$
  \Lambda^p(\tau) \subset \bigwedge^p \Lambda^1(\tau)
  $$
  Since $\tau$ is a codimension 1 face, it has two faces of maximal
  dimension contianing it, call them $\eta_1$ and $\eta_2$, $\ul\eta_1
  = \ul\tau \cup \{i\}$ and $\ul\eta_2 = \ul\tau \cup \{j\}$. Any
  function $f \in A^1(\tau)$ detemines an affine function on $\eta_1$, 
  and since $\eta$ is of maximal dimension, this linear function
  extends uniquely to an affine function on $\Stbar(\tau)$. It follows
  that
  $$
  \bar \Lambda^1(\tau) = \Lambda^1(\tau) \oplus \Q j, \qquad \bar \Lambda^p(\tau) =
  \bigwedge^p \Lambda^1(\tau) \oplus \Q j \wedge \left( \bigwedge^{p-1}
    \Lambda^1(\tau) \right)
  $$
  By Lemma~\ref{c-deriv} we have for any $a \in \bar \Lambda^{p-1}(\tau)$
  $$
  c^p_\tau(j \wedge a) = c^1_\tau(j) \wedge a + (-1)^{p-1}
  c^{p-1}_\tau(a) \wedge j = c^1_\tau(j) \wedge a,
  $$
  from which it is clear that that $c^p_\tau$ cannot vanish on the
  elements of the second direct summand, and we can conclude.
\end{proof}

\begin{cor}
  \label{toric-tubular}
  Let $f: X \to S$ be an snc degeneration, $Y=f^{-1}(0)$, and let
  $g: X_\sigma \to \A^1$ be a morphism of toric varities as in the
  statement of Proposition~\ref{toric}. Assume that a tubular
  neighbourhood of a stratum $Y_\sigma$ of $Y$ is isomorphic to a
  tubular neighbourhood of a stratum of $X_\tau \subset
  g^{-1}(0)$. Then $\Lambda^\bu(\sigma)$ is isomorphic to the algebra
  of translation-invariant rational differential forms on
  $\Stbar(\tau)$.
\end{cor}

\section{The complex $K^{\bu}(\Delta_X, \Lambda^p)$}
\label{s:comb}

In this section we will study a certain complex that will play the
crucial role in the proof of Theorem~A. Let us outline our strategy,
which is inspired by \cite[Section~5]{ikmz} (though since Theorem~A is
about cohomology and not homology, our set up is dual; we also do not
have to deal with sedentarity and infinite faces). We construct
resolutions $K^{\bu}(\sigma, \Lambda^p)$ of vector spaces
$\Lambda^p(\sigma)$, functorial in $\sigma$, and use them to define a
complex $K^{\bu}(\Delta_X, \Lambda^p)$ quasi-isomorphic to the complex
$C^{\bu}(\Delta_X, \Lambda^p)$ of simplicial cochains with
coefficients in $\Lambda^p$. We then show that
$C^{\bu}(\Delta_X, \Lambda^p)$ has the same cohomology as a certain
subcomplex $s\Dd^\bu_p(Y)$ of $sD^\bu_p(Y)$, the $2p$-th line of the
Steenbrink spectral sequence.

For the benefit of a reader familiar with the proof of the main result
of \cite{ikmz} we discuss technical peculiarities that pertain to our
set up. Firstly, the definition of the complex
$K^{\bu}(\Delta_X, \Lambda^p)$ uses only cohomology classes that lie
in the span of cycle classes of the strata of $Y$. Its differential
takes into account the multiplicities of the components of the central
fibre as opposed to \cite{ikmz}, where a unimodular triangulation of
the tropical limit is chosen first. This complex plays a role similar
to that of the complex $K^{(p)}_{\bu,\bu}$ of \cite[5.3-4]{ikmz}, but
in order to show that it is quasi-isomorphic to
$C^{\bu}(\Delta_X,\Lambda^p)$ we use an inductive argument in
Lemma~\ref{h0m} instead of the purity of the weight filtration on the
cohomology of the complement of a hyperplane arrangement that was used
in \cite{ikmz}. Secondly, the treatment of signs in the computations
of the differentials of the zeroth and first pages of the spectral
sequence associated to the filtration $F^\Delta$ on the cochain
complex $K^{(p)}_\bu$ in Section~5.4 of \cite{ikmz} contains gaps. In
particular, the components of the differential of the residue complex
should have signs that depend on the destination face. Filtration $G$
on the cochain complex $K^\bu(\Delta_X, \Lambda^\bu)$ is the analogue
of $F^\Delta$ in our set up, and in our Lemma~\ref{acyclic} the
correct signs for the differential in the complex $S^\bu_n(\sigma)$,
the dual of the residue complex of \cite{ikmz}, are computed. These
signs are then used to derive correct signs for the differentials and
for the coordinate components of the map $N$ in the proofs of
Lemmas~\ref{e0-of-ss}, \ref{e1-of-ss}.

\subsection{Complexes $\Mm_q^{\bu}(Y)$ and $\Dd_p^{\bu,\bu}(Y)$}
\label{s:complexes}

For any $k$ denote $Y_{\neq k}=\sum_{i \neq k} Y_i$ and denote
$Y^\circ_k$ the divisor $Y_k \cap Y_{\neq k}$ in $Y_k$. More
generally, if $Y_\sigma \subset Y$ is a stratum of $Y$, denote
$$
Y_\sigma^\circ = \bigcup_{j \notin \ul\sigma} Y_j \cap Y_\sigma.
$$

\begin{defn}
  For any stratum $Y_\sigma \subset Y$ denote
  $$
  \Hh^\bu(Y_\sigma) = \Q
  \suchthat{\Gys^\tau_\sigma(1_{Y_\tau})}{\tau \in\St(\sigma)}
  \subset H^\bu(Y_\sigma).
  $$  
\end{defn}

More generally, denote
$$
\Hh^\bu(Y^{(p)}) = \bigoplus\limits_{|\ul\sigma|=p+1} \Hh^\bu(Y_\sigma)
$$
Note that $\Hh^\bu(Y^{(p)})$ is a subring of
$H^\bu(Y^{(p)})$.  Denote $\Mm^\bu_p(X, Y)$ the subcomplex of
$M^\bu_p(X, Y)$:
$$
0 \to \Hh^0(Y^{(p)}) \xrightarrow{\gamma^{(p)}} \Hh^2(Y^{(p-1)})
\to  \ldots \to \Hh^{2p-2}(Y^{(1)})
\xrightarrow{\gamma^{(1)}} \Hh^{2p}(X) \to 0
$$
Similarly, denote $\Dd_p^{\bu,\bu}(Y)$ the complex
$$
0 \to \sigma_{\leq p}\Mm^\bu_p(X,Y) \xrightarrow{\rho^{(p)}} \sigma_{\leq
  p}\Mm^\bu_{p+1}(X,Y) \xrightarrow{\rho^{(p+1)}} \sigma_{\leq
  p}\Mm^\bu_{p+2}(X,Y) \to \ldots
$$

\begin{lemma}
  \label{h0m}
  For any variety $X$ and an snc divisor $Y \subset X$, for all
  $p \geq 0$ the complex $\Mm^\bu(X, Y)$ doesn't
  have cohomology in degree $> 0$.
\end{lemma}

\begin{proof}
  The proof is by induction on $p$. Indeed, $\Mm^\bu_0(X, Y)$ consists
  a single term and so the statement of the Lemma is true for it. It
  follows from the definition of $\Mm^\bu_p(X,Y)$ that the following
  sequence is exact
  $$
  0 \to \Mm^\bu_p(X, Y_{\neq k}) \to \Mm^\bu_p(X,Y) \to
  \Mm^\bu_{p-1}(Y_k, Y_k^\circ) \to 0.
  $$
  By induction hypothesis the cohomology of the first and third terms
  is concentrated in degree 0, therefore, the same is true for 
  $\Mm^\bu_p(X, Y)$.
\end{proof}

\subsection{Complex $K^{\bu}(\sigma, \Lambda^p)$: differential}
\label{sec:k-diff}

Consider a face $\sigma \subset \Delta_X$. The space $\lambda^1(\sigma)$
has a number of distinguished bases: for each $k \in \ul\sigma$, the
vectors of the form $dx_{ik} = i - k \in \lambda^1(\sigma), i \neq k$
are linearly independent. These bases have dual ones in
$T\la \sigma \ra$: for any $i \in \ul\sigma$ let
$\d_i \in T\la \sigma \ra$ be the derivation determined by
$$
\d_i(i) = 1 \qquad \d_i(j) = 0, j \neq i
$$
and let $\d_{ik} = \d_i - \d_k$. In other words, $\d_{ik}$ is the
contraction with the vector going from vertex $k$ to the vertex $i$.

\begin{lemma}
  \label{deriv}
  If $\tau \supset \sigma, i, k \in \ul\sigma$ then for all
  $a \in \lambda^\bu(\tau), (\d_{ik} a)|_{\sigma} = \d_{ik}
  (a|_\sigma)$.
\end{lemma}

\begin{proof}
  Straightforward.
\end{proof}

\begin{defn}
  For any face $\sigma$ define let the term $K^r(\sigma, \Lambda^p)$
  to be the subspace of
  $$
  \bigoplus\limits_{\tau \supset \sigma} \Hh^{2r}(Y_\tau) \otimes
  \lambda^{p-r}(\tau)
  $$
  spanned by tuples of primitive tensors
  $\sum_{\tau \supset \sigma} a_\tau \otimes b_\tau$ that satisfy the
  equations
  $$
  \res_\beta(a_\alpha) = a_\beta, \qquad b_\beta|_\alpha = b_\alpha,
  $$
  for all $\beta \supset^1 \alpha$. Pick a vertex $k \in \ul\sigma$
  and define the differential 
  $$
  d'_k: \bigoplus_{\tau \supset \sigma} \Hh^{2r}(Y_\tau) \otimes
  \lambda^{p-r}(\tau) \to \bigoplus_{\tau \supset \sigma}
  \Hh^{2r+2}(Y_\tau) \otimes \lambda^{p-r-1}(\tau)
  $$
  on the primitive tensors
  by the formula
  $$
  d'_k(a_\tau \otimes b_\tau) = \sum_{i \in \ul\tau \wo \ul\sigma} N_i
  \Gys^\tau_{\6_i \tau}(a_\tau) \otimes (\d_{ik} b_\tau)|_{\6_i \tau}
  + \sum_{i \in \ul\sigma} N_i \res_\tau(\Gys^\tau_{\6_i
    \tau}(a_\tau)) \otimes \d_{ik} b_\tau
  $$
  where the first term belongs to
  $$
  \bigoplus_{i \in \ul\tau} \Hh^{2r+2}(Y_{\6_i \tau}) \otimes
  \lambda^{p-r}(\6_i \tau)
  $$
  and the second belongs to
  $$
  \Hh^{2r+2}(Y_{\tau}) \otimes \lambda^{p-r}(\tau).
  $$
\end{defn}

\begin{lemma}
   \label{ind-k}
   For any $k, k' \in \ul\sigma$ and any $v \in K^r(\sigma, \Lambda^p)$,
   $d'_k(v) = d'_{k'}(v)$.
 \end{lemma}

 \begin{proof}
  Consider a sum of primitive tensors
  $$
  v = \sum_{\tau \supset \sigma} a_\tau \otimes b_\tau \in
  K^r(\sigma, \Lambda^p),
  $$
  then we have
  \begin{multline*}
    (d'_k-d'_{k'})(v)_\alpha = \sum_{\eta: \alpha = \6_i \eta}
    N_i \Gys^\eta_\alpha (a_\eta) \otimes (\d_{k'k} b_\eta)|_{\alpha} +
    \sum_{i \in \ul\alpha} N_i \res_\alpha(\Gys^\alpha_{\6_i
      \alpha}(a_\alpha)) \otimes \d_{k'k} b_{\alpha} \\
    = \left(\sum_{\eta: \alpha = \6_i \eta} N_i \Gys^\eta_\alpha
      (\res_\eta(a_\alpha)) 
      + \sum_{i \in \ul\alpha} N_i \res_\alpha(\Gys^\alpha_{\6_i
        \alpha}(a_\alpha)) \right)\otimes \d_{k'k} b_{\alpha}, \\
  \end{multline*}
  which vanishes by Corollary~\ref{sum-res-gys}.
\end{proof}

In view of the Lemma, we will denote any of the differentials $d'_k$
as simply $d'$.

\begin{lemma}
  \label{preserves}
  $d'(K^r(\sigma, \Lambda^p)) \subset K^{r+1}(\sigma, \Lambda^p)$.
\end{lemma}

\begin{proof}
  To perform the computations below let us fix a vertex
  $k \in \ul\sigma$.  It suffices to prove that
  $d'_k v \in K^{r+1}(\sigma, \Lambda^p)$ for all sums of primitive
  tensors
  $$
  v = \sum_{\tau \supset \sigma} a_\tau \otimes b_\tau \in
  K^r(\sigma, \Lambda^p).
  $$

  Take two faces
  $\alpha, \beta \supset \sigma, \alpha = \6_l \beta$. Then
  \begin{multline*}
    (d' v)_\alpha = \sum_{\eta: \alpha = \6_i \eta, \eta \neq \beta}
    N_i \Gys^\eta_\alpha (a_\eta) \otimes (\d_{ik} b_\eta)|_{\alpha} + \\
    + N_l \Gys^\beta_\alpha(a_\beta) \otimes (\d_{lk} b_\beta)|_{\alpha} + \sum_{j
      \in \ul\alpha} N_j \res_\alpha(\Gys^\alpha_{\6_j \alpha}(a_\alpha))
    \otimes \d_{jk} b_{\alpha},
  \end{multline*}
  \begin{multline*}
    (d' v)_\beta = \sum_{\eps: \beta = \6_i \eps} N_i \Gys^\eps_\beta
    (a_\eps)
    \otimes (\d_{ik} b_\eps)|_{\beta}  +  \\
    + \sum_{j \in \ul\alpha} N_j \res_\beta(\Gys^\beta_{\6_j
      \beta}(a_\beta)) \otimes \d_{jk} b_{\beta} +
    N_l \res_\beta(\Gys^\beta_{\alpha}(a_\beta)) \otimes \d_{lk}
    b_\beta.
   \end{multline*}
   Since $\Delta_X$ is assumed to be a simplicial complex, there is a
   one-to-one correspondence between faces enumerating the summands of
   the first terms of $d(v)_\beta$ and $d(v)_\alpha$.  All
   such pairs of faces $\eps, \eta$ satisfy
   $$
   \eta = \6_l \eps, \qquad \6_i\eps=\beta, \qquad \6_i \eta=\alpha.
   $$
   We have
   \begin{gather*}
     \res_\beta(\Gys^\eta_\alpha(a_\eta)) = \Gys^\eps_\beta(\res_\eps(a_\eta)))
     =
     \Gys^\eps_\beta(a_\eps),
   \end{gather*}
   since $v \in K^r(\sigma, \Lambda^p)$ and by
   Corollary~\ref{cor-res-gys}, and we also have
   $$
   (\d_{ik} b_\eps)|_\alpha = (\d_{ik}(b_\eps|_\eta))|_\alpha=
   (\d_{ik} b_\eta)|_\alpha
   $$
   by Lemma~\ref{deriv}.
   
   Again by Corollary~\ref{cor-res-gys}, for all $j \in \ul\alpha$
   $$
   \res_\beta\res_\alpha(\Gys^\alpha_{\d_j \alpha}(a_\alpha)) =
   \res_\beta(\Gys^\beta_{\6_j \beta}(\res_\beta(a_\alpha))) =
   \res_\beta(\Gys^\beta_{\6_j \beta}(a_\beta)),
   $$
   and by Lemma~\ref{deriv},
   $$
   (\d_{jk} b_\beta)|_\alpha = \d_{jk} (b_\beta|_\alpha) = \d_{jk} b_\alpha.
   $$

   Finally, the terms
   $$
   \Gys^\beta_\alpha(a_\beta) \otimes (\d_{lk} b_\beta)|_\alpha
   \textrm{ and } \res_\beta \Gys^\beta_\alpha(a_\beta) \otimes
   \d_{lk} b_\beta
   $$
   clearly satisfy the identities from the definition of $K^r(\sigma,
   \Lambda^p)$. We conclude that $d'_k(v) \in K^{r+1}(\sigma,
   \Lambda^p)$. 
 \end{proof}

 \begin{lemma}
  \label{is-diff}
  $d' \circ d' = 0$.
\end{lemma}

\begin{proof}
  Let $v \in K^r(\sigma, \Lambda^p)$ be a sum of primitive tensors as
  before.
  
  \begin{multline*}
    d'_k(d'_k v) = \sum_{\tau \supset \sigma} \left(\sum\limits_{j \in
        \ul{\6_i \tau}} \sum\limits_{i \in \ul\tau} N_i N_j \Gys^{\6_i \tau}_{\6_i
        \6_j \tau} \Gys^\tau_{\6_i \tau}
      a_\tau \otimes (\d_{jk} \d_{ik} b_\tau)|_{\6_j \6_i \tau} + \right. \\
    + \sum_{j \in \ul\tau} \sum\limits_{i \in \ul\tau}  N_i N_j \Gys^\tau_{\6_j
      \tau}\res_\tau \Gys^\tau_{\6_i \tau}
    a_\tau \otimes (\d_{jk} \d_{ik} b_\tau)|_{\6_j\tau}  + \\
    + \sum_{j \in \ul{\6_i\tau}} \sum\limits_{i \in \ul\tau} N_i N_j \res_{\6_i
      \tau} \Gys^{\6_i \tau}_{\6_j \6_i \tau} \Gys^\tau_{\6_i \tau}
    a_\tau \otimes \d_{jk} (\d_{ik} b_\tau)|_{\6_i\tau}  + \\
    \left.+ \sum_{j \in \ul\tau} \sum\limits_{i \in \ul\tau} N_i N_j \res_\tau
      \Gys^\tau_{\6_j \tau} \res_\tau \Gys^\tau_{\6_i \tau}
      a_\tau \otimes \d_{jk} \d_{ik} b_\tau\right). \\
  \end{multline*}
  Then by Corollary~\ref{cor-res-gys}
  $$  
  \Gys^\tau_{\6_j \tau}\res_\tau \Gys^\tau_{\6_i \tau} =
  \res_{\6_j \tau} \Gys^{\6_i\tau}_{\6_j\6_i \tau} \Gys^\tau_{\6_i
    \tau} = \res_{\6_j \tau} \Gys^{\6_j \tau}_{\6_i \6_j \tau}
  \Gys^\tau_{\6_j
    \tau}\\
  $$
  and
  \begin{multline*}
    \res_\tau \Gys^\tau_{\6_j \tau} \res_\tau \Gys^\tau_{\6_i \tau}
    a_\tau = \res_\tau \res_{\6_j \tau} \Gys^{\6_i \tau}_{\6_j \6_i
      \tau} \Gys^\tau_{\6_i \tau}
    a_\tau = \\
    = \res_\tau \res_{\6_i \tau} \Gys^{\6_j \tau}_{\6_j \6_i \tau}
    \Gys^\tau_{\6_j \tau} a_\tau = \res_\tau \Gys^{\tau}_{\6_i \tau}
    \res_\tau
    \Gys^\tau_{\6_j \tau} a_\tau,\\
  \end{multline*}
  and therefore by Lemma~\ref{deriv} the sum in parenthesis vanishes
  for each $\tau \supset \sigma$.
  
\end{proof}

\subsection{Complex $K^{\bu}(\sigma, \Lambda^p)$: resolution}

It is clear from the Definition~\ref{lambdap-defn} that
$$
\Lambda^p(\sigma) \cong K^0(\sigma, \Lambda^p)
$$

For any pair of  faces $\sigma, \tau, \sigma \subset \tau$ define the
map 
$$
\kappa^\tau_\sigma: \lambda^\bu(\tau) \to \lambda^{\bu-l}(\sigma)
\qquad a \mapsto (\d_{i_l,\tau} \circ \ldots \circ \d_{i_1,\tau})(a)|_{\sigma},
$$
where
$\ul\tau \wo \ul\sigma = \set{i_1, \ldots, i_l}, i_1 < \ldots < i_l$.

Denote $\lambda^\bu(\tau; \sigma)$ the subspace of
differential forms in $\lambda^\bu(\tau)$ such that their
restriction to $\la \sigma \ra$ vanishes, in particular, for any
such form $a$ and any derivation $\d \in T\la \sigma \ra$,
$\d a = 0$. Define 
$$
F^m \lambda^p(\tau) = \bigoplus_{l=0}^{m-1} \lambda^{p-l}(\tau;
\sigma) \wedge \lambda^l(\sigma), \qquad
F^m K^r(\sigma, \Lambda^p) = \bigoplus\limits_{\tau \supset \sigma}
\Hh^{2r}(Y_\tau) \otimes F^m \lambda^{p-r}(\tau).
$$
If $a \in \lambda^{p-l}(\tau; \sigma) \wedge \lambda^l(\sigma)$ and
$k \in \ul\sigma, i \in \ul\tau \wo \ul\sigma$ then
$\d_{ik} a \in \lambda^{p-l-1}(\tau; \sigma) \wedge \lambda^l(\tau)$,
and therefore the differential $d'$ respects the filtration $F^\bu$,
so it is an increasing filtration of $(K^\bu(\sigma, \Lambda^p), d')$
by subcomplexes. Since $\bar \Lambda^\bu(\sigma)$ and
$\Lambda^\bu(\sigma)$ embed into
$\bigoplus\limits_{\tau \supset \sigma} \lambda^\bu(\sigma)$, both
graded vector spaces inherit the filtration.

\begin{lemma}
  \label{parity}
  Let $J' \subset J$ and let $i \notin J'$.  We have
  \begin{itemize}    
  \item[-] $\sgn(i, J' \cup \{i\})\sgn(i, J \wo J') = \sgn(i, J),$ 
  \item[-] $\dfrac{\prod_{j \in J'}(-1)^{|J|}\sgn(j, J)}{\prod_{j \in
        J'} (-1)^{|J|-1}\sgn(j, J \wo \{i\})} = \sgn(i, J' \cup \{i\})$.
  \item \item[-]
    $\dfrac{\prod_{j \in J' }(-1)^{|J|-1}\sgn(j, J \wo
      \{i\})}{\prod_{j \in J' \cup \{i\}}
      (-1)^{|J|}\sgn(j, J)} = \sgn(i, J \wo J')$.
  \end{itemize}
\end{lemma}

\begin{proof}
  The first statement is straigtforward.

  To prove the second statement, notice that
  $$
  \begin{array}{ll}
    (-1)^{|J|}\sgn(j, J) = -(-1)^{|J|-1}\sgn(j, J \wo \{i\}
    & \textrm{ for } j < i, \\
    (-1)^{|J|}\sgn(j, J) = (-1)^{|J|-1}\sgn(j, J \wo \{i\})
    & \textrm{ for } j > i.\\
  \end{array}
  $$
  Therefore, the left hand side expression is equal $(-1)^{l}$, where
  $l$ is the number of elements of $J'$ less than $i$. The statement
  follows.

  The left hand side in the third statement differs from left hand
  side in the second statement by $\sgn(i, J)$. So we can conclude by
  combining second and the first statement.
\end{proof}

\begin{lemma}
  \label{reso}
  The cohomology of the  complex $K^\bu(\sigma, \Lambda^p)$ is
  concentrated in degree 0 for all $p \geq 0$.
\end{lemma}

\begin{proof}

  Clearly
  $\gr^m_F \lambda^p(\tau) \cong \lambda^{p-m+1}(\tau; \sigma) \wedge
  \lambda^{m-1}(\sigma)$. If $|\ul\tau| = |\ul\sigma| + p-m+1$, the
  space $\lambda^{p-m+1}(\tau; \sigma)$ is of dimension 1 and the map
  $\kappa^\tau_\sigma$ establishes an isomorphism between
  $\gr^m_F \lambda^p(\tau)$ and $\lambda^{m-1}(\sigma)$.  Observe that
  if $|\ul\tau| < |\ul\sigma| + p-m+1$ then
  $\gr^m_F \lambda^p(\tau) = 0$, and if
  $|\ul\tau| > |\ul\sigma| + p-m+1$ then any element
  $a_\tau \otimes b_\tau \in F^m \lambda^p(\tau)$ is a sum of elements
  of the form
  $a_\eta \otimes b_\tau|_{\eta} \in F^m \lambda^{p-m+1}(\eta),
  a_\eta|_\tau = a_\tau,$ where $\eta \subset \tau$ is a face such
  that $|\ul\eta| = |\ul\sigma| + p-m+1$. Therefore, the maps
  \begin{multline*}
    \pi_{m,r}: \gr^m_F K^r(\sigma, \Lambda^p) \to \bigoplus_{\tau \supset^{p-r-m+1}
      \sigma} \Hh^{2r}(Y_\tau) \otimes \lambda^{m-1}(\sigma) \cong \Mm^r_{p-m+1} (Y_\sigma, Y_\sigma^\circ) \otimes \lambda^{m-1}(\sigma),\\
    \sum_{\tau \supset \sigma} a_\tau \otimes b_\tau \mapsto \sum_{\tau
      \supset^{p-r-m+1} \sigma} N_\tau a_\tau \otimes
    \kappa^\tau_{\sigma}(b_\tau)\\
  \end{multline*}
  are isomorphisms.  In fact, for a fixed $m$ the maps $\pi_{m,r}$
  define a morphism of complexes, and hence a quasi-isomorphism: for any
  $v=\sum\limits_{\tau \supset \sigma} a_\tau \otimes b_\tau \in
  K^r(\sigma, \Lambda^p)$,
  \begin{multline*}
    \pi_{m,r+1}(d' v) = \sum_{\tau \supset^{p-r-m} \sigma} \sum_{i \in
      \ul\tau \wo \ul\sigma} N_i N_{\6_i \tau} \Gys^\tau_{\6_i
      \tau} (a_\tau) \otimes \kappa^{\6_i \tau}_\sigma ((\d_{ik} b_\tau)|_{\6_i \tau}) = \\
    = \sum_{\tau \supset^{p-m-r} \sigma} \sum_{i \in \ul\tau \wo
      \ul\sigma}N_\tau \sgn(i,\ul\tau \wo \ul\sigma) \Gys^\tau_{\6_i
      \tau} a_\tau \otimes \kappa^\tau_\sigma (b_\tau),
  \end{multline*}
  Since $M_{p-m+1}^\bu(Y_\sigma, Y_\sigma^\circ)$ has cohomology
  concentrated in degree 0, by Lemma~\ref{h0m} then so does
  $\gr^m_F K^r(\sigma, \Lambda^p)$.
\end{proof}

\begin{defn}
  Let
  $$
  K^{q,r}(\Delta_X, \Lambda^p) = \bigoplus_{|\ul\sigma|=q+1} K^r(\sigma,
  \Lambda^p)
  $$
  and for any $\alpha \supset \sigma$ let $\pi^\sigma_\alpha$ be the
  natural projection
  $$
  \pi^\sigma_\alpha: \bigoplus_{\tau \supset \sigma} \Hh^{2r}(Y_\tau)
  \otimes \lambda^{p-r}(\tau) \to \bigoplus_{\beta \supset \alpha}
  \Hh^{2r}(Y_\beta) \otimes \lambda^{p-r}(\beta).
  $$
  Define the horizontal differential $d'': K^{q,r}(\Delta_X, \Lambda^p) \to
  K^{q+1,r}(\Delta_X, \Lambda^p)$ 
  $$
  d''(x) = \sum_{\alpha:\6_i \alpha=\sigma} \sgn(i, \ul\alpha)
  \pi^\sigma_\alpha.
  $$
\end{defn}

The inclusions $\Lambda^\bu(\sigma) \hookrightarrow K^0(\sigma,
\Lambda^\bu)$ extend to the inclusions
$$
\bigoplus_{|\ul\sigma|=q} \Lambda^\bu(\sigma) = C^q(\Delta_X, \Lambda^\bu)
\hookrightarrow K^{q,0}(\Delta_X, \Lambda^\bu)
$$
for each $q \geq 0$ and 
it is clear from the definition of $d''$ that they define a morphism
of complexes. Moreover, by Lemma~\ref{reso}
$$
C^0(\Delta_X, \Lambda^p) = \Ker d'.
$$
It follows that the induced morphism $C^\bu(\Delta_X, \Lambda^\bu) \to
sK^\bu(\Delta_X, \Lambda^\bu)$ is a quasi-isomorphism.

\subsection{Logarithm of monodromy}

Let $\sigma$ be a face of $\Delta_X$,
$\ul\sigma=\set{i_1, \ldots, i_l}, l \geq 2$, and assume that
$j,k \in \ul\sigma$ are top two vertices of $\sigma$ in the
orientation order.  Define
$N_{\sigma}: K^r(\6_k \sigma, \Lambda^p) \to K^r(\sigma,
\Lambda^{p-1})$ to be the map
$$
N_\sigma(v) = \sum_{\tau \supset \sigma} (-1)^{|\ul\tau|}a_\tau \otimes \d_{jk} b_\tau,
$$
where $v$ is a sum of primitive tensors
$\sum_{\tau \supset \6_k\sigma} a_\tau \otimes b_\tau$ (the meaning of
the sign will become apparent in the proof of
Lemma~\ref{n-equiv}). This map is always well-defined when $p=1$ and
is well-defined for $p > 1$ if $\Lambda^\bu$ is regular at
$\sigma$. Define the map
$$
N: K^{q,r}(\Delta_X, \Lambda^\bu) \to K^{q+1,r}(\Delta_X,
  \Lambda^{\bu-1}), \qquad
  N(v) = (-1)^{q+1}\sum_{|\ul\sigma|=q+2} N_\sigma,
$$
which in turn induces the map
$N: sK^\bu(\Delta_X, \Lambda^\bu) \to sK^{\bu+1}(\Delta_X,
\Lambda^{\bu-1})$, when it is well-defined.

\begin{lemma}
  \label{n-is-cobound}
  The map on the cohomology
  $$
  H^\bu(\Delta_X, \Lambda^1) \to H^{\bu+1}(\Delta_X, \Lambda^0)
  $$
  induced by the map $N$ coincides with the coboundary map associated
  with the short exact sequence of sheaves (\ref{seq:lambda'}). If
  $\Lambda^\bu$ is regular at every face $\sigma$ of $\Delta_X$ then
  the statement is also true about the map
  $$
  H^\bu(\Delta_X, \Lambda^{p+1}) \to H^{\bu+1}(\Delta_X, \Lambda^p)
  $$
  induced by $N$ for all $p > 0$.
\end{lemma}

\begin{proof}
  The proof uses an idea similar to the one used in the proof of
  \cite[Proposition~3.5]{jsr-epi}.

  Since the inclusion
  $C^\bu(\Delta_X, \Lambda^\bu) \to sK^\bu(\Delta_X, \Lambda^\bu)$ is
  a quasi-isomorphism, it suffices to check the statement on the
  complex $C^\bu(\Delta_X, \Lambda^\bu)$.

  Take a cocycle
  $a = (a_\sigma) \in C^q(\Delta_X, \Lambda^p)$, where $a_\sigma \in
  \Lambda^p(\sigma)$. If
  $\tilde a = (\tilde a_\sigma) \in C^q(\Delta_X, A^p)$ is some lifting of the
  cocycle $a$ then
  $$
  d \tilde a \in \Ker ( C^{q+1}(\Delta_X, A^p) \to
  C^{q+1}(\Delta_X, \Lambda^p))=\Im (h: C^{q+1}(\Delta, \Lambda^{p-1})
  \to C^{q+1}(\Delta_X, A^p)).
  $$
  We will pick $\tilde a$ in such a way that $d \tilde a = h(N a)$.

  Recall that for any face $\tau$ elements of $\lambda^1(\tau)$ are
  translation-invariant differential forms with rational coefficients
  on
  $\tau \subset T(\tau)\otimes \R \subset H^0(\St(\tau), A^1)^*
  \otimes \R$. The elements of $A^1(\tau)$ are tautologically
  identified with linear functions on $H^0(\St(\tau), A^1)^*$ or with
  affine functions on $e(\tau) \subset H^0(\St(\tau), A^1)^*$. More
  generally, elements of $\bigwedge^p A^1(\tau)$ can be identified
  with translation-invariant differential forms with rational
  coefficients on the linear subspace spanned by
  $e(\tau) \subset H^0(\St(\tau), A^1)^* \otimes \R$. If $p > 1$ we
  assume from now on that $\Lambda^\bu$ regular at any face of
  $\Delta_X$: 
  $\Lambda^p(\tau) = \bigwedge^p \Lambda^1(\tau)$. Recall that
  $1_\tau$ is the function that is constantly 1 on $\tau$, we will use
  the same notation for the corresponding differential form. The
  natural inclusion
  $$
  \Lambda^p(\tau) \hookrightarrow A^p(\tau)
  $$
  is induced by wedging with $1_\tau$. Note that this is well-defined,
  since $\Lambda^1(\tau)$ is the quotient of $A^1(\tau)$ be the
  subspace spanned by $1_\tau$.

  For each element $a_\sigma \in \Lambda^p(\sigma)$ pick the unique
  lifting $\tilde{a}_\sigma \in A^p(\sigma)$ such that
  $\d_j a_\sigma = 0$, where $j \in \ul\sigma$ is maximal with respect
  to the orientation ordering. Since $da=0$, we have that
  $d\tilde{a} \in C^{q+1}(\Delta_X, \Lambda^{p-1}) \subset
  C^{q+1}(\Delta_X, A^p)$. Pick a face $\tau \subset \Delta_X$. Let
  $j < k$ be two topmost vertices of $\tau$ in the orientation
  ordering. We have
  $$
  \d_k \tilde{a}_{\6_i \tau} = 0, \textrm{ for all } i \in \ul\tau, i
  \neq k, \qquad \d_j \tilde{a}_{\sigma} = 0,
  $$
  where $\sigma = \6_k \tau$.

  The elements of the form $1_\tau \wedge b \in A^p(\tau)$ are
  characterized by the property
  $$
  \d_i (1_\tau \wedge b) = \d_j (1_\tau \wedge b),
  $$
  for any vertices $i,j \in \ul\tau$, moreover, for any vertex $i
  \in ul\tau$,
  $$
  \d_i (1_\tau \wedge b) = b|_{\la \tau \ra}.
  $$
  In particular,  
  $$
  \d_k (d\tilde{a})_\tau = \sum_{i \in \tau} \sgn(i, \ul\tau) \d_k
  a_{\6_i \tau} = (-1)^{|\ul\tau|} \d_k \tilde{a_\sigma} =
  (-1)^{|\ul\tau|} (\d_k - \d_j) \tilde{a}_\sigma.
  $$
  But since $\tilde{a}_\sigma|_{\la e(\sigma) \ra} = a_\sigma$, we
  have
  $$
  (d\tilde{a})_\tau = (-1)^{|\ul\tau|} 1_\tau \wedge \d_{jk}
  a_\sigma|_{\tau} = (Na)_\tau,
  $$
  and we conclude.  
\end{proof}

\section{Cohomology of $\Lambda^p$}
\label{s:cohomology}

\subsection{Proof of Theorems~A and A'}

Let $\alpha, \sigma$ be faces such that $\6_i \alpha = \sigma$. Denote
the unique vector in $\lambda^1(\alpha; \sigma)$ such that
$\d_{ik} dz_{ik,\sigma} = 1$ for some (equivalently, any)
$k \in \ul\sigma$ as $dz_{i,\sigma}$. Clearly,
$$
\sum_{i \in \ul\sigma} dz_{i,\6_i \sigma} = 0.
$$

Let $G^\bu$ be the filtration on the total complex
$sK^\bu(\Delta_X, \Lambda^p))$ induced by the filtration
$$
G^m K^r(\sigma, \Lambda^p) = F^{|\ul\sigma|+r-m} K^r(\sigma, \Lambda^p).
$$

\begin{lemma}
  \label{e0-of-ss}
  The zeroth page of the
  spectral sequence associated to the filtration $G$ on
  $sK^\bu(\Delta_X, \Lambda^p))$ has the following form
  $$
  \EG{p}^{i,j}_0(Y) = \bigoplus_{r \geq 0} \bigoplus_{|\sigma|=i+j-r+1}
  \bigoplus_{\tau \supset^{p-r-j} \sigma} \Hh^{2r}(Y_\tau) \otimes
  \lambda^j(\sigma)
  $$  
  with the differential $d_0: \EG{p}^{i,j}_0 \to \EG{p}^{i,j+1}_0$
  defined as follows on primitive tensors
  $a_\tau \otimes b_\tau \in \Hh^{2r}(Y_\tau) \otimes
  \lambda^j(\sigma)$:
  $$
  d_0(a_\tau \otimes b_\tau) = \sum_{l \in \ul\tau \wo \ul\sigma} a_\tau \otimes
  \sgn(l, \ul\tau)(dz_{l,\sigma} \wedge b_\tau)
  $$
  where $k$ is the maximal vertex in $\ul\sigma$ with respect to the
  orientation ordering.
\end{lemma}

\begin{proof}
  The first statement follows from the fact that
  \begin{multline*}
  \EG{p}^{i,j}_0(Y) = \gr^i_G sK^{i+j}(\Delta_X, \Lambda^p) =  \bigoplus_{r \geq 0} \bigoplus_{r+|\ul\sigma|-1=i+j} \gr^{i}_G
  K^r(\sigma, \Lambda^p) = \\
  = \bigoplus_{r \geq 0} \bigoplus_{|\ul\sigma|=i+j-r+1} \gr^{j+1}_F
  K^r(\sigma, \Lambda^p).
  \end{multline*}

  The contribution of the vertical differential
  $d': K^r(\sigma, \Lambda^p) \to K^{r+1}(\sigma, \Lambda^p)$ to
  $$
  d_0: \bigoplus_{r \geq 0} \bigoplus_{|\ul\sigma|=i+j-r+1}
  \gr^{j+1}_F K^r(\sigma, \Lambda^p) \to 
  \bigoplus_{r \geq 0} \bigoplus_{|\ul\sigma|=i+j-r+2}
  \gr^{j+2}_F K^r(\sigma, \Lambda^p)
  $$
  is zero because $d'$ preserves the filtration $F$.

  We identify $\gr^{j+1}_F K^r(\sigma, \Lambda^p)$ with
  $\bigoplus_{\tau \supset^{p-r-j} \sigma} \Hh^{2r}(Y_\tau) \otimes
  \lambda^j(\sigma)$ via the isomorphism
  $\id \otimes \kappa^\tau_{\sigma,k}$ where $k$ is the vertex in
  $\ul\sigma$ that is maximal with respect to the orientation
  ordering. Let us compute the contribution of the horizontal
  differential
  $d'': K^{q,r}(\Delta_X, \Lambda^p) \to K^{q+1,r}(\Delta_X,
  \Lambda^p)$.

  Consider a sum of primitive tensors
  $v = \sum_{\tau \supset \sigma} a_\tau \otimes b_\tau \in
  K^r(\sigma, \Lambda^p)$. Take some $l \in \ul\tau \wo \ul\sigma$ and
  let $\6_l \alpha = \sigma$. Then the images of $v$ and $d'' v$ in
  $\gr^{j+1}_F K^r(\sigma, \Lambda^p)$ and
  $\gr^{j+2}_F K^r(\alpha, \Lambda^p)$ are, respectively,
  $$
  \sum_{\tau \supset^{p-r-j} \sigma} a_\tau \otimes
  \kappa^\tau_{\sigma,k}(b_\tau), \qquad \sgn(l,\ul\alpha) \sum_{\tau
    \supset^{p-r-j-1} \alpha} a_\tau \otimes
  \kappa^\tau_{\alpha,k}(b_\tau),
  $$
  with
  $$
  \kappa^\tau_{\alpha,k}(b_\tau) = \sgn(l,\ul\tau \wo \ul\sigma)
  dz_{lk,\sigma} \wedge \kappa^\tau_{\sigma,k}(b_\tau).
  $$
  Since
  $$
  \sgn(l,\ul\alpha)\sgn(l,\ul\tau \wo \ul \sigma) = \sgn(l, \ul\tau)
  $$
  by Lemma~\ref{parity}, the statement of the lemma follows.
\end{proof}

For any face $\tau$ and a number
$n \leq |\ul\sigma|$ consider the complex $S^\bu_n(\tau)$
$$
0 \to \bigoplus_{\sigma \subset^n \tau} \lambda^0(\sigma)
\to \bigoplus_{\sigma \subset^{n-1} \tau} \lambda^1(\sigma)
\to \ldots \to \lambda^n(\tau) \to 0
$$
with the differential sending $a \in \lambda^p(\sigma)$ to
$\sum_{\6_l \alpha=\sigma} \sgn(l, \ul\sigma) dz_{l,\sigma} \wedge a$.
  
\begin{lemma}
  \label{acyclic}
  For any $\tau$ and any $n < |\ul\tau|$
  $$
  H^i(S^\bu_n(\tau)) = \left\{\begin{array}{ll}
                              \Q & i = 0,\\
                              0, & i > 0.
                            \end{array}\right.
  $$  
  Moreover, the unique up to scalar $0$-th cohomology class is
  represented by a cocycle $a = (a_\sigma)_{\sigma \subset^n \tau}$ with
  $$
  a_\sigma = \prod_{i \in \ul\sigma} (-1)^{|\ul\tau|}\sgn(i, \ul\tau)
  $$
\end{lemma}

\begin{proof}
  The proof is by induction on $n$. Let's check the base of induction:
  clearly, $S^0_0(\tau) = \lambda^0(\tau)$ and $S^i_0(\tau)=0$ for
  $i > 0$. For the induction step, let $k$ be some vertex in
  $\ul\sigma$ and consider the exact sequence
  $$
  0 \to S^\bu_n(\tau,\6_k \tau) \to S^\bu_n(\tau) \to
  S^\bu_{n-1}(\6_k \tau) \to 0,
  $$
  where $S^\bu_n(\tau,\6_k \tau)$ is the subcomplex with the terms
  $$
  S^i_n(\tau,\6_k \tau) = \bigoplus\limits_{\tau \supset^{n-i}
    \sigma, \sigma \not\subset \6_k\sigma} \lambda^i(\sigma).
  $$
  Since $S^\bu_{n-1}(\6_k \tau)$ only has cohomology $\Q$ in degree
  0, and vanishing cohomology in all other degrees by induction
  hypothesis, in order to prove that $S^\bu_{n}(\tau)$ has the same
  property it would suffice to show that the complex
  $S^\bu_n(\tau, \6_k \tau)$ is acyclic.

  Notice that for any faces $\alpha, \beta, \eta \subseteq \sigma$
  such that $\ul\alpha \cup \ul\beta = \ul\eta$ and
  $\ul\alpha \cap \ul\beta = \{k\}$, and any numbers
  $m_1, m_2, m_1 \leq |\ul\alpha|, m_2 \leq |\ul\beta|$,
  $$
  S^\bu_{m_1+m_2}(\eta, \6_k \sigma) \cong S^\bu_{m_1}(\alpha, \6_k \alpha)
  \otimes S^\bu_{m_2}(\beta, \6_k \beta)
  $$
  (recall that the tensor product in the right hand side is the total
  complex of the double complex
  $(S^\bu_{m_1}(\alpha, \6_k \alpha) \otimes S^\bu_{m_2}(\beta, \6_k
  \beta), d_\alpha, d_\beta)$). In particular, denoting
  $\alpha_l \subset \sigma$ the 1-faces with $\ul\alpha_l=\{l,k\}$, we
  have
  $$
  S^\bu_n(\sigma, \6_k \sigma) = \bigoplus_{\tau \subset \6_k \sigma,
    |\ul\tau| = n} \bigotimes_{l \in \ul{\6_k \sigma}}
  S^\bu_{1_{\ul\tau}(l)}(\alpha_l, l),
  $$
  where $1_{\ul\tau}: \ul{\6_k \sigma} \to \{0,1\}$ is the indicator
  function of $\ul\tau$. The complex
  $S^\bu_1(\alpha_l, l)$ has the form $\Q \to \Q$ with
  non-trivial differential, and so is acyclic, and
  $S^\bu_0(\alpha_i, \6_k \alpha)$ has unique non-trivial term $\Q$ in
  degree 0.  Therefore, since there is at least one $l \in \ul\tau$,
  the tensor product is acyclic, and so is
  $S^\bu_n(\sigma, \6_k\sigma)$.

  For the last statement of the Lemma, direct computation gives
  \begin{multline*}
  (da)_\sigma = \sum_{i \in \ul\sigma} \left(\prod_{j \in \ul\6_i
      \sigma} (-1)^{|\ul\tau|} \sgn(j, \ul\tau)\right) \sgn(i,
  \ul\tau) dz_{i, \6_i\sigma} = \\ = (-1)^{|\ul\tau|} \prod_{j \in
    \ul\sigma} (-1)^{|\ul\tau|} \sgn(j, \ul\tau) (\sum_{i \in \ul\sigma} dz_{i,
    \6_i\sigma}) = 0.
  \end{multline*}
  So $a$ is indeed a cocycle representing a non-trival class in
  $H^0(S^\bu_n(\sigma))$.
\end{proof}

\begin{lemma}
  \label{e1-of-ss}
  The first page of the spectral sequence associated to
  the filtration $G$ on $sK^\bu(\Delta_X, \Lambda^p)$ consists of a single
  row:
  $$
  (\EG{p}_1^{\bu,0}(Y), d_1) = (sD^\bu_p(Y), d).
  $$
\end{lemma}

\begin{proof}
  By Lemma~\ref{e0-of-ss},
  $$
  \EG{p}_0^{i,j}(Y) = \bigoplus_{r = 0}^{\min\{i,p\}} \bigoplus_{|\tau|=p -2r + i
    +1} \Hh^{2r}(Y_\tau)\otimes S^j_{p-r}(\tau),
  $$
  so by Lemma~\ref{acyclic}, $\EG{p}_1^{i,j}(Y)=0$ unless $j=0$, and
  $$
  \EG{p}_1^{i,0}(Y) = \bigoplus_{r = 0}^{\min\{i,p\}}
  \bigoplus_{|\tau|=p-2r+i+1}  \Hh^{2r}(Y_\sigma) = \bigoplus_{r = 0}^{\min\{i,p\}}
    \Hh^{2r}(Y^{(p-2r+i+1)}) = sD^i_p(Y).
  $$

  Fix $p$ and denote for brevity
  $L^\bu = sK^\bu(\Delta_X, \Lambda^p))$.  The differential
  $d_1: H^{i+j}(\gr^i_G L^\bu) \to H^{i+j+1}(\gr^{i+1}_G L^\bu)$ is the
  boundary morphism on the cohomology associated to the short exact
  sequence
  $$
  0 \to \gr^{i+1}_G L^\bu \to G^i L^\bu/G^{i+2} L^\bu \to \gr^i_G L^\bu \to 0
  $$
  If $v$ is an element of $G^{i+j} L$ that represents a cohomology
  class $[v] \in H^{i+j}(\gr^i_G L^\bu)$ then
  $dv \in G^{i+j+1}L^{i+1}$ and $dv$ represents $d_1([v])$ in
  $H^{i+j+1}(\gr^{i+1}_G L^\bu)$.

  Take
  $$
  [v]=\sum_{r \geq 0} \sum_{|\tau|=p-2r+i+1} a_\tau \in \bigoplus_{r \geq 0}
  \bigoplus_{|\tau|=p-2r+i+1} \Hh^{2r}(Y_\tau) \cong \Hh^i( \gr^i_G L^\bu),
  $$
  We identify 
  $$
  H^i( \gr^i_G L^\bu) = \bigoplus_{|\tau|=p-2r+i+1} \Hh^{2r}(Y_\tau)
  \otimes S^0_{p-r}(\tau) \textrm{ and }  \bigoplus_{|\tau|=p-2r+i+1} \Hh^{2r}(Y_\tau)
  $$
  via the isomorphism that sends the cocycle $a$ from the statement of
  Lemma~\ref{acyclic}, to 1.  Denote $\vol^\tau_\sigma$ the unique
  element of $\lambda^{p-r}(\tau;\sigma)$ such that
  $\kappa^\tau_\sigma(\vol^\tau_\sigma) = 1_\sigma$.  Then the
  cohomology class of $[v]$ can be represented by a sum
  $$
  v = \sum_{r \geq 0} \sum_{|\tau|= p-2r+i+1 \atop \sigma
    \subset^{p-r} \tau } v_{\sigma,\tau} \in G^i L, \qquad v_{\sigma,\tau} \in
  K^r(\sigma, \Lambda^p), 
  $$
  where
  $$
  (v_{\sigma,\tau})_\tau = \prod_{l \in \ul\sigma}
  (-1)^{|\tau|}\sgn(l, \ul\tau) \cdot a_\tau \otimes \vol^\tau_\sigma.  
  $$
  Then
  $$
  d'(v_{\sigma,\tau}) = \prod_{l \in \ul\sigma} (-1)^{|\tau|}\sgn(l,
  \ul\tau)  \cdot \sum_{l \in \ul\tau \wo
    \ul\sigma} N_l \Gys^\tau_{\6_l \tau} (a_\tau) \otimes
  \d_{lk}\vol^\tau_\sigma|_{\6_l\tau} + \sum_{l \in \ul\sigma} N_i
  \res_\tau \Gys^\tau_{\6_l \tau} (a_\tau) \otimes
  \d_{lk}\vol^\tau_{\sigma}.
  $$
  Since $\vol^\tau_\sigma \in \lambda^{p-r}(\tau;\sigma)$ we have that
  $\d_{lk}\vol^\tau_\sigma= 0$ for all $l \in \ul\sigma$, and the
  summands in the second sum vanish. And since for
  $l \in \ul\tau \wo \ul\sigma$ we have
  $$
  \kappa^{\6_l \tau}_{\sigma,k} (\d_{lk}\vol^\tau_\sigma) =
  \sgn(l,\ul\tau \wo \ul\sigma)
  \kappa^\tau_{\sigma,k}(\vol^\tau_\sigma)
  $$
  (as in the proof of Lemma~\ref{reso}), we have that the image of
  $d'(a_\tau \otimes \vol^\tau_\sigma)$ in $\gr^{i+1}_G L^{i+1}$ is
  $$
  \sum_{l \in \ul\tau \wo \ul\sigma} N_l \sgn(l,\ul\tau \wo \ul\sigma)
  \dfrac{\prod_{l' \in  \ul\sigma}
    (-1)^{|\ul\tau|}\sgn(l',\ul\tau)}{\prod_{l' \in  \ul\6_l\sigma} (-1)^{|\ul\tau|-1}\sgn(l',\ul\tau)}
  \Gys^\tau_{\6_i \tau} (a_\tau) = \sum_{l \in \ul\tau \wo \ul\sigma} N_l \sgn(l,\ul\tau)
  \Gys^\tau_{\6_i \tau} (a_\tau)
  $$
  by Lemma~\ref{parity}.
  
  Similarly, since $\vol^\beta_\alpha|_\sigma = \vol^\tau_\sigma$ for
  any pair of faces $\beta \supset \alpha$ such that $\6_i \beta=\tau$
  and $\6_i \alpha = \sigma$ for some $i \notin \ul\tau$,
  $$
  d''(v)_\beta = \sgn(l,\ul\alpha) a_\beta \otimes
  \vol^\beta_\alpha = \sgn(l,\ul\alpha) a_\tau|_{Y_\beta} \otimes
  \vol^\beta_\alpha.
  $$
  Therefore, the image of
  $d''(v)_\beta$ in $H^{i+1}(\gr^{i+1} L^\bu)$
  is then
  $$
  \sgn(l,\alpha) \dfrac{\prod_{l \in
      \ul\sigma} (-1)^{|\tau|} \sgn(l, \ul\tau)}{\prod_{l \in
      \ul\alpha} (-1)^{|\beta|} \sgn(l,\ul\beta)} a_\tau|_{Y_\alpha} =
   \sgn(l,\ul\beta) a_\tau|_{Y_\alpha}
  $$
  by Lemma~\ref{parity}.
  
  We observe that $d_1=d'+d''$ coincides with the differential of the
  complex $sD^\bu_p$.
\end{proof}

Coboundary morphism $N$ associated to the exact sequence
(\ref{seq:lambda'}) is analogous to the eigenwave morphism in tropical
geometry, introduced by Mikhalkin and Zharkov in \cite{mz14}. See
Proposition~3.5 \cite{jrs17} for the comparison between a coboundary
morphism of a sequence analogous to \ref{seq:lambda'} in tropical
geometry. We will now show that under the isomorphism from
Lemma~\ref{e1-of-ss} the coboundary morphism corresponds to the
morphism $N$ on the weight spectral sequence for the limit mixed Hodge
structure induced by the logarithm of the monodromy morphism.

\begin{lemma}
  \label{n-equiv}
  The map
  $$
  N: \EG{p}_1^{i,0}(Y) \to \EG{p-1}_1^{i+1,0}(Y)
  $$
  induced by the morphism $N: sK^r(\Delta_X, \Lambda^p)
  \to sK^{r+1}(\Delta_X, \Lambda^{p-1})$ is identity.
\end{lemma}

\begin{proof}
  As in the proof of Lemma~\ref{e1-of-ss} take
  $$
  [v] = \sum_{r \geq 0} \sum_{|\tau|=p-2r+i+1} a_\tau \in H^i(\gr^i_G
  sK^\bu(\Delta_X, \Lambda^p)) \cong \EG{p}_1^{i,0}(Y)
  $$
  represented by an element
  $$
  v = \sum_{r\geq 0}\sum_{|\tau|= p-2r+i+1 \atop \sigma \subset^{p-r}
    \tau} v_{\sigma,\tau}, \ \  v_{\sigma,\tau} \in K^r(\sigma,
  \Lambda^p), \ \  (v_{\sigma,\tau})_\tau = \left( \prod_{l \in
      \ul\sigma} (-1)^{|\ul\tau|} \sgn(l, \ul\tau) \right) a_\tau \otimes
  \vol^\tau_\sigma,
  $$
  Take a face $\alpha \subset \tau$,  let $j,k$ be the two topmost
  vertices in $|\ul\alpha|$, and denote $\sigma = \6_k \alpha$, then
  since
  $\d_{j k} \vol^\tau_\sigma = \sgn(j, \ul\tau \wo \ul\sigma)
  \vol^\tau_\alpha$, we have
  $$
  N_\alpha(v_{\sigma,\tau}) = (-1)^{|\ul\tau|+|\ul\sigma|+1}\prod_{l \in
      \ul\sigma} (-1)^{|\ul\tau|} \sgn(l, \ul\tau) \sgn(j, \ul\tau \wo
    \ul\sigma) \cdot a_\alpha \otimes \vol^\tau_\alpha.
  $$
  The image of $N_\alpha(v_{\sigma,\tau})$ in
  $H^{i+1}(\gr^{i+1}_G sK^\bu(\Delta_X, \Lambda^{p-1})$ is thus by Lemma~\ref{parity} 
  \begin{multline*}
  (-1)^{|\ul\tau|+|\ul\sigma|+1}\dfrac{\prod_{l \in \ul\sigma} (-1)^{|\ul\tau|} \sgn(l, \ul\tau)
  }{\prod_{l \in \ul\alpha} (-1)^{|\ul\tau|} \sgn(l, \ul\tau)} \sgn(j,
  \ul\tau \wo \ul\sigma) \cdot a_\alpha =\\
  =(-1)^{|\ul\sigma|+1}\sgn(j, \ul\tau) \sgn(j, \ul\tau \wo
  \ul\sigma) a_\alpha 
  = (-1)^{|\ul\sigma|+1}\sgn(j, \ul\alpha) a_\alpha,
  \end{multline*}
  which is equal to just $a_\alpha$, since $j$ is maximal in
  $\ul\alpha$.
\end{proof}

\begin{thm}
  \label{thm-a}
  Let $f: X \to S$ be a unipotent degeneration and assume that there
  exists a cohomologically K\"ahler class in $H^2(Y)$. Then for all
  $p,q \geq 0$, there exists a map 
  $$
  H^q(\Delta, \Lambda^p) \to \gr^{2p}_W H^{p+q}(X_\infty).
  $$
  Assume that $\Lambda^\bu$ is regular at every face of
  $\Delta_X$. Then the map commutes with the logarithm of monodromy
  morphism $N$ on the right, and the coboundary morphism of the short
  exact sequence (\ref{seq:lambda'}) on the left.
\end{thm}

\begin{proof}
  The proof is a conjunction of Lemmas~\ref{h0m}, \ref{reso},
  \ref{e0-of-ss}, \ref{acyclic}, \ref{e1-of-ss}, \ref{n-equiv} and
  Fact~\ref{hl-module}.
\end{proof}

\begin{thm}
  \label{hodge}
  Assume that there exists a combinatorial Lefschetz class
  $\omega \in H^2(Y)$, then the morphism
  $$
  H^\bu(\Delta, \Lambda^\bu) \to \gr^W_{2\bu} H^\bu(X_\infty)
  $$
  constructed in Theorem~A is injective and
  $\dim H^{n-q}(\Delta_X, \Lambda^{n-p}) = \dim H^q(\Delta_X,
  \Lambda^p)$. If $\Lambda^\bu$ is regular at any face
  $\sigma \subset \Delta_X$ and the morphism $N$ is well-defined then
  $$
  N^{p-q}: H^q(\Delta_X, \Lambda^p) \to H^p(\Delta_X, \Lambda^q)
  $$
  is an isomorphism.
\end{thm}

\begin{proof}
  By Fact~\ref{hl-module} $(K^{\bu,\bu} \otimes \R, N, L_\omega)$ is a
  Hodge-Lefschetz module, since $\omega|_{Y_\sigma}$ is a Lefschetz
  class for all strata $Y_\sigma$. Since $\omega$ is also
  combinatorial, the operator $L_\omega$ preserves the subcomplex
  $\Kk^{\bu,\bu}$, so it is a Hodge-Lefschetz submodule of
  $(K^{\bu,\bu}, N, L_\omega)$.  By Fact~\ref{harmonic}, the inclusion
  $\Ker \Box \cap \Kk^{\bu,\bu} \hookrightarrow \Ker \Box$ induces the
  inclusion of cohomology of complexes
  $s\Dd^\bu_{2p} \hookrightarrow sD^\bu_{2p}$ for all $p$. The last
  statement of the theorem follows from the definition of
  combinatorial Lefschetz classes and Fact~\ref{harmonic}.
\end{proof}

\begin{prop}
  If $X \to S$ is a degeneration of curves such that the central fibre
  has at least one double point, then any cohomologically K\"ahler
  class is combinatorial Lefschetz.
\end{prop}

\begin{proof}
  Immediate, since $\Hh^2(Y^{(1)}) = H^2(Y^{(1)})$.
\end{proof}

If $X' \to S^*$ is a degeneration of polarized Abelian varieties over
a punctured disc with multiplicative reduction, then by
\cite{kunnemann98}, there exists a smooth projective model
$f: X \to S$ with a central fibre such that its irreducible components
are toric varieties.

\begin{prop}
  If $X \to S$ is a K\"unemann-Mumford degeneration with toric central
  fibre then
  $$
  H^i(X_\infty) \cong \bigoplus_{p+q=i} H^q(\Delta_X, \Lambda^p)
  $$
  and
  $$
  N^{q-p}: H^q(\Delta_X, \Lambda^p) \to H^q(\Delta_X, \Lambda^p)
  $$
  is an isomorphism.  
\end{prop}

\begin{proof}
  The irreducible components of the central fibre are toric and
  intersect along the irreducible components of the toric boundary,
  therefore, $\Hh^i(Y^{(j)}) \cong H^i(Y^{(j)})$ and
  $\Dd^{\bu}(Y) \cong D^\bu(Y)$, which implies the first
  statement. Any cohomologically K\"ahler class satisfies the
  conditions of Theorem~\ref{hodge}, which implies the second
  statement.
\end{proof}

The morphism constructed in Theorem~A can be more concretely described
with the help of the notion of PL metrized virtual line bundle
(intruduced in \cite{kt02}, see also \cite{tony}).

\begin{defn}
  A \emph{virtual line bundle} is a torsor under the sheaf
  $\Lambda^1$. If $L$ is a virtual line bundle, a \emph{piece-wise
    linear metrization} (or \emph{PL metrization}) of $\omega$ is a
  section of $h \in H^0(L \otimes \bar \Lambda^1)$. The
  \emph{curvature of $h$ at $\sigma$} is the following element of
  $H^2(Y_\sigma)$
  $$
  c_1(L,h)_\sigma = \sum_{i \in \Stbar(\sigma)} \tilde{h}_\sigma(i)
  N_i c_1(\Oo(Y_i))|_{Y_\sigma}
  $$
  where $\tilde{h}_\sigma$ is a lifting of $h$ to a section an element of
  $\bar A^1(\sigma)$. We denote
  $c_1(L,h) = \sum_{\sigma} c_1(L,h)_\sigma \in \bigoplus_\sigma
  H^2(Y_\sigma)$.
\end{defn}

Clearly, if $(L, h)$ and $(L',h')$ are two PL metrized virtual line
bundles then
$$
c_1(L \otimes L', h + h') = c_1(L, h) + c_1(L, h').
$$

\begin{prop}
  Let $a \in C^1(\Delta, \Lambda^1)$ and let $L$ be the corresponding
  virtual line bundle.

  If $Na = 0 \in H^2(\Delta_X, \Lambda^0)$ then $L$ admits a PL
  metrization and for any PL metrization $h$ the class
  $\sp(c_1(L,h)) \in H^2(X_\infty)$ is equal to the image of $a$ in
  $\gr^W_2H^2(X_\infty)$.
\end{prop}

\begin{proof}
  Consider the exact sequence
  $$
  \ldots H^1(\Delta_X, A^1) \to H^1(\Delta_X, \Lambda^1) \xrightarrow{N}
  H^2(\Delta_X, \Lambda^0) \to \ldots
  $$
  Since $a \in \Ker N$, it can be lifted to a cocycle
  $\tilde{a} \in C^1(\Delta_X, A^1) \hookrightarrow C^1(\Delta_X, \bar
  A^1)$. Since the sheaf $\bar A^1$ is flabby,
  $$
  0 \to C^0(\Delta_X, \bar A^1) \to C^1(\Delta_X, \bar A^1) \to
  C^2(\Delta_X, \bar A^1) \to \ldots
  $$
  is exact in degrees $> 0$. In particular, there exists a cochain
  $\tilde{b} \in C^0(\Delta_X, \bar A^1)$ such that
  $d\tilde{b}=\tilde{a}$. The diagram of simplicial cochain  complexes
  arising from the morphism $\bar A^1 \to \bar \Lambda^1$
  $$
  \xymatrix{
     0 \ar[r] & C^0(\Delta_X, \bar A^1) \ar[d] \ar[r] & C^1(\Delta_X,
     \bar A^1) \ar[r] \ar[d] & \ldots \\
     0 \ar[r] & C^0(\Delta_X, \bar \Lambda^1) \ar[r] & C^1(\Delta_X,
     \bar \Lambda^1) \ar[r] & \ldots \\
    }
  $$
  is commutative, and therefore there exists a cochain
  $b \in C^0(\Delta_X, \bar \Lambda^1)$ such that $db=a$. Now observe
  that
  $(C^\bu(\Delta_X, \bar \Lambda^1),d) \cong (K^{\bu,0}(\Delta_X,
  \Lambda^1), d'')$ and that the data of a cochain $b$ defines a PL
  metric $h$ on the virtual line bundle $L$ corresponding to $a$, with
  $c_1(L,h)=d''b$. The cocyle
  $(0, -d''b) \in K^{1,0}(\Delta_X, \Lambda^1) \oplus
  K^{0,1}(\Delta_X, \Lambda^1)$ is then cohomologous to the image of
  $a$ in $sK^1(\Delta_X, \Lambda^1)$. It follows that the image of $a$
  in $\gr^W_2 H^2(X_\infty)$ is $(0, -c_1(L,h))=\sp(-c_1(L,h))$.

  If $h'$ is another metrization of $L$ then there exists a cochain
  $b'$ such that $d''b=a$ and $d'b' = c_1(L,h')$, $a$ is
  cohomologous to $(0, -c_1(L,h')) = \sp(-c_1(L,h))$.
\end{proof}

\subsection{Superforms on dual intersection complexes}

Let $O \subset e(\St(\sigma))$ be an open set that intersects
$e(\sigma)$ non-trivially. Call two germs $\alpha$ and $\beta$ of
$(p,q)$-superforms defined in a neighbourhood of $O$ in $T(\sigma)$
equivalent if 
$$
\alpha|_{e(\tau)} = \beta|_{e(\tau)}
$$
for any face $\tau \supset \sigma$.  For any open set
$U \subset \St(\sigma)$ that intersects $\sigma$ non-trivially define
the set of superforms on $U$ to be the set of equivalence classes of
germs of superforms in a neighbourhood of $e(U) \subset
T(\sigma)$. It is easy to check that differentials $d'$ and $d''$ map
equivalence classes of germs of superforms to equivalence classes,
same is true about the map $N: A^{p,q} \to A^{p-1,q+1}$.

If $\tau \supset \sigma$ then the inclusion
$\St(\tau) \subset \St(\sigma)$ induces a natural map
$r: T(\tau) \to T(\sigma)$. Clearly, $e(\St(\tau))$ is mapped to
$e(\St(\sigma))$ under this map. We call the pullback of a superform
$\alpha$ defined on an open set $U \subset \St(\sigma)$ along this map
its \emph{restriction} to $r(U) \cap \St(\tau)$.

\begin{defn}[Superforms on $\Delta_X$]
  \label{sf-deltax}
  For an open set $U \subset \Delta_X$, let $\Sigma_U$ be the
  collection of faces $\sigma$ such that
  $U \cap \mathring{\sigma} \neq \emptyset$.  A \emph{$(p,q)$-form
    $\alpha$} on an open subset $U \subset \Delta_X$ is a collection
  $(\alpha_\sigma)_{\sigma \in \Sigma_U}$, where $\alpha_\sigma$ is a
  germ of a $(p,q)$-form on a neighbourhood of
  $U \cap \mathring{\sigma}$ in $\St(\sigma)$, such that whenever
  $\sigma \subset \tau$, the restriction of $\alpha_\sigma$ to
  $U \cap \St(\tau)$ is $\alpha_\tau$. We denote the sheaves of
  $(p,q)$-forms on $\Delta_X$ as $\Aa^{p,q}_X$.
\end{defn}

\begin{lemma}  
  The definition of superforms on $\Delta_X$ is invariant under
  subdivisions induced by blow-ups of the strata of $X_0$.
\end{lemma}

\begin{proof}  

  Let $\sigma$ be a face of $\Delta_X$ and let $Y_\sigma$ be the
  stratum being blown up, giving rise to a new degeneration
  $X' \to X$. For any face $\tau \subset \Delta_{X'}$ such that
  $\tau \subset \sigma$ we have by Proposition~\ref{subdiv}
  $H^0(\St(\tau), A^1) \cong H^0(\St(\sigma), A^1)$ and therefore
  $T(\tau) \cong T(\sigma)$. But the definition of $\Lambda^1$ is
  invariant under subdivisions by Proposition~\ref{subdiv}.  
\end{proof}

If $\alpha$ is a local section of $\Aa^{p,q}$ on a subset
$U \subset \St(\sigma)$ then superforms $d'f$ and $d''f$ are germs of
superforms in the neighbourhood of
$e(U) \subset e(\St(\sigma)) \subset T(\sigma)$ and therefore define
local sections of sheaves $\Aa^{p+1,q}$ and $\Aa^{p,q+1}$. This
defines differentials $d',d''$ on all sheaves $\Aa^{p,q}$; maps $J$
and $N$ can be defined similarly.

\begin{prop}
  \label{resolutions}
  \, \noindent
  \begin{enumerate}
  \item for all $p \geq 0$
    $\Ker \{ d'': \Aa^{p,0}_X \to \Aa^{p,1}_X \} \cong \bigwedge \Lambda^1_X
    \otimes \R$;
  \item\label{poinc}
    $\Im \{ d'': \Aa^{p,q} \to \Aa^{p,q+1} \} = \Ker \{ d'':
    \Aa^{p,q+1} \to \Aa^{p,q+2} \}$ for all $p \geq 0$;
  \end{enumerate}

  In particular, if $\Lambda^\bu$ is regular at every face then the
  complex $(H^0(\Delta_X, \Aa^{p,\bu}), d'')$ computes the cohomology
  of the sheaves $\Lambda^p$.
\end{prop}

\begin{proof}
  The first statement is straightforward as soon as one observes that
  $d''$-closed $(p,0)$-superforms on $\St(\sigma)$ for some face
  $\sigma$ can be identified with translation-invariant $p$-forms on
  $T(\sigma)$, which in turn can be identified with the sections
  $H^0(\St(\sigma), \Lambda^p)$.

  Further, it follows from Definition~\ref{sf-deltax} that the second
  statement reduces to the corresponding statement about superforms on
  a vector space. Therefore it is true by \cite[Lemma~1.10]{lag12} or
  \cite[Theorem~2.16]{jell16}.
\end{proof}

\begin{prop}
  \label{q-iso}
  Assume that $\Lambda^\bu$ is regular at all faces
  $\sigma \subset \Delta_X$. Then for any $p \geq 0$ there exists a
  quasi-isomorphism from the complex
  $(H^0(\Delta_X,\Aa^{p,\bu}_X),d'')$ to the complex of singular
  chains $C^\bu(\Delta_X, \Lambda^p \otimes \R)$. If $\Lambda^\bu$ is
  not regular at all faces, then such quasi-isomorphism still exists
  for $p=1$.
\end{prop}

\begin{proof}
  Let $\omega$ be a $p,q$-superform defined in a neighbourhood of a
  face $\sigma \subset \Delta_X$. By definition, $\omega$ is in a
  neighbourhood of $\sigma$ by a germ $\tilde{\omega}_\sigma$ of a
  $p,q$-form in a neighbourhood of
  $e(\Stbar(\sigma)) \subset T_\sigma\Delta_X$. The latter can also be
  regarded as a $\Lambda^p$-valued $q$-form, which we will denote
  $\bar \omega_\sigma$. Define
  $$
  I_p(\omega) = \sum_{|\ul\sigma| = p+1} \int_{\sigma} \bar
  \omega_{\sigma}
  $$
  This map defines a morphism of complexes by the Stokes theorem for
  superforms:
  $$
  I_p(d''\omega) = \sum_{|\ul\tau| = p+2} \int_{\tau} d\bar
  \omega_{\sigma} = \sum_{|\ul\tau| = p+2} \int_{\partial \tau} d\bar
  \omega_{\sigma} = \sum_{|\ul\tau| = p+2} \sum_{\sigma:
    \ul\sigma=\ul\tau \wo \{i\}} \sgn(u,\ul\tau) \int_{\sigma}
  \omega_{\sigma}
  $$
  
  Since this map defines a quasi-isomorphism between two resolutions
  of the sheaf $\Lambda^p$, it induces an isomorphism
  $H^q(\Aa^{p,\bu}(\Delta_X), d'')) \to H^q(C^\bu(\Delta_X,
  \Lambda^p_X \otimes \R))$ for any $q \geq 0$.
\end{proof}

\subsection{The monodromy morphism on superforms}

We will call \emph{local coordinates} at $\sigma$ any set of affine
functions $x_1, \ldots, x_n$ on $T(\sigma)$ such that
$dx_1, \ldots, dx_n$ form a base of a cotangent space of $T(\sigma)$
at any (equivalently, every) point of $T(\sigma)$.  Let $\sigma$ be a
face of $\Delta_X$, assume that $\Lambda^\bu$ is regular at $\sigma$
and let $x_1, \ldots, x_m$ be local coordinates at $\sigma$. Then
$$
x_0 = 1 - \sum_{i=1}^m x_i
$$
is an affine function on $T(\sigma)$ which is a restriction of a
linear function on $H^0(\St(\sigma), A^1)^*$ and $x_0, \ldots, x_m$ is
a basis of $H^0(\St(\sigma), A^1)$. 

\begin{lemma}
  \label{ap-lp-from-sf}
  For any $p > 0$ the space $H^0(\Stbar(\sigma), A^p \otimes \R)$ is
  isomorphic to the span of sections in
  $H^0(\St(\sigma), \Aa^{p-1,0})$ of the form
  $$
  \sum_{i \in I} \sgn(i, I) x_i d'x_{I \wo \{i\}}, \qquad x_0 d'x_J,
  $$
  for all multi-indices $I, J, |I|=p, |J|=p-1$. Moreover, if
  $a \in H^0(\St(\sigma), A^p)$, $\pi: A^p \otimes \R \to \Lambda^p \otimes \R$ is the
  natural projection and $a$ is represented by a superform $\eta$ then
  $\pi(a)$ is represented by the superform $d'\eta$.
\end{lemma}

\begin{proof}
  Recall that the sections in $H^0(\St(\sigma), \Lambda^p \otimes \R)$
  correspond to translation-invariant differential $p$-forms on
  $T(\sigma)$, sections in $H^0(\St(\sigma), A^p \otimes \R)$
  correspond to translation-invariant $p$-forms on
  $H^0(\St(\sigma), A^1 \otimes \R)^*$, and the morphism $\pi$ is
  induced by the restriction of forms to $T(\sigma)$.

  If $(p-1,0)$-superform $\eta$ on $T(\sigma)$ is given by a formula
  from the statement of the Lemma, for some coefficients
  $a_{I,i}, b_i$, then it can be lifted to a unique $(p-1,0)$-form
  $\tilde{\eta}$. Further, $d'\tilde{\eta}$ is a $d''$-closed
  $(p,0)$-superform, which is the same as a translation-invariant
  $p$-form, and so gives rise to a section of $A^p \otimes \R$ over
  $\St(\sigma)$. Clearly this map establishes a bijective
  correspondence between two natural bases of the subspace
  $H^0(\St(\sigma), A^1)$ and the space of translation-invariant
  $p$-forms on $H^0(\St(\sigma), A^1)^*$, and so is an isomorphism.
\end{proof}

\begin{prop}
  Assume that $\Lambda^\bu$ is regular at any face
  $\sigma \subset \Delta_X$. Then for any $p \geq 0$, there exists a
  distinguished triangle in the derived category of sheaves on
  $\Delta_X$
  $$
  \Lambda^p \to A^{p+1} \to \Lambda^{p+1} \xrightarrow{N} \Lambda^p[1]
  $$
  where the last morphism is given by the morphism
  $N: \Aa^{p+1,q} \to \Aa^{p, q+1}$.
\end{prop}

\begin{proof}
  It suffices to show that the 0-th cohomology of the cone complex
  $\Cone(N)$ is isomorphic to $A^p$, its cohomology in all other
  degrees vanishes since this is the case for the source and the
  destination of the morphism of complexes $N$.

  Let $x_1, \ldots, x_m$ be local coordinates at a face $\sigma$. If
  $\beta \in H^0(\St(\sigma), \Aa^{p,0}_{X})$ is a $d''$-closed
  superform then it is of the form $\sum_{|I|=p} b_I d'x_I$ for some
  constants $b_I$, and $N\beta$ is of the form
  $$
  \sum_{|I|=p-1} \sum_{j=1}^m (-1)^p \sgn(j, I \cup \{j\}) b_{I \cup \{j\}}
  d'x_I \wedge d''x_j.
  $$
  If $\alpha = \sum_{|I|=p-1} f_I d'x_I$ then $d''\alpha = N\beta$
  implies
  $$
  \dfrac{\partial f_I}{\partial x_j} = \sgn(j, I \cup \{j\})
  b_{|I| \cup \{j\}}
  $$
  for all $I$ such that $|I|=p-1$ and all $j, 1 \leq j \leq
  n$. Therefore,
  $$
  f_I = \left\{
    \begin{array}{ll}
      \sum_{j=1}^m \sgn(j, I \cup \{j\}) b_{I \cup \{j\}} x_j + a_I, & j \notin I,\\
      0, & j \in I.\\
    \end{array}\right.,
  $$
  for some constant $a_I$. To conclude observe that the space of forms
  $\alpha \in \Aa^{p-1, 0}$ such that there exists a form
  $\beta \in \Aa^{p,0}$ such that $d''\alpha = N \beta$ coincides with
  the space of froms from the statement of Lemma~\ref{ap-lp-from-sf}.
\end{proof}

\begin{cor}
  For any $p,q \geq 0$, $N$ is the coboundary morphism in the long exact
  sequence
  \begin{multline*}
    0 \to H^q(\Lambda^p) \to H^q(A^{p+1}) \to H^q(\Lambda^p)
    \xrightarrow{N}
    H^{q+1}(\Lambda^p) \to \\
    \to H^{q+1}(A^{p+1}) \to H^{q=1}(\Lambda^{p+1}) \xrightarrow{N}
    H^{q+2}(\Lambda^p) \to H^{q+2}(A^{p+1}) \to \ldots
  \end{multline*}
  associated to the short exact sequenece (\ref{seq:lambda'}).
\end{cor}

\section{Sheaves $\Lambda^p$ on Kulikov degenerations of K3 surfaces}
\label{s:kulikov}

\subsection{Singular affine structure and sheaf $\Lambda^1$}

A \emph{Kulikov degeneration $f: X \to S$} is a unipotent snc
degeneration such that the central fibre is reduced and $K_X=0$. By a
theorem of Kulikov, Persson and Pinkham \cite{kul77, pp} any unipotent
snc degeneration $f': X' \to S$ of K3 surfaces can be made a Kulikov
degeneration $f: X \to S$ after a bimerorphic modification
$X' \dashrightarrow X$ over $S$. Depending on the unipotency rank of
the monodromy, $\Delta_X$ is homeomorphic to a point, to an interval,
or to a 2-sphere. In the latter case, the case of maximally unipotent
monodromy, or \emph{Type~III degeneration}, $Y_i$ are rational
surfaces, double curves are smooth rational curves, and since
$\Delta_X$ is a manifold, the double curves on a given component $Y_i$
form a cycle.

An \emph{anticanonical pair} is the data of a smooth projective
surface $V$ and a divisor $D \in |-K_V|$ that is a sum of rational
curves $D_i$ forming a cycle and intersecting normally. By adjunction
formula one easily sees that $(Y_i, \sum Y_{ij})$ are anticanonical
pairs. To an anti-canonical pair one can associate a polygon, called
its \emph{pseudo-fan}, endowed with a singular affine structure
(\cite[Section~1.2]{ghk15}, \cite[Section~3]{eng18},
\cite[Section~8]{aet19}).

An \emph{affine structure} on a manifold $V$ is a flat torsion-free
connection on its tangent bundle. An \emph{integral affine structure}
is the data of an affine structure together with a flat $\Z$-local
system $T^\Z V \subset TV$ which spans $TV$. For any point $p \in V$
it is possible to find a local coordinate system $x_1, \ldots, x_n$ on
a neighbourhood $U$ of $p$, so that the action of $\nabla$ on
differential forms coincides with the de Rham differential. If
additionally $T^\Z V$ is chosen, one can choose such coordinate system
in a way that $\partial/\partial x_1, \ldots, \partial/\partial x_n$
generate $T^\Z V$. Given two such local charts, a transition function
between them belongs to $\SL_n(\Z) \rtimes \R^n$. Providing an atlas
of charts with such transition functions is equivalent to defining an
integral affine structure.

By a \emph{singular affine structure} on a surface we understand an
affine structure on a complement of finitely many points. A point is
called a singularity of a given singular affine structure if the
affine structure cannot be extended to this point.

\begin{defn}[Singular affine structure on a pseudo-fan]
  \label{affine-k3}
  A pseudo-fan of an anticanonical pair $(V,D)$ is a triangulated
  surface with boundary endowed with a singular integral affine
  structure on the interior as follows. The triangles are identified
  with triangles in $\R^2$, spanned by vectors corresponding to the
  irreducible components of $D$, all vectors originating in a fixed
  point $o$. The gluing is defined by identifying pairs of adjacent
  triangles, possibly with repetitions (in case $D$ is irrecudible)
  with pairs of triangles of lattice volume 1, spanned by three
  vectors $g_i, g_j, g_k$ (corresponding to a chain of irreducible
  components $D_i, D_j, D_k$) such that
  $$
  g_i + g_k = d_j \cdot g_j,
  $$
  where $d_j = -D_j^2$ if $D_k$ is a smooth rational curve and
  $d_j=-D_j^2+2$ if $D_j$ is a nodal rational curve.
\end{defn}

Consider a Kulikov degeneration $f: X \to S$.  If $o \in \Delta_X$,
the closed star $\Stbar(o)$ can be identified with the pseudo-fan of
the anti-canonical pair $(Y_o, \sum_i Y_{oi})$. By
\cite[Proposition~3.10]{eng18}, the singular integral affine
structures on stars of all vertices of $\Delta_X$ glue together to
singular affine structure on the whole of $\Delta_X$.

We are now going to compare this singular affine structure with the
one defined by the sheaf $\Lambda^1$ (since sheaves $\Lambda^1$ were
defined over $\Q$, we will not concern ourselves with integrality in
this discussion). In view of Definition~\ref{sheaves-non-snc} we can
pass to a subdivision of $\Delta_X$ by blowing up some strata, so that
we can ensure that there are no double curves that intersect in more
than one point. Since such blow-ups introduce components of
multiplicity $> 1$, we will have to modify the definition of the
affine structure on a pseudo-fan to take the multiplicities into
account. Note that we do not have to consider the degenerate case when
the divisor in an anticanonical pair is a nodal curve, since such
pairs do not occur as irreducible components of the central fibre of
Kulikov models.

\begin{defn}[Singular affine structure on $\Delta_X$]
  \label{affine-k3-subdiv}
  Let $f: X \to S$ be a degeneration of K3 surfaces obtained from a
  Kulikov model by a sequence of blow-ups of strata. For any
  irreducible component $Y_o$ of the central fibre the singular affine
  structure on $\St(o)$ is defined by specifying affine structure on
  stars of all 1-dimensional faces containing $o$ as follows. Let
  $\sigma, \eta \subset X$ be a pair of adjacent triangles,
  $\ul\sigma=\set{o,i,j}, \ul\eta=\set{o, j, k},$ define affine
  structure on $\sigma \cup \eta$ by identifying $\sigma$ and $\eta$
  with a pair of triangles in $\R^2$, of lattice volume
  $1/N_\sigma, 1/N_\eta$ respectively, spanned by three vectors
  $g_i, g_j, g_k$ such that
  $$
  N_ig_i + N_kg_k = d_{oj}\cdot N_jg_j,
  $$
  where $d_{oj} = -Y_{oj}^2$ on $Y_o$.
\end{defn}

\begin{lemma}
  \label{selfints}
  For any pair of adjacent triangles $\sigma, \eta$ as above
  $$
  N_o d_{jo} + N_j d_{oj} = N_i + N_k.
  $$  
\end{lemma}

\begin{proof}
  Follows immediately from writing out the restriction of
  $\sum N_i Y_i$, which is a principal divisor, to $Y_{oj}$.
\end{proof}

\begin{prop}
  Definition~\ref{affine-k3-subdiv} is consistent with respect to
  subdivisions of $\Delta_X$ induced by blow-ups of strata.
\end{prop}

\begin{proof}
  Let $g_i, g_j, g_k$ be vectors corresponding to double curves
  $Y_{oi}, Y_{oj}, Y_{ok}$. We need to consider the following
  modifications that will affect the subdivision of $\Delta_X$:
  \begin{enumerate}            
  \item blow up of any of $Y_{oi}, Y_{oj}, Y_{ok}$;
  \item blow up of any of $Y_{ij}, Y_{jk}$;
  \item blow up of any of the triple points $Y_{oij}, Y_{ojk}$.
  \end{enumerate}
  In all cases, denote $Y_e$ the exceptional divisor.

  \emph{Case (i)}. Let $Y_{ok}$ be the center of the blow-up for
  definiteness. Then $N_e = N_o + N_k, g_e = N_k/N_e g_k$ spans a
  triangle of volume $N_\eta \cdot N_k/N_e = 1/(N_o N_j N_e)$ with
  $g_j$ and
  $$
  N_j d_{oj} g_j = N_i g_i + N_k g_k = N_i g_i + N_e g_e,
  $$
  implying that the singular affine structure on the union of $\sigma$
  with the triangle $oje$ is the restriction of the singular affine
  structure on the union of $\sigma$ and $\eta$. The case when
  $Y_{oi}$ is the center of the blow-up is treated similarly.

  Let $Y_{oj}$ be the center of the blow-up. Then
  $N_e = N_o + N_j, g_e= N_j/N_e g_j$ spans triangles of volume
  $1/(N_o N_i N_e)$ and $1/(N_o N_k N_e)$ with $g_i$ and $g_k$,
  respectvely. We have $d_{oe}=d_{oj}$ and
  $$
  N_e d_{oe} g_e = N_j d_{oj} g_j  = N_i g_i + N_k g_k.
  $$

  Let us now check that the affine structure near $e$ is the
  restriction of the affine structure on $\sigma \cup \eta$ before
  subdivision. The exceptional divisor $Y_e$ is a projective bundle
  over $Y_{oj}$, the projectivisation of
  $\Oo_{Y_{oj}}(d_{oj}) \oplus \Oo_{Y_{oj}}(d_{jo}))$, and
  $Y_{ei}, Y_{ek}$ are its fibres, so $d_{ei}=d_{ek} = 0$. The curves
  $Y_{oe}, Y_{ej}$ are its sections and we have $d_{eo}=
  -d_{ej}=d_{jo} - d_{oj}$.

  We have the following vectors originating in $e$:
  $$
  h_i = g_i - g_e, h_j=g_j - g_e, h_k = g_k - g_e, h_o = - g_e
  $$

  Clearly,
  $$
  N_o h_o + N_j h_j = -N_o g_e + N_j g_j - N_j g_e = N_j g_j - N_e g_e
  = 0 = N_i \cdot 0 \cdot h_i = N_k \cdot 0 \cdot h_k,
  $$
  so the affine structure on the subdivided $\sigma$ as per
  Definition~\ref{affine-k3-subdiv} coincides with the flat affine
  structure on $\sigma$ before subdivision, and similarly for $\eta$.

  To check that the subdivision does not affect the affine structure
  on the adjacent halves of triangles $\sigma, \tau$, 
  apply Lemma~\ref{selfints} and the identity for $N_e
  d_{oe}g_e$ already established, so we get
  \begin{multline*}
    N_i h_i + N_k h_k = N_i g_i + N_k g_k - (N_i + N_k) g_e = N_e
    d_{oe} g_e - (N_od_{jo} + N_jd_{oj}) g_e = \\
    = (N_od_{oj} + N_j d_{oj} -N_od_{jo} - N_jd_{oj}) g_e=-N_o d_{eo}
    g_e = N_o d_{eo} h_o,
  \end{multline*}
  and analogously,
  \begin{multline*}
    N_i h_i + N_k h_k = N_i g_i + N_k g_k - (N_i + N_k) g_e = N_j
    d_{oj} g_j - (N_ed_{je} + N_jd_{ej}) g_e = \\
    = N_jd_{oj}g_j - N_j d_{jo} g_j  - N_jd_{ej}  g_e= N_j d_{ej}
    (g_j - g_e)  = N_j d_{ej} h_j.
  \end{multline*}

  \emph{Case (ii)}. Let $Y_{jk}$ be the center of the blow-up for
  definiteness, and let $Y_e$ be the exceptional divisor. The strict
  transform $Y'_o$ of $Y_o$ is a blow up of $Y_o$ in the triple point
  $Y_{ojk}$ with the exceptional curve $Y_{oe}$. Clearly, $d_{oe}=1$,
  and $N_e = N_j + N_k$. Denote $d'_{oj}, d'_{ok}$ the negative
  self-intersection numbers of the strict transforms of the double
  curves $Y_{oj}, Y_{ok}$. Clearly,
  $d'_{oj}=d_{oj}+1, d'_{ok}=d_{ok}+1$. Then putting
  $g_e = (N_j g_j + N_k g_k)/N_e$, we get 
  $$
  N_j g_j + N_k g_k = N_e g_e = N_e d_{oe} g_e,
  $$
  and
  $$
  N_i g_i + N_e g_e = N_i g_i + N_j g_j + N_k g_k = N_j (d_{oj} + 1)
  g_j = N_j d'_{oj} g_j,
  $$
  and observe that the subdivision does not affect the affine
  structure. 
  
  \emph{Case (iii)}. Let $Y_{ojk}$ be the center of the blow-up for
  definiteness. Denote as before $Y_e$ the exceptional divisor and
  $d'_{oj}, d'_{ok}$ the negative self-intersection numbers of the
  strict transforms of the corresponding double curves. Clearly,
  $N_e = N_o + N_j + N_k$, $d'_{oj}=d_{oj}+1$, and $Y_{eo}, Y_{ej},
  Y_{ek}$ is a triangle of lines on $Y_e \cong \P^2$.

  Let $g_e = (N_j g_j + N_k g_k) / N_e$. Then similarly to the
  previous case,
  $$
  N_j g_j + N_k g_k = N_e g_e = N_e d_{oe} g_e,
  $$  
  and
  $$
  N_i g_i + N_e g_e  = N_ig_i + N_j g_j + N_k g_k = N_j d'_{oj} g_j.
  $$

  Let $h_o=-g_e, h_j=g_j-g_e, h_k=g_k-g_e$, then the fact that the
  affine structure on the subdivided triangle $ojk$ is the same as the
  affine structure on it before subdivision folows from the equality
  $$
  N_o h_o + N_j h_j + N_k h_k = 0
  $$
  and the fact that double curves on $Y_e$ have self-intersection 1.
\end{proof}

Consider an anticanonical pair $(V,D)$, then a blow up at a point
$p \in D$ gives rise to an anticanonical pair $(V', D')$ where $D'$ is
the log pull-back of $D$. Denote $D'_i$ the strict transforms of
irreducible components $D_i$ of $D$ and let $D_e$ be the exceptional
divisor; denote $d_k, d'_k$ the self-intersection numbers of the
components on $V$, $V'$, respectively. If $p$ is a smooth point of
$D_j$ then $D' = \sum D'_i$ and $d'_i=d_i, i\neq j$ and $d'_j=d_j-1$
and $(V', D')$ is called an \emph{internal blow-up}. If $p$ is a node,
then the exceptional disivor $D_e$ is a summand of $D'$, $d'_e=1$ and
$d'_i=d_i+1$ for two components $D_i$ that contain $p$. In this case
$(V',D')$ is called a \emph{corner blow-up}.

By Proposition~2.12\cite{friedman15}, for any anti-canonical pair
$(V,D)$ there exists a sequence of corner blow-ups resulting in a pair
$(V',D')$ and a toric variety $\bar V$ with toric boundary $\bar D$
such that $(V',D')$ is obtained from $(\bar V, \bar D)$ by a sequence
of internal blow-ups. Notice that both Definition~\ref{affine-k3} and
Definition~\ref{affine-k3-subdiv} define a singular affine structure
on $\R^2$, which can be regarded as the cone over the boundary of a
pseudo-fan of a toric model or of $\Stbar(i)$.

\begin{prop}
  \label{lambda1-from-affine}
  If $f: X \to S$ is a Kulikov degeneration of K3 surfaces, then the
  sheaf $\Lambda^1$ on $\Delta_X$ coincides with the sheaf of affine
  functions with rational coefficients with respect to the singular
  affine structure given by Definition~\ref{affine-k3-subdiv}.
\end{prop}

\begin{proof}
  If $\alpha$ is a face of dimension 2 then by definition of the sheaf
  $\Lambda^1$, $\Lambda^1(\sigma)$ coincides with the space of affine
  functions on the interior of $\sigma$.

  Consider two trianges $\sigma$ and
  $\eta$ which share an edge $\alpha$. Let $\ul\sigma=\set{o,i,j},
  \ul\beta=\set{o,j,k}$. Then 
  $$
  \Lambda^1(\sigma) = \suchthat{ f \in \bar A^1(\alpha) }{ f(o)=0
    \textrm{ and } N_i f(i) Y_i + N_j f(j) Y_j + N_k f(k) Y_k
    \sim_{Y_\alpha} 0 }
  $$
  The condition on $f$ can be rewriten as
  $$
  N_i f(i) Y_i.Y_j + N_k f(k) Y_k.Y_j = N_i f(i) + N_k f(k) = N_j
  (-Y_j^2) f(j)
  $$
  which coincides with the condition imposed on $f$ to be affine in a
  neighbourhood of $\alpha$ with respect to the affine structure
  defined in Definition~\ref{affine-k3-subdiv}.

  Let $o$ be a vertex in $\Delta_X$. We perform enough blow-ups of
  strata so that $\Delta_X$ is a simplicial complex and so that the
  anti-canonical pair $(Y_o, \sum Y_{oi})$ can be obtained from a
  toric anti-canonical pair $(\bar Y_o, \sum \bar Y_{oi})$ by a
  sequence of internal blow-ups. We have
  $$
  \Lambda^1(o) = \suchthat{ f \in \bar A^1(o) }{ f(o)=0
    \textrm{ and } \sum N_i f(i) \Oo(Y_i)|_{Y_o}
    \sim_{Y_o} 0 }, 
  $$  
  and we need to check that $\Lambda^1(o)$ consists of affine
  functions with rational coefficients (with respect to the affine
  structure from Definition~\ref{affine-k3-subdiv}) that vanish at
  $o$.

  If $(Y_o, \sum Y_{oi})$ is toric, this statement is true by
  \cite[Proposition~3.9]{eng18} and Corollary~\ref{toric-tubular}.
  The star of $i$ is identified with a polygon in $\R^2$ with its
  standard integral affine structure, spanned by vectors
  $e_j, j \in \Stbar^0(i)$ and the elements of $\bar \Lambda^1(i)$ are
  identified with the functions that are piece-wise linear on the
  triangles that belong to $\Stbar(i)$, up to a constant. Such
  functions are completely determined by their values in points
  $j \in \Stbar^0(j)$.

  Let us analize how an internal blow-up affects the affine structure
  on the pseudo-fan of a toric canonical pair. A sequence of blow-ups
  of smooth points $Y_{ij}$ introduces a shearing transformation (see
  \cite[Section 3]{eng18}) with the effect that functions in a
  neighbourhood of $i$ that are affine with respect to the resulting
  singular affine structure must be constant along the vector
  $e_j$. If one performs internal blow-up on two components $Y_{ij},
  Y_{ik}$ such that $e_j, e_k$ are not collinear, then there are no
  non-constant affine functions.

  The general case follows from the following Claim.

  \emph{Claim}. Let $h: (Y',D') \to (Y,D)$ be an internal blow-up of a
  smooth point $p \in D_i$, let $I$ be the set of irreducible
  components of $D$ (and $D'$) and let $\set{N_i}_{i\in I}$ be a set
  of integers. Define the following sets
  $$
  \Aff_{Y'} = \suchthat{ f: I \to \Q }{ \sum N_i f(i) D'_i \sim_{Y'} 0
  },
  \Aff_Y = \suchthat{ f: I \to \Q }{ \sum N_i f(i) D_i  \sim_Y 0 },
  $$
  which we will also interpret as sets of functions on $\R^2$ that
  are linear on the cones spanned by $e_1, \ldots, e_j$. Then 
  $$
  f \in \Aff_{Y'} \textrm{ if and only if } f \in \Aff_Y,
  f(i) = 0.
  $$

  Indeed, since
  $$
  \sum N_i f(i) D'_i = h^*(\sum N_i f(i) D_i)) - N_i E
  $$
  and $E$ is linearly independent from the image of $\Pic(Y)$ in
  $\Pic(Y)$, the above expression can vanish if and only if both
  summands vanish.

\end{proof}

We will now show that the conditions of Theorem~\ref{hodge} are
satisfied at least for some Type~III Kulikov models.

\begin{prop}
  \label{comb-lef}
  Let $Y$ be an snc surface such that its irreducible components are
  rational. Then a combinatorial class $\omega \in H^2(Y)$ is
  combinatorial Lefschetz if and only if $(\omega|_{Y_i})^2 > 0$ and
  $\omega|_{Y_\sigma} \neq 0$ for all irreducible components $Y_i$ and
  double curves $Y_\sigma$.
\end{prop}

\begin{proof}
  Necessity is immediate.

  The class $\omega$ satisfying the conditions in the statement of the
  proposition clearly restricts to a Lefschetz class on double curves,
  so we only need to check whether its restrictions to $Y_i$ are
  Lefschetz. A rational surface has cohomology classes of type
  $(0,0), (1,1), (2,2)$ only, therefore any class with non-zero square
  satisfies the Lefschetz property. The Hodge-Riemann bilinear
  relations hold immediately on $H^{0,0}(Y_i)$ and $H^{2,2}(Y_i)$ by
  the positivity of $\omega^2$, and they hold on $H^{1,1}$ since the
  intersection form is negative definite on
  $\Ker L_\omega \subset H^2(Y_i)$ by Hodge index theorem.
\end{proof}

By the triple point formula \cite[2.1]{kul77}, for all double curves
in a Type~III Kulikov model we have
$$
(Y_{ij})^2|_{Y_i} + (Y_{ij})^2|_{Y_i} = -2.
$$

\begin{prop}
  \label{minus-one}
  Let $f:X \to S$ be a Type~III Kulikov model such that all double
  curves have square $-1$. Then there exists a combinatorial 
  Lefschetz class on $Y=f^{-1}(0)$.
\end{prop}

\begin{proof}
  Let $\omega_i = \sum_{j \in \Stbar^0(i)}
  c_1(\Oo_{Y_i}(Y_{ij}))$. Then $\omega_i|_{Y_{ij}} =
  \omega_j|_{Y_{ij}}$ since 
  $$
  \omega_i . Y_{ij} = 2 + Y_{ij}^2 = \omega_j . Y_{ij},
  $$
  so there exists a class $\omega \in H^2(Y)$ such that
  $\omega|_{Y_i}=\omega_i$. Clearly, this class is combinatorial and
  $$
  \omega_i^2 = \sum_{\sigma: i \in \sigma} 2 +
  \sum_{j \in \Stbar^0(i)} -1 > 0 \textrm{ and } \omega_i|_{Y_\sigma}
  \neq 0,
  $$
  so it is combinatorial Lefschetz by Proposition~\ref{comb-lef}.
\end{proof}

\subsection{Positive $(1,1)$-superforms}

Let $V$ be a vector space, a $(p,p)$-superform $\eta$ is called
\emph{symmetric} if $J\eta=\eta$. A symmetric $(p,p)$-superform $\eta$
is called \emph{(strongly) positive} \cite{lag12} if
$$
\eta = \sum_j f_j \cdot \alpha_{1j} \wedge J\alpha_{1j} \wedge \ldots
\wedge \alpha_{pj} \wedge J\alpha_{pj},
$$
where $f_j \geq 0$ and $\alpha_{ij} \in \Aa^{1,0}(V)$. A section
$\omega \in H^0(\Delta_X, \Aa^{1,1})$ is called \emph{positive} if its
restriction to each $\Stbar(\sigma)$ gives rise to a positive
(1,1)-superform on $e(\Stbar(\sigma)) \subset T(\sigma)$.

One can observe that the data of a positive symmetric
$(1,1)$-superform $\omega$ is equivalent to that of a pseudo-metric
tensor on $\Delta_X$.

Recall that a function $f: U \to \R$ on a convex domain
$U \subset \R^n$ is convex is it is continuous and all sublevel sets
$U_c := \suchthat{ x \in U}{f(x) \leq c }$ are convex. It is
\emph{strictly convex} if $U_c$ are strictly convex, that is, if every
line segment that lies in $U_c$ is contained in the interior of $U_c$
except maybe its endpoints. Let $\Sigma \subset \R^n$ be a polyhedral
complex, then a PL function on $\Sigma$ is called \emph{strictly
  conex} if it is convex and additionally its restriction to a
neighbourhood of any polyhedron in $\Sigma$ is not linear.

This definition makes sense on $\Delta_X$ when $X$ is a Kulikov
degeneration. In this case we call a function on $\St(i)$ strictly
convex if its restrictions to $\St(\sigma)$ are convex for all
$\sigma$ containing $i$ --- where $\St(\sigma)$ is identified with a
pair of triangles in $\R^2$ using the affine structure on $\Delta_X$.

\begin{prop}
  \label{form-plmetric}
  Let $X \to S$ be a Kulikov degeneration.  If $\omega$ is a symmetric
  $d''$-closed $(1,1)$-superform on $\Delta_X$, let $(a_\sigma)$ be a
  cocycle that corresponds to $\omega$ under the quasi-isomorphism
  from Proposition~\ref{q-iso} and let $L$ be the corresponding
  virtual line bundle. Assume that away from the singularities of
  affine structure $\omega$ is positive definite, then there exists a
  PL convex metrization of $L$.
\end{prop}

\begin{proof}
  A PL metrization of $L$ is a collection of sections
  $h_i \in H^0(\St(i), \bar \Lambda^1)$ such that
  $h_i - h_j = a_\sigma$ for all one-dimensional faces
  $\sigma \subset \Delta_X, \ul\sigma={i,j}$. We will regard $h_i$ as
  PL functions vanishing at 0 and linear on faces contaning $i$.

  For each vertex $i$ take a function $f_i \in \Aa^{0,0}(\St(i))$ such
  that $d'd''f_i = \alpha|_{\St(i)}$ and
  $$
  \dfrac{\partial f_i}{\partial x_1}(i)=\dfrac{\partial f_i}{\partial x_2}(i) = 0,
  $$
  where $x_1, x_2$ is some coordinate system at $i$. To define $h_i$,
  suffices to define them on vertices that belong to $\Stbar(i)$:
  $$
  h_i(j) = f_i(j).
  $$
  Since $d'd''f_i=d'd''f_j=0$, the functions $f_i$ and $f_j$ have the
  same Hessian matrix on $\Stbar(i) \cap \Stbar(j)$ and we have that
  $$
  f_i - f_j = \dfrac{\partial f_i}{\partial x_1} x_1 + \dfrac{\partial
    f_i}{\partial x_2} x_2 + c
  $$
  for some constant $c \in \R$. But the linear part of the difference
  is precisely the value of the cocycle $a_{ij}$.
\end{proof}

\begin{cor}
  Let $T$ be a positive supercurrent on $\Delta_X$. Then there exists
  a virtual line bundle $L$ on $\Delta_X$ and a PL convex metrization
  $h$.
\end{cor}

\begin{proof}
  By \cite[Theorem~1.13]{lag12} there exist non necessarily smooth
  functions $f_i$ such that $d'd''f_i = T$ on $\Stbar(i)$ in the sense
  of currents. Pick some such functions and define, as in the proof of
  Proposition~\ref{form-plmetric},
  $$
  h_i(j) = f_i(j)
  $$
  It is then clear from the definition that
  $$
  a_{ij} = [h_j - h_i]
  $$
  forms a cocyle in $H^1(\Delta_X, \Lambda^1)$ that gives rise to a
  virtual line bundle $L$ with the PL metrization defined by the
  functions $h_i$.
\end{proof}

\begin{prop}
  \label{kahler-convex}
  Let $f: X \to S$ be a Kulikov degeneration, let $L$ be a virtual
  line bundle on $\Delta_X$ and $h$ be a PL metric on $L$. If $h$ is
  stritcly convex  then $c_1(L,h)_i \in H^2(Y_i)$ is combinatorial 
  Lefschetz.
\end{prop}

\begin{proof}
  Let $h_i \in H^0(\Stbar(i), \bar\Lambda^1)$ be a trivialization of
  $h$ in $\Stbar(i)$. As before, we will identify $h_i$ with
  piece-wise linear functions such that $h_i(i) = 0$.  We need to show
  that the  class
  $$
  c_1(L,h)_i = \sum_{j \in \Stbar^0(i)} N_i h(i) c_1(\Oo_{Y_j}(Y_{ij})) \in H^2(Y_i)
  $$
  has a positive square. 
  $$
  c_1(L,h)_i^2 = \sum_{o \in \Stbar^0(i)} N_j^2 h(j)^2 Y_{ij}^2 +
  2\sum_{\stackrel{\exists
      \sigma\subset\Stbar(i)}{\ul\sigma=\set{i,j,k}}} N_j N_k h(j) h(k).
  $$
  Since $c_1(L,h)_i$ does not change when a linear function is added
  to $h_i$, we may assume that $h_i$ strictly positive on all vertices
  of $\Stbar^0(i)$. Since $h_i$ is strictly convex, from the
  definition of affine structure on $\Delta_X$ for any triple of
  adjacent vertices $j,k,l \in \Stbar^0(i)$ we have
  $$
  N_j h(j) + N_k (Y_{ik}^2) \cdot h(k) + N_l h(l) > 0,
  $$
  where $Y_{ik}^2$ is the self-intersection of the curve $Y_{ik}$ on
  $Y_i$, and hence
  $$
  N_k N_j h(k) h(j) + N_k^2 (Y_{ik}^2) \cdot (h(k))^2 + N_k N_l h(k)
  h(l) > 0,
  $$
  summing up these expressions for all triples $j,k,l$ of adjacent
  vertices in $\Stbar^0(i)$ we obtain the expression for $c_1(L,h)_i$,
  and therefore $c_1(L,h)_i > 0$.
\end{proof}

\subsection{Simple affine structure singularities}

Let $X \to S$ be a Kulikov degeneration and let $Y_i$ be an
irreducible component of the central fibre. Take an anticanonical pair
$$
(Y_i, D), \qquad D = \sum_{\exists \sigma: \ul\sigma=\set{i,j}} Y_j \cap
Y_i,
$$
and assume that it is obtained from a toric anticanonical pair
$(\bar Y_i, \bar D)$ by a single blow-up of a smooth point
$p \in Y_{ij} \subset \bar D$. Applying
Definition~\ref{affine-k3-subdiv} one observes that the monodromy matrix
of the affine structure around $i$ is
$$
\left(
  \begin{array}{cc}
    1 & 1 \\
    0 & 1 
  \end{array}
\right)
$$
in the basis consisting of vectors $ij, ik$. This is the simplest
affine structure singularity.

Now assume that all non-toric irreducible components $Y_i$ become
toric after blowing down on every one of them of a single exceptional
curve $C_i$ that intersects the double locus of $Y$ in one point. Let
$X'$ be the blow up in $X$ of all such curves $C_i$ on non-toric
irreducible components $Y_i$ of the central fibre $Y$. The morphisms
$$
H^q(\Delta_{X'}, \Lambda^p_{X'}) \cong \gr^W_{2p} H^{p+q}(X'_\infty)
$$
are isomorphisms since $\Hh^{2\bu}(Y'_\sigma) = H^{2\bu}(Y'_\sigma)$
for any stratum $Y'_\sigma$ of $Y$. We will show in this section that
the cohomology of $\Lambda^0, \Lambda^1, \Lambda^2$ does not change
when one passes to $\Delta_{X'}$. In the paper \cite{hessian} I show
that $\Delta_X$ admits a K\"ahler superform and therefore by
Propositions~\ref{form-plmetric} and \ref{kahler-convex} the morphisms
$$
H^q(\Delta_X, \Lambda^1_X) \cong \gr^W_{2p} H^{p+q}(X_\infty)
$$
are isomorphisms too. This is parallel to the results of Ruddat
\cite{ruddat2010log} that show that affine cohomology recovers full
nearby fibre cohomology in case of toric degenerations obtained from
manifolds with simple affine structure singularities.

\begin{prop}
  \label{k3-coho}
  The sheaf $\Lambda^\bu_{X'}$ is regular at any face 
  $\sigma \subset \Delta_{X'}$ and there exists a surjective morphism
  $$
  H^1(\Delta_X, \Lambda^1_X) \to H^1(\Delta_{X'}, \Lambda^1_{X'}).
  $$
  and an injective morphism
  $$
  H^0(\Delta_X, \Lambda^1_X) \hookrightarrow H^0(\Delta_{X'},
  \Lambda^1_{X'}).
  $$
  For $p \in \set{0,2}$ and any $q$ there exist isomorphisms
  $$
  H^q(\Delta_X, \Lambda^p_X) \cong H^q(\Delta_{X'}, \Lambda^p_{X'}) .
  $$
  
\end{prop}

\begin{proof}
  Since the modification $X' \to X$ is performed in a similar manner
  in a neighbourhood of each non-toric component $Y$, so below we will
  concentrate on the situation and fix notation for the strata near
  one such component, $Y_i$.
  
  Let $Y_j$ be the irreducible component of $Y$ that contains the
  double curve $Y_{ij}$ that $C_i$ intersects, and let $Y_k$ and $Y_l$
  be the irreducible components of $Y$ that intersect $Y_{ij}$ in
  triple points.

  We may assume that $Y_j$ is toric and that the double curves on it
  are irreducible components of the toric boundary divisor, otherwise,
  reduce to this situation by doing a base change with respect to a
  finite cover $S \to S$ of sufficiently high degrre and resolving the
  arising singularities.

  Denote $Y'$ the central fibre of $X'$ and $Y'_i$ the strict
  transforms of the strata of $Y$, and denote $Y'_e$ the exceptional
  divisor, whose multiplicity is $N_e = N_i + N_j$. Note that
  $Y'_i \cong Y_i$ and that $Y'_j$ is a blow-up of $Y_j$ in the triple
  point $C \cap Y_{ij}$. The normal bundle of $C$ in $X$ is a direct
  sum of the normal bundle of $C$ in $Y_i$, which is $\Oo_C(-1)$, and
  the restriction of the normal bundle of $Y_i$ in $X$ to $C$ which is
  equal to
  $$
  \Oo_{Y_i}(-\sum_{o \in \Stbar^0(i) \wo \{i\}} N_o Y_o)|_{C} = \Oo_C(-1).
  $$
  It follows that $Y'_e \cong \P^1 \times \P^1$.

  Since the stars of vertices $k$ and $l$ and the corresopnding stata
  do not change when passing from $\Delta_X$ to $\Delta_{X'}$, the
  sheaves $\Lambda^i, i=1,2$ also do not change near these vertices,
  and we only need to compute the sections of these sheaves near faces
  $i, j, e, ij, j, ie, je$.

  We have
  \begin{multline*}
  \Oo_{X'}(Y'_e)|_{Y'_{ie}} = \Oo_{Y'_{ie}}(-1), \ \   
  \Oo_{X'}(Y'_i)|_{Y'_{ie}} = \Oo_{Y_{ie}}(-1), \ \  \Oo_{X'}(Y'_j)|_{Y'_{ie}} =
  \Oo_{Y'_{ie}}(1),  \\
  \Oo_{X'}(Y'_e)|_{Y'_{je}} = 0, \ \ 
  \Oo_{X'}(Y'_i)|_{Y'_{je}} = \Oo_{Y'_{je}}(1), \ \  \Oo_{X'}(Y'_j)|_{Y'_{je}} =
    \Oo_{Y'_{ie}}(-1).
  \end{multline*}
  From the considerations above one gets immediately that
  \begin{multline*}
  \Lambda^1_{X'}(e) = 0, \qquad 
  \Lambda^1(ie) = \Lambda^1(je) = \suchthat{ f: \set{i,j,e} \to \Q }{
    f(i) = f(j) } / \const.
  \end{multline*}
  Any element of $\bar \Lambda^2(ie)$ or $\bar \Lambda^2(je)$ is
  proportional to $i \wedge j$ but
  $$
  c^2_{ie}(i \wedge j) = c^1_{ie}(j) \otimes i \qquad c^2_{je}(i
  \wedge j) = c^1_{ie}(i) \otimes j
  $$
  and both expressions are non-zero. Therefore,
  $$
  \Lambda^2(e) = \Lambda^2(ie) = \Lambda^2(je) = 0.
  $$
  First of all, we have
  $$
  \Oo_{X'}(Y'_k)|_{Y'_{ij}} = \Oo_{X'}(Y_k)|_{Y_{ij}} =
  \Oo_{Y_{ij}}(1),
  \qquad \Oo_{X'}(Y'_l)|_{Y'_{ij}} = \Oo_X(Y_l)|_{Y_{ij}} =
  \Oo_{Y_{ij}}(1).
  $$

  Since $Y'_j$ is a blow-up of $Y_j$, $Y'_i=Y_i$ and $Y'_{ij}=Y_{ij}$,
  we have
  \begin{eqnarray*}
  \Oo_{X'}(Y'_i)|_{Y_{ij}} = \Oo_{Y'_{ij}}((Y'_{ij})|_{Y'_i}^2) =
    \Oo_{Y_{ij}}((Y_{ij})|_{Y_i}^2) = \Oo_{X}(Y_i)|_{Y_{ij}},  \\
  \Oo_{X'}(Y'_j)|_{Y_{ij}} = \Oo_{Y'_{ij}}((Y'_{ij})|_{Y'_j}^2) =
    \Oo_{Y_{ij}}((Y_{ij})|_{Y_j}^2 - 1) = \Oo_{X}(Y_j)|_{Y_{ij}}
    \otimes \Oo_{Y_{ij}}(-1).
  \end{eqnarray*}
  
  If we identify elements of $\Lambda^1_X(ij)$,
  resp. $\Lambda^1_{X'}(ij)$, with maps $f$, resp. $f'$, on finite
  sets $\set{i,j,k,l}$, resp. $\set{i,j,k,l,e}$, up to constant maps
  then we see that there is a natural inclusion
  $$
  \iota: \Lambda^1_X(ij) \hookrightarrow \Lambda^1_{X'}(ij), \ \  [f] \mapsto [f'],
  \textrm{ where } f'(e) =  f(i) + f(j),\ f'|_{\set{i,j,k,l}}
    = f|_{\set{i,j,k,l}}.
  $$
  The sections of $\Lambda^1(ij)$ are of the form $f' + g$ where $f
  \in \Lambda^1_X(ij)$ and $g$ is a function supported on
  $\set{i,j,e}$ such that 
  $$
  (g(e) \Oo_{X'}(Y'_e) + g(i)\Oo_{X'}(Y'_i) + g(j)
  \Oo_{X'}(Y'_j))|_{Y'_{ij}} = 0
  $$

  To compute $\Lambda^2(ij)$ note that any element of
  $\bar \Lambda^2(ij)$ is represented by a tensor of the form
  $i \wedge a$ where $a$ is a linear combination of $e,k,l$ and that
  $$
  c^2_{ij}(i \wedge a) = c^1_{ij}(a) \otimes i.
  $$
  It follows that
  $$
  \Lambda^2_{X'}(ij) = \suchthat{ i \wedge a }{ a|_{ij} = 0, a \in \Lambda^1_{X'}(ij)}.
  $$
  In particular, $\dim \Lambda^2_{X'}(ij) = 1$.
  \vspace{1ex}
  



  From the computations of $\Lambda^2_{X'}(ij)$ we conclude that
  $\Lambda^2_X$ is a constant sheaf in the neighbourhood of $i$ and
  $j$.

  Summing up, the sections of the sheaf $\Lambda^1_{X'}$ near vertices
  $i$ and $j$ are sums of pull-backs of sections on
  $\Delta_X \subset \Delta_{X'}$ that are affine with respect to a
  certain non-singular affine structure under a certain natural
  projection map $\Delta_{X'} \to \Delta_X$ that collapses the
  triangle $ije$ onto the its edge $ij$, and sections supported
  outside $\Delta_X$. The sheaf $\Lambda^2_{X'}$ is a push-forward of
  the constant sheaf along the open embedding of the complement of the
  boundary of the space $\Delta_{X'}$ into $\Delta_{X'}$. It is also
  clear from the computations that $\Lambda^\bu$ is regular at every
  face of $\Delta_{X'}$.

  Since $\Delta_X$ is homotopy equivalent to $\Delta_{X'}$, there
  exist isomorophisms
  $H^\bu(\Delta_X, \Lambda^p_X) \cong H^\bu(\Delta_{X'}, \Lambda^p_{X'})$
  for $p=0,2$. Is is left to treat the $\Lambda^1$ case.
  
  Consider an inclusion $C^\bu(\Delta_X, \Lambda^1) \hookrightarrow
  C^\bu(\Delta_{X'}, \Lambda^1_{X'})$ of complexes
  $$
  \xymatrix{
    0 \ar[r] & \bigoplus\limits_{u \in \Delta_X} \Lambda^1(u) \ar[r]^d\ar[d]^{\iota} &
    \bigoplus\limits_{\stackrel{|\sigma|=2}{\sigma \subset \Delta_X}} \Lambda^1(\sigma)
    \ar[r]^d\ar[d]^{\iota} & 
    \bigoplus\limits_{\stackrel{|\tau|=3}{\tau \subset \Delta_X}} \Lambda^1(\tau)
    \ar[r]\ar[d]^{\iota} &
    0 \\
    0 \ar[r] & \bigoplus\limits_{v \in \Delta_{X'}} \Lambda^1(v) \ar[r]^d&
    \bigoplus\limits_{\stackrel{|\alpha|=2}{\alpha \subset \Delta_{X'}}} \Lambda^1(\alpha)
    \ar[r]^d & 
    \bigoplus\limits_{\stackrel{|\eta|=3}{\eta \subset \Delta_{X'}}} \Lambda^1(\eta)
    \ar[r] &
    0 
  }  
  $$  
  Since there are no coboundaries in degree 0 the first map is an
  inclusion. 

  For the last statement we only need to analyze the cocycles with
  coefficients in $\Lambda^1_X$ and $\Lambda^1_{X'}$ near the glued in
  triangles, since away from them the sheaves are isomorphic. We keep
  the notation for vertices and strata from the computations above.

  Let $(a_\alpha) \in \oplus_{|\ul\alpha|=2} \Lambda^1_X(\tau)$ be a
  cocycle. Since a cocycle must satisfy in particular
  $$
  d(a_\alpha)|_{\set{i,j,e}} = 0 
  $$
  but
  $$
  \Lambda^1(ij)|_{\set{i,j,e}} \cap \Lambda^1(ie) =
  \Lambda^1(ij)|_{\set{i,j,e}} \cap \Lambda^1(je) = 0,
  $$
  we have then $a_{ij}|_{\set{i,j,e}} = 0$. In particular, $(a_\alpha)
  \in \Im \iota$.
    
\end{proof}

Note that since
$H^0(\Delta_{X'},\Lambda^1_{X'}) \cong \gr^W_2 H^1(X'_\infty) = 0$,
$\iota: H^0(\Delta_X,\Lambda^1_X) \to H^0(\Delta_{X'},\Lambda^1_{X'})
$ is an isomorphism.


\footnotesize

\bibliography{coho}

\begin{thebibliography}{dFKX17}

\bibitem[ABW13]{abw13}
Donu Arapura, Parsa Bakhtary, and Jaros{\l}aw W{\l}odarczyk.
\newblock Weights on cohomology, invariants of singularities, and dual
  complexes.
\newblock {\em Mathematische Annalen}, 357(2):513--550, 2013.

\bibitem[AET19]{aet19}
Valery Alexeev, Philip Engel, and Alan Thompson.
\newblock Stable pair compactification of moduli of {K3} surfaces of degree 2.
\newblock {\em arXiv preprint arXiv:1903.09742}, 2019.

\bibitem[Ber99]{berkovich99}
Vladimir~G Berkovich.
\newblock Smooth p-adic analytic spaces are locally contractible.
\newblock {\em Inventiones mathematicae}, 137(1):1--84, 1999.

\bibitem[Ber09]{berk09}
Vladimir~G Berkovich.
\newblock A non-archimedean interpretation of the weight zero subspaces of
  limit mixed {Hodge} structures.
\newblock In {\em Algebra, Arithmetic, and Geometry}, pages 49--67. Springer,
  2009.

\bibitem[BJ17]{bj17}
S{\'e}bastien Boucksom and Mattias Jonsson.
\newblock Tropical and non-archimedean limits of degenerating families of
  volume forms.
\newblock {\em Journal de l’{\'E}cole polytechnique—Math{\'e}matiques},
  4:87--139, 2017.

\bibitem[CGP19]{cgp2}
Melody Chan, Soren Galatius, and Sam Payne.
\newblock Topology of moduli spaces of tropical curves with marked points.
\newblock {\em arXiv preprint arXiv:1903.07187}, 2019.

\bibitem[CGP21]{cgp1}
Melody Chan, S{\o}ren Galatius, and Sam Payne.
\newblock Tropical curves, graph complexes, and top weight cohomology of $m_g$.
\newblock {\em Journal of the American Mathematical Society}, 34(2):565--594,
  2021.

\bibitem[CL12a]{cld}
Antoine Chambert-Loir.
\newblock Differential forms and currents on berkovich spaces.
\newblock {\em arXiv preprint arXiv:1204.6277}, 2012.

\bibitem[CL12b]{dacl12}
Antoine Chambert-Loir.
\newblock Differential forms and currents on berkovich spaces.
\newblock {\em arXiv preprint arXiv:1204.6277}, 2012.

\bibitem[dCM02]{dcm02}
Mark Andrea~A de~Cataldo and Luca Migliorini.
\newblock The hard {Lefschetz} theorem and the topology of semismall maps.
\newblock {\em Annales scientifiques de l'Ecole normale sup{\'e}rieure},
  35(5):759--772, 2002.

\bibitem[Del71]{dhodge2}
Pierre Deligne.
\newblock Th{\'e}orie de {Hodge}: {II}.
\newblock {\em Publications Math{\'e}matiques de l'IH{\'E}S}, 40:5--57, 1971.

\bibitem[dFKX17]{dfkx}
Tommaso de~Fernex, J{\'a}nos Koll{\'a}r, and Chenyang Xu.
\newblock The dual complex of singularities.
\newblock In {\em Higher Dimensional Algebraic Geometry: In honour of Professor
  Yujiro Kawamata's sixtieth birthday}, pages 103--129. Mathematical Society of
  Japan, 2017.

\bibitem[Eng18]{eng18}
Philip Engel.
\newblock Looijenga's conjecture via integral-affine geometry.
\newblock {\em Journal of Differential Geometry}, 109(3):467--495, 2018.

\bibitem[Fri15]{friedman15}
Robert Friedman.
\newblock On the geometry of anticanonical pairs.
\newblock {\em arXiv preprint arXiv:1502.02560}, 2015.

\bibitem[GA90]{gna90}
Francisco Guill{\'e}n and V~Navarro Aznar.
\newblock Sur le th{\'e}oreme local des cycles invariants.
\newblock {\em Duke Mathematical Journal}, 61(1):133--155, 1990.

\bibitem[GHK15]{ghk15}
Mark Gross, Paul Hacking, and Sean Keel.
\newblock Mirror symmetry for log calabi-yau surfaces {I}.
\newblock {\em Publications Math{\'e}matiques de l'{IH{\'E}S}}, 122(1):65--168,
  2015.

\bibitem[Got22]{goto22}
Keita Goto.
\newblock On the two types of affine structures for degenerating kummer
  surfaces-non-archimedean vs gromov-hausdorff limits.
\newblock {\em arXiv preprint arXiv:2203.14543}, 2022.

\bibitem[GS10]{gs2}
Mark Gross and Bernd Siebert.
\newblock Mirror symmetry via logarithmic degeneration data, {II}.
\newblock {\em Journal of Algebraic Geometry}, 19(4):679--780, 2010.

\bibitem[Gub16]{gubler16}
Walter Gubler.
\newblock Forms and currents on the analytification of an algebraic variety
  (after chambert-loir and ducros).
\newblock In {\em Nonarchimedean and tropical geometry}, pages 1--30. Springer,
  2016.

\bibitem[Hat02]{hatcher}
Allen Hatcher.
\newblock {\em Algebraic Topology}.
\newblock Cambridge University Press, 2002.

\bibitem[IKMZ19]{ikmz}
Ilia Itenberg, Ludmil Katzarkov, Grigory Mikhalkin, and Ilia Zharkov.
\newblock Tropical homology.
\newblock {\em Mathematische Annalen}, 374(1-2):963--1006, 2019.

\bibitem[Ive86]{iversen}
Birger Iversen.
\newblock {\em Cohomology of sheaves}.
\newblock Springer Science \& Business Media, 1986.

\bibitem[Jel16]{jell16}
Philipp Jell.
\newblock A poincar{\'e} lemma for real-valued differential forms on berkovich
  spaces.
\newblock {\em Mathematische Zeitschrift}, 282(3):1149--1167, 2016.

\bibitem[Jel19]{jell19}
Philipp Jell.
\newblock Tropical {Hodge} numbers of non-archimedean curves.
\newblock {\em Israel Journal of Mathematics}, 229(1):287--305, 2019.

\bibitem[JRS17]{jrs17}
Philipp Jell, Johannes Rau, and Kristin Shaw.
\newblock Lefschetz (1, 1)-theorem in tropical geometry.
\newblock {\em arXiv preprint arXiv:1711.07900}, 2017.

\bibitem[JSS19]{jss19}
Philipp Jell, Kristin Shaw, and Jascha Smacka.
\newblock Superforms, tropical cohomology, and {Poincar{\'e}} duality.
\newblock {\em Advances in Geometry}, 19(1):101--130, 2019.

\bibitem[Kas84]{kashiwara}
Masaki Kashiwara.
\newblock The {Riemann}-{Hilbert} problem for holonomic systems.
\newblock {\em Publications of the Research Institute for Mathematical
  Sciences}, 20(2):319--365, 1984.

\bibitem[KLSV18]{klsv18}
J{\'a}nos Koll{\'a}r, Radu Laza, Giulia Sacc{\`a}, and Claire Voisin.
\newblock Remarks on degenerations of hyper-k{\"a}hler manifolds.
\newblock {\em Annales de l'Institut Fourier}, 68(7):2837--2882, 2018.

\bibitem[KS06]{ks06}
Maxim Kontsevich and Yan Soibelman.
\newblock Affine structures and non-{Archimedean} analytic spaces.
\newblock {\em Progress in mathematics}, 244:321, 2006.

\bibitem[KT02a]{kt}
Maxim Kontsevich and Yury Tschinkel.
\newblock Non-archimedean k{\"a}hler geometry.
\newblock {\em Unpublished note}, 2002.

\bibitem[KT02b]{kt02}
Maxim Kontsevich and Yury Tschinkel.
\newblock Non-archimedean {K{\"a}hler} geometry.
\newblock {\em Unpublished manuscript}, 2002.

\bibitem[Kul77]{kul77}
Viktor~S Kulikov.
\newblock Degenerations of {K3} surfaces and {Enriques} surfaces.
\newblock {\em Mathematics of the USSR-Izvestiya}, 11(5):957, 1977.

\bibitem[K{\"u}n98]{kunnemann98}
Klaus K{\"u}nnemann.
\newblock Projective regular models for abelian varieties, semistable
  reduction, and the height pairing.
\newblock {\em Duke mathematical journal}, 95(1):161--212, 1998.

\bibitem[KX16]{kx16}
J{\'a}nos Koll{\'a}r and Chenyang Xu.
\newblock The dual complex of calabi--yau pairs.
\newblock {\em Inventiones mathematicae}, 205(3):527--557, 2016.

\bibitem[Lag11]{lag11}
Aron Lagerberg.
\newblock $l^2$-estimates for the $d$-operator acting on super forms.
\newblock {\em arXiv preprint arXiv:1109.3983}, 2011.

\bibitem[Lag12a]{lagerberg12}
Aron Lagerberg.
\newblock Super currents and tropical geometry.
\newblock {\em Mathematische Zeitschrift}, 270(3):1011--1050, 2012.

\bibitem[Lag12b]{lag12}
Aron Lagerberg.
\newblock Super currents and tropical geometry.
\newblock {\em Mathematische Zeitschrift}, 270(3):1011--1050, 2012.

\bibitem[Liu17]{liu17}
Yifeng Liu.
\newblock Monodromy map for tropical {Dolbeault} cohomology.
\newblock {\em arXiv preprint arXiv:1704.06949}, 2017.

\bibitem[Mau20]{mauri20}
Mirko Mauri.
\newblock The dual complex of log calabi--yau pairs on mori fibre spaces.
\newblock {\em Advances in Mathematics}, 364:107009, 2020.

\bibitem[MN15]{mn15}
Mircea Musta\c{t}a and Johannes Nicaise.
\newblock Weight functions on non-archimedean analytic spaces and the
  kontsevich–soibelman skeleton.
\newblock {\em Algebraic Geometry}, 2(3):365--404, 2015.

\bibitem[MZ14]{mz14}
Grigory Mikhalkin and Ilia Zharkov.
\newblock Tropical eigenwave and intermediate {Jacobians}.
\newblock In {\em Homological mirror symmetry and tropical geometry}, pages
  309--349. Springer, 2014.

\bibitem[NX16]{nicaise16}
Johannes Nicaise and Chenyang Xu.
\newblock The essential skeleton of a degeneration of algebraic varieties.
\newblock {\em American Journal of Mathematics}, 138(6):1645--1667, 2016.

\bibitem[NXY19]{nicxuyu19}
Johannes Nicaise, Chenyang Xu, and Tony~Yue Yu.
\newblock The non-archimedean {SYZ} fibration.
\newblock {\em Compositio Mathematica}, 155(5):953--972, 2019.

\bibitem[OO18]{oo18}
Yuji Odaka and Yoshiki Oshima.
\newblock Collapsing {K3} surfaces, tropical geometry and moduli
  compactifications of {Satake}, {Morgan-Shalen} type.
\newblock {\em arXiv preprint arXiv:1810.07685}, 2018.

\bibitem[Pay08]{payne}
Sam Payne.
\newblock Analytification is the limit of all tropicalizations.
\newblock {\em arXiv preprint arXiv:0805.1916}, 2008.

\bibitem[PP81]{pp}
Ulf Persson and Henry Pinkham.
\newblock Degeneration of surfaces with trivial canonical bundle.
\newblock {\em Annals of Mathematics}, pages 45--66, 1981.

\bibitem[PS08]{ps08}
Chris Peters and Joseph Steenbrink.
\newblock {\em Mixed {Hodge} structures}, volume~52.
\newblock Springer Science \& Business Media, 2008.

\bibitem[Rud10]{ruddat2010log}
Helge Ruddat.
\newblock Log hodge groups on a toric calabi-yau degeneration.
\newblock {\em Mirror Symmetry and Tropical Geometry, Contemporary
  Mathematics}, 527:113--164, 2010.

\bibitem[Sch73]{schmid}
Wilfried Schmid.
\newblock Variation of {Hodge} structure: the singularities of the period
  mapping.
\newblock {\em Inventiones mathematicae}, 22(3):211--319, 1973.

\bibitem[Sol18]{soldatenkov}
Andrey Soldatenkov.
\newblock Limit mixed {Hodge} structures of hyperk\"ahler manifolds.
\newblock {\em arXiv preprint arXiv:1807.04030}, 2018.

\bibitem[SRJ18]{jsr-epi}
Kristin Shaw, Johannes Rau, and Philipp Jell.
\newblock Lefschetz (1, 1)-theorem in tropical geometry.
\newblock {\em {\'E}pijournal de G{\'e}om{\'e}trie Alg{\'e}brique}, 2, 2018.

\bibitem[Ste76]{ste76}
Joseph Steenbrink.
\newblock Limits of {Hodge} structures.
\newblock {\em Inventiones mathematicae}, 31(3):229--257, 1976.

\bibitem[Ste95]{ste95}
Joseph Steenbrink.
\newblock Logarithmic embeddings of varieties with normal crossings and mixed
  {Hodge} structures.
\newblock {\em Mathematische Annalen}, 301(1):105--118, 1995.

\bibitem[Sus18]{sus18}
Dmitry Sustretov.
\newblock Gromov-{Hausdorff} limits of flat {Riemannian} surfaces.
\newblock {\em arXiv preprint arXiv:1802.03818}, 2018.

\bibitem[Sus22]{hessian}
Dmitry Sustretov.
\newblock {Hessian} metrics with distribution coefficients on a 2-sphere.
\newblock {\em arXiv preprint arXiv:2212.10640}, 2022.

\bibitem[Thu07]{thuillier07}
Amaury Thuillier.
\newblock G{\'e}om{\'e}trie toro{\"\i}dale et g{\'e}om{\'e}trie analytique non
  archim{\'e}dienne. application au type d’homotopie de certains sch{\'e}mas
  formels.
\newblock {\em manuscripta mathematica}, 123(4):381--451, 2007.

\bibitem[Tos20]{tosatti20}
Valentino Tosatti.
\newblock Collapsing {Calabi-Yau} manifolds.
\newblock {\em arXiv preprint arXiv:2003.00673}, 2020.

\bibitem[Voi03]{voisin}
Claire Voisin.
\newblock {\em Hodge Theory and Complex Algebraic Geometry, Volume 2},
  volume~77.
\newblock Cambridge University Press, 2003.

\bibitem[Yam21]{yamamoto}
Yuto Yamamoto.
\newblock Tropical contractions to integral affine manifolds with
  singularities.
\newblock {\em arXiv preprint arXiv:2105.10141}, 2021.

\bibitem[Yu15]{tony}
Tony~Yue Yu.
\newblock Tropicalization of the moduli space of stable maps.
\newblock {\em Mathematische Zeitschrift}, 281(3):1035--1059, 2015.

\end{thebibliography}

\vspace{3ex}
\noindent {\sc Department of Mathematics\\
  KU Leuven\\
  Celestijnenlaan 200B \\
  B-3001 Leuven (Heverlee)\\
  Belgium\\}
{\tt dsustretov.math@gmail.com\\}
{\tt dmitry.sustretov@kuleuven.be\\}

\end{document}